\date{February 16th, 2007}

\documentclass[12pt]{amsart}
\usepackage{latexsym,amsmath,amsfonts,amscd,amssymb}
\usepackage{graphics}
\newtheorem{theorem}{Theorem}[section]
\newtheorem{lemma}[theorem]{Lemma}
\newtheorem{proposition}[theorem]{Proposition}
\newtheorem{corollary}[theorem]{Corollary}

\newtheorem{definition}[theorem]{Definition}

\newtheorem{counterexample}[theorem]{Counter-example}
\newtheorem{example}[theorem]{Example}

\theoremstyle{remark}
\newtheorem{remark}[theorem]{Remark}

\newcommand{\la}{\langle}
\newcommand{\ra}{\rangle}
\newcommand{\too}{\longrightarrow}

\newcommand{\inc}{\hookrightarrow}
\newcommand{\bd}{\partial}
\newcommand{\x}{\times}
\newcommand{\ox}{\otimes}

\newcommand{\Vol}{\operatorname{Vol}}

\newcommand{\GL}{\operatorname{GL}}

\newcommand{\Int}{\operatorname{Int}}
\newcommand{\Hom}{\operatorname{Hom}}

\newcommand{\im}{\operatorname{im}}

\newcommand{\id}{\operatorname{id}}
\newcommand{\Hol}{\operatorname{Hol}}
\newcommand{\supp}{\operatorname{supp}}

\newcommand{\cA}{{\mathcal A}}

\newcommand{\cC}{{\mathcal C}}

\newcommand{\cM}{{\mathcal M}}

\newcommand{\cL}{{\mathcal L}}

\newcommand{\cS}{{\mathcal S}}
\newcommand{\cT}{{\mathcal T}}

\newcommand{\cV}{{\mathcal V}}

\newcommand{\PP}{{\mathbb P}}
\newcommand{\QQ}{{\mathbb Q}}
\newcommand{\RR}{{\mathbb R}}
\newcommand{\TT}{{\mathbb T}}

\newcommand{\ZZ}{{\mathbb Z}}

\renewcommand{\a}{\alpha}
\renewcommand{\b}{\beta}
\renewcommand{\d}{\delta}
\newcommand{\g}{\gamma}

\renewcommand{\o}{\omega}

\renewcommand{\t}{\tau}

\newcommand{\G}{\Gamma}

\newcommand{\fro}{{\mathfrak{o}}}

\begin{document}

\title{Ergodic solenoidal homology}


\author[V. Mu{\~n}oz]{Vicente Mu{\~n}oz}
\address{Departamento de Matem\'aticas,
CSIC, Serrano 113 bis, 28006 Madrid, Spain}

\address{Facultad de
Matem\'aticas, Universidad Complutense de Madrid, Plaza de Ciencias
3, 28040 Madrid, Spain}

\email{vicente.munoz@imaff.cfmac.csic.es}

\author[R. P{\'e}rez Marco]{Ricardo P{\'e}rez Marco}
\address{CNRS, LAGA UMR 7539, Universit\'e Paris XIII \\ 
99, Avenue J.-B. Cl\'ement, 93430-Villetaneuse, France}

\address{UCLA, Mathematics Dept., Box 951555, Los Angeles, CA
90095-1555, USA}

\email{ricardo@math.univ-paris13.fr}

\thanks{First author supported through MCyT grant MTM2004-07090-C03-01
(Spain) and NSF grant DMS-0111298 (US). 
Second author supported by CNRS (UMR 7539) and NSF grant
DMS-0202494.}

\begin{abstract}
We define generalized currents associated with immersions of
abstract solenoids with a transversal measure. We realize
geometrically the full real homology of a compact manifold with
these generalized currents, and more precisely with immersions of
minimal uniquely ergodic solenoids. This makes precise and geometric
De Rham's realization of the real homology by only using a
restricted geometric subclass of currents. These generalized
currents do extend Ruelle-Sullivan and Schwartzman currents. We
extend Schwartzman theory beyond dimension $1$ and provide a unified
treatment of Ruelle-Sullivan and Schwartzman theories via Birkhoff's
ergodic theorem for the class of immersions of 
controlled solenoids. We develop some intersection theory of these
new generalized currents that explains why the realization theorem
cannot be achieved only with Ruelle-Sullivan currents.
\end{abstract}

\maketitle

\tableofcontents

\newpage

\section{Introduction} \label{sec:introduction}

We consider a smooth compact oriented manifold $M$ of dimension
$n\geq 1$. Any closed oriented submanifold $N\subset M$ of dimension
$0\leq k\leq n$ determines a homology class in $H_k(M, \ZZ)$. This
homology class in $H_k(M,\RR)$, as dual of De Rham cohomology, is
explicitly given by integration of the restriction to $N$ of
differential $k$-forms on $M$. Unfortunately, because of topological
reasons dating back to Thom \cite{Thom}, not all integer homology
classes in $H_k(M,\ZZ )$ can be realized in such a way by a compact
submanifold. Geometrically, we can realize any class in $H_k(M,
\ZZ)$ by topological $k$-chains. The real homology $H_k(M,\RR)$
classes are only realized by formal combinations with real
coefficients of $k$-cells. This is not satisfactory for various
reasons. In particular, for diverse purposes it is important to have
an explicit realization, as geometric as possible, of real homology
classes.

\medskip

The first contribution in this direction came in 1957 from the work
of S. Schwartzman \cite{Sc}. Schwartzman showed how, by a limiting
procedure, one-dimensional curves embedded in $M$ can define a real
homology class in $H_1(M,\RR)$. More precisely, he proved that this
happens for almost all curves solutions to a differential equation
admitting an invariant ergodic probability measure. Schwartzman's
idea is very natural. It consists on integrating $1$-forms over
large pieces of the parametrized curve and normalizing this integral
by the length of the parametrization. Under suitable conditions, the
limit exists and defines an element  of the dual of $H^1(M,\RR )$,
i.e. an element of $H_1(M, \RR )$. This procedure is equivalent to
the more geometric one of closing large pieces of the curve by
relatively short closing paths. The closed curve obtained defines an
integer homology class. The normalization by the length of the
parameter range provides a class in $H_k(M,\RR )$. Under suitable
hypothesis, there exists a unique limit in real homology when the
pieces exhaust the parametrized curve, and this limit is independent
of the closing procedure. In sections \ref{sec:calibrating},
\ref{sec:Schwartzman-cycles} and \ref{sec:k-schwartzman}, we study
the different aspects of the Schwartzman procedure, that we extend
to higher dimension.

\medskip

Later in 1975, D. Ruelle and D. Sullivan \cite{RS} defined, for
arbitrary dimension $0\leq k\leq n$, geometric currents by using
oriented $k$-laminations embedded in $M$ and endowed with a
transversal measure. They applied their results to Axiom A
diffeomorphisms. In a later article Sullivan \cite{Su} extended
further these results and their applications. The point of view of
Ruelle and Sullivan is also based on duality. The observation is
that $k$-forms can be integrated on each leaf of the lamination and
then all over the lamination using the transversal measure. This
makes sense locally in each flow-box, and then it can be extended
globally by using a partition of unity.  The result only depends on
the cohomology class of the $k$-form. In section
\ref{sec:Ruelle-Sullivan} we review and extend Ruelle-Sullivan
theory.

\medskip

It is natural to ask wether it is possible to realize every real
homology class using a Ruelle-Sullivan current. A first result, with
a precedent in \cite{HM}, confirms that this is not the case (see
section \ref{sec:Ruelle-Sullivan}).

\begin{theorem} \label{thm:11}
Homology classes with non-zero self-intersection cannot be
represented by Ruelle-Sullivan currents.
\end{theorem}

More precisely, for each Ruelle-Sullivan lamination with a
non-atomic transversal measure, we can construct a smooth
$(n-k)$-form which provides the dual in De Rham cohomology. Using
it, we prove that the self-intersection of a Ruelle-Sullivan current
is $0$, therefore it is not possible to represent a real homology
class in $H_k(M,\RR )$ with non-zero self-intersection. This
obstruction only exists when $n-k$ is even. This may be the
historical reason behind the lack of results on the representation
of an arbitrary homology class by Ruelle-Sullivan currents.

\medskip

Therefore, in order to hope to  represent every real homology class
we must first enlarge the class of Ruelle-Sullivan currents. This is
done by considering immersions of abstract oriented solenoids. We
define a $k$-solenoid to be a Hausdorff compact space foliated by
$k$-dimensional leaves with finite dimensional transversal structure
(see the precise definition in section \ref{sec:minimal}). For these
oriented solenoids we can consider $k$-forms that we can integrate
provided that we are given a transversal measure invariant by the
holonomy group. We define an immersion of a solenoid $S$ into $M$ to
be a regular map $f: S\to M$ that is an immersion in each leaf. If
the solenoid $S$ is endowed with a transversal measure $\mu$, then
any smooth $k$-form in $M$ can be pulled back to $S$ by $f$ and
integrated. The resulting numerical value only depends on the
cohomology class of the $k$-form. Therefore we have defined a
current that we denote by $[f,S_\mu ] \in H_k(M,\RR )$ and that we
call a generalized current. Using these generalized currents, the
above mentioned obstruction disappears. Our main result is:

\begin{theorem}\textbf{\em (Realization Theorem)} \label{thm:1.2}
Every real homology class in $H_k(M,\RR )$ can be realized by a
generalized current $[f,S_\mu]$ where $S_\mu$ is an oriented,
minimal, uniquely ergodic solenoid.
\end{theorem}

Minimal and uniquely ergodic solenoids are defined later on. This
result gives a geometric version and makes precise De Rham's
realization theorem of homology classes by abstract currents, i.e.
forms with coefficients distributions.

\medskip

But the space of solenoids is large, and we would like to realize
the real homology classes by a minimal class of solenoids enjoying
good properties. We are first naturally led to topological
minimality. As we prove in section \ref{sec:minimal}, the spaces of
$k$-solenoids is inductive and therefore there are always minimal
$k$-solenoids. However, the transversal structure and the holonomy
group of minimal solenoids can have a rich structure. In particular,
such a solenoid may have many different transversal measures, each
one yielding a different generalized current for the same immersion
$f$. Also when we push Schwartzman ideas beyond $1$-homology for
some nice classes of solenoids, we see that in general, even when
the immersion is an embedding, the generalized current does not
necessarily coincide with the Schwartzman homology class of the
immersion of each leaf (actually not even this Schwartzman class
needs to be well defined). Indeed the classical literature lacks of
information about the precise relation between Ruelle-Sullivan and
Schwartzman currents. One would naturally expect that there is some
relation between the generalized currents and the Schwartzman
current (if defined) of the leaves of the lamination. We study this
problem in section \ref{sec:1-schwartzman} for $1$-dimensional
currents and in section \ref{sec:k-schwartzman} in general. The main
result is that there is such relation for the class of minimal,
ergodic and controlled solenoids (see definition in section
\ref{sec:k-schwartzman}) for which the transversal structure is well
behaved. A controlled solenoid has a trapping region (see definition
\ref{def:trapping}), and the holonomy group is generated by a single
map. Then the bridge between generalized currents and Schwartzman
currents of the leaves is provided by Birkhoff's ergodic theorem as
explained in sections \ref{sec:Schwartzman-cycles} and
\ref{sec:k-schwartzman}.

\begin{theorem}\label{thm:1.3}
Let $S_\mu$ be a controlled and minimal solenoid endowed with an
ergodic transversal measure $\mu$. Let $f: S_\mu \to M$ be an
immersion of $S_\mu$ into $M$. Then for $\mu$-almost all leaves
$l\subset S_\mu$, the Schwartzman homology class of $f(l)\subset M$
is well defined and coincides with the generalized current
$[f,S_\mu]$.
\end{theorem}

We are particularly interested in uniquely ergodic solenoids, with
only one ergodic transversal measure. As is well known, in this
situation we have uniform convergence of Birkhoff's sums, which
implies the stronger result:

\begin{theorem}\label{thm:1.4}
Let $S_\mu$ be a minimal and uniquely ergodic solenoid. Let $f:
S_\mu \to M$ be an immersion of $S_\mu$ into $M$. Then for all
leaves $l\subset S_\mu$, the Schwartzman homology class of
$f(l)\subset M$ is well defined and coincides with the generalized
current $[f,S_\mu]$.
\end{theorem}

We also make a thorough study of Riemannian solenoids. We identify
transversal measures with the class of measures that desintegrate as
volume along leaves (daval measures), we define Schwartzman measures
as limits of $k$-volume normalized measures along leaves, and we
prove a canonical decomposition of measures into a daval measure and
a singular part, corresponding to the classical Lebesgue
decomposition on manifolds (see section \ref{sec:riemannian}).

\medskip

Numerous arguments and auxiliary results in this article are part of
the folklore or have appeared earlier in the literature on
foliations or laminated spaces. We do not give systematic
references. The reader can consult the monograph \cite{Candel} for a
comprehensive bibliography and classical results.

\bigskip

\noindent \textbf{Acknowledgements.} \ 
The authors are grateful to Alberto Candel, Etienne Ghys and Jaume Amor\'os
for their comments and interest on this work. In particular, Etienne
Ghys early pointed out on the impossibility of realization in
general of integer homology classes by embedded manifolds.

The first author wishes to acknowledge Universidad Complutense de
Madrid and Institute for Advanced Study at Princeton for their
hospitality and for providing excellent working conditions.  The second author
thanks Jean Bourgain and the IAS at Princeton for their hospitality
and facilitating the collaboration of both authors.

\newpage

\section{Minimal solenoids} \label{sec:minimal}

We first define abstract solenoids, which are the main tools in
this article.

\begin{definition}\label{def:W}
  Let $0\leq r,s\leq \omega$, and let $k,l\geq 0$ be two integers.
  A foliated manifold (of dimension $k+l$, with $k$-dimensional
  leaves, of regularity $C^{r,s}$) is a smooth
  manifold $W$ of dimension $k+l$ endowed with an atlas
  $\cA=\{ (U_i,\varphi_i)\}$ whose changes of
  charts
 $$
 \varphi_{ij}=\varphi_i\circ \varphi_j^{-1}: \varphi_j(U_i\cap U_j)
 \to \varphi_i(U_i\cap U_j) \ \, ,
 $$
are of the form $\varphi_{ij}(x,y)=(X_{ij}(x,y), Y_{ij}(y))$,
where $Y_{ij}(y)$ is of class $C^s$ and $X_{ij}(x,y)$ is of class
$C^{r,s}$.

A flow-box for $W$ is a pair $(U,\varphi)$ consisting of an open
subset $U\subset W$ and a map $\varphi:U\to \RR^k\x \RR^l$ such that
$\cA \cup \{(U,\varphi)\}$ is still an atlas for $W$.
\end{definition}

Clearly an open subset of a foliated manifold is also a foliated
manifold.

Given two foliated manifolds $W_1$, $W_2$ of dimension $k+l$, with
$k$-dimensional leaves, and of regularity $C^{r,s}$, a regular map
$f:W_1\to W_2$ is a continuous map which is locally, in
flow-boxes, of the form $f(x,y)=(X(x,y),Y(y))$, where $Y$ is of
class $C^{s}$ and $X$ is of class $C^{r,s}$.

A diffeomorphism $\phi:W_1\to W_2$ is a homeomorphism such that
$\phi$ and $\phi^{-1}$ are both regular maps.

\begin{definition}\label{def:k-solenoid}\textbf{\em ($k$-solenoid)}
Let  $0\leq r,s\leq \omega$, and let $k,l\geq 0$ be two integers.
A pre-solenoid of dimension $k$, of class $C^{r,s}$ and
transversal dimension $l$ is a pair $(S,W)$ where $W$ is a
foliated manifold and $S\subset W$ is a compact subspace which is
a collection of leaves.

Two pre-solenoids $(S,W_1)$ and $(S,W_2)$ are equivalent if there
are open subsets $U_1\subset W_1$, $U_2\subset W_2$ with $S\subset
U_1$ and $S\subset U_2$, and a diffeomorphism $f: U_1\to U_2$ such
that $f$ is the identity on $S$.

A $k$-solenoid of class $C^{r,s}$ and transversal dimension $l$
(or just a $k$-solenoid, or a solenoid) is an equivalence class of

pre-solenoids.
\end{definition}

We usually denote a solenoid by $S$, without making explicit
mention of $W$. We shall say that $W$ defines the solenoid
structure of $S$.

\begin{definition}\label{def:flow-box}
\textbf{\em (Flow-box)} Let $S$ be a solenoid. A flow-box for $S$
is a pair $(U,\varphi)$ formed by an open subset $U\subset S$ and
a homeomorphism
 $$
 \varphi : U\to D^k\times K(U) \, ,
 $$
where $D^k$ is the $k$-dimensional open ball and $K(U)\subset
\RR^l$, such that there exists a foliated manifold $W$ defining
the solenoid structure of $S$, $S\subset W$, and a flow-box
$\hat\varphi:\hat U\to \RR^k\x\RR^l$ for $W$, with $U=\hat U\cap
S$,  $\hat\varphi(U)= D^k\x K(U)\subset \RR^k\x\RR^l$ and
$\varphi=\hat\varphi_{|U}$.

The set $K(U)$ is the transversal space of the flow-box. The
dimension $l$ is the transversal dimension.
\end{definition}

As $S$ is locally compact, any point of $S$ is contained in a flow-box $U$ whose closure
$\overline U$ is contained in a bigger flow-box. For such flow-box, $\overline U
\cong \overline{D}^k \times \overline K(U)$, where $\overline D^k$
is the closed unit ball, $\overline K(U)$ is some compact subspace of
$\RR^l$, and $U=D^k \times  K(U) \subset
\overline{D}^k \times \overline K(U)$. We might call these flow-boxes \emph{good}.
All flow-boxes that we shall use are of this type so we shall not
append any appelative to them.

We refer to $k$ as the dimension of the solenoid and we write
 $$
 k=\dim S \, .
 $$
Note that, contrary to manifolds, this dimension in general does not
coincide with the topological dimension of $S$. The local structure
and compactness imply that solenoids are metrizable topological
spaces. The Hausdorff dimension of the transversals $K(U)$ is well
defined and the Hausdorff dimension of $S$ is well defined, and
equal to
 $$
 \dim_H S=k+\max_i \dim_H K(U_i) \leq k+l<+\infty \, .
 $$

When the transversals of flow-boxes $K(U)\subset \RR^l$ are open
sets of $\RR^l$ we talk about full transversals. In this case the
solenoid structure is a $(k+l)$-dimensional compact manifold
foliated by $k$-dimensional leaves.

\begin{remark}
The definition of solenoid admits various generalizations. We
could focus on intrinsic changes of charts in $S$ with some
transverse Whitney regularity but without requiring a local
diffeomorphism extension. Such a definition would be more general,
but it is not necessary for our purposes. The present definition
balances generality and simplicity.

Another alternative generalization would be to avoid any
restrictive transversal assumption beyond continuity, and allow
for transversals of flow-boxes any topological space $K(U)$. But a
fruitful point of view is to regard the theory of solenoids as a
generalization of the classical theory of manifolds. Therefore it
is natural to restrict the definition only allowing finite
dimensional transversal spaces. For an alternative approach see
\cite{MoSch}.
\end{remark}

\medskip

\begin{definition}\label{def:diffeomorphisms}
\textbf{\em (Diffeomorphisms of solenoids)} Let $S_1$ and $S_2$ be
two  $k$-solenoids of class $C^{r,s}$ with the same transversal
dimension. A $C^{r,s}$-diffeomorphism $f:S_1\to S_2$ is the
restriction of a $C^{r,s}$-diffeomorphism $\hat{f}:W_1\to W_2$ of
two foliated manifolds defining the solenoid structures of $S_1$ and
$S_2$, respectively.

A homeomorphism of solenoids is a diffeomorphism of class $C^{0,0}$.
\end{definition}

This defines the category of smooth solenoids of a given regularity.
Note that we have the subcategory of smooth solenoids with full
transversals, and we have a forgetting functor into the category of
smooth manifolds.

\begin{definition} \label{def:leaf}\textbf{\em (Leaf)}
A leaf of a $k$-solenoid $S$ is a leaf $l$ of any foliated manifold
$W$ inducing the solenoid structure of $S$, such that $l\subset S$.
Note that this notion is independent of $W$.
%
\end{definition}

Note that $S\subset W$ is the union of a collection of leaves.
Therefore, for a leaf $l$ of $W$ either $l\subset S$ or $l\cap
S=\emptyset$.

\medskip

Observe that when the transversals of flow-boxes $K(U)$ are totally
disconnected then the leaf-equivalence coincides with path
connection equivalence, and the leaves are the path connected
components of $S$.

\begin{definition} \label{def:oriented-solenoid}\textbf{\em (Oriented solenoid)}
An oriented solenoid is a solenoid $S$ such that there is a
foliated manifold $W\supset S$ inducing the solenoid structure of
$S$, where $W$ has oriented leaves.
\end{definition}

It is easy to see that $S$ is oriented if and only if there is an
orientation for the leaves of $S$ such that there is a covering by
flow-boxes which preserve the orientation of the leaves.

Notice that we do not require $W$ oriented. For example, we can
foliate a Mo\"ebius strip and create an oriented solenoid.

\begin{definition}\label{def:space-solenoids}
We define $\cS^{r,s}_{k,l}$ to be the space of $C^{r,s}$
$k$-solenoids with transversal dimension $l$.
\end{definition}

\begin{proposition}\label{prop:sub-solenoids}
Let $S_0\in \cS_{k,l}^{r,s}$ be a solenoid. A non-empty compact
subset $S$ of $S_0$ which is a union of leaves is a $k$-solenoid
of class $C^{r,s}$ and transversal dimension $l$.
\end{proposition}

\begin{proof}
Let $W$ be a $C^{r,s}$-foliated manifold defining the solenoid
structure of $S_0$. Then $S\subset W$ and $W$ defines a
$C^{r,s}$-solenoid structure for $S$.

Note that the flow-boxes of $S_0$ give, by restriction to $S$,
flow-boxes for $S$.
\end{proof}

\begin{corollary} \label{cor:connected-component}
Connected components of solenoids $\cS_{k,l}^{r,s}$ are in
$\cS_{k,l}^{r,s}$.
\end{corollary}

\begin{theorem}\label{thm:inductive}
The space $(\cS_{k,l}^{r,s}, \subset )$ ordered by inclusion is an
inductive set.
\end{theorem}

\begin{proof}
Let $(S_n) \subset \cS_{k,l}^{r,s}$ be a nested sequence of
solenoids, $S_{n+1} \subset S_n$. Define
 $$
 S_\infty = \bigcap_n S_n \, .
 $$
Then $S_\infty$ is a non-empty compact subset of $S_1$  as
intersection of a nested sequence of such sets. It is also a union
of leaves since each $S_n$ is so. Therefore by
proposition \ref{prop:sub-solenoids}, it is an element of $\cS_{k,l}^{r,s}$.
\end{proof}

\begin{corollary}\label{cor:minimal-elements}
The space $\cS_{k,l}^{r,s}$ has minimal elements.
\end{corollary}

\begin{proposition}\label{prop:dense}
If $S\in \cS_{k,l}^{r,s} $ is minimal then $S$ is connected and
all leafs of $S$ are dense.
\end{proposition}

\begin{proof}
Each connected component of $S$ is a solenoid, thus by minimality
$S$ must be connected.

Also the closure ${\overline L}$ of any leaf $L\subset S$ is a
non-empty compact set union of leaves. Thus it is a solenoid and by
minimality we must have ${\overline L} =S$.
\end{proof}

\section{Topological transversal structure of solenoids}\label{sec:topological}

\begin{definition}\label{def:transversal}\textbf{\em (Transversal)}
Let $S$ be a $k$-solenoid. A local transversal at a point $p\in S$
is a subset $T$ of $S$ with $p\in T$, such that there is a flow-box $(U, \varphi)$
of $S$ with $U$ a neighborhood of $p$ containing $T$ and such that
 $$
 \varphi (T)=\{0\}\x K(U) \, .
 $$

A transversal $T$ of $S$ is a compact subset of $S$ such that for
each $p\in T$ there is an open neighborhood $V$ of $p$ such that
$V\cap T$ is a local transversal at $p$.
\end{definition}

If $S$ is a $k$-solenoid of class $C^{r,s}$, then any transversal
$T$ inherits an  $l$-dimensional $C^s$-Whitney structure.

We clearly have:

\begin{proposition}\label{prop:two-transversals}
The union of two disjoint transversals is a transversal.
\end{proposition}

\begin{definition}\label{def:global-transversal}
A transversal $T$ of $S$ is a global transversal if all leaves
intersect $T$.
\end{definition}

The next proposition is clear.

\begin{proposition}\label{prop:global-transversals}
The union of two disjoint transversals, one of them global, is a
global transversal.
\end{proposition}

\begin{proposition}\label{prop:global-transversal}
If $S$ is minimal then all transversals are global. Moreover, if $S$
is minimal then any local transversal intersects all leaves of $S$.
\end{proposition}

\begin{proof}
It is enough to see the second statement, since it implies the
first. Let $U$ be a flow-box and $T=\varphi^{-1}(\{0\}\x K(U))$ a
local transversal. By proposition \ref{prop:dense}, all leaves
intersect $U$ and therefore they intersect $T$.
%
\end{proof}

Observe that the definition of solenoid with regular transverse
structure says that $S$ is always embedded in a $(k+l)$-dimensional
manifold $W$. Therefore $S\subset W$ has an interior and a boundary
relative to $W$. These sets do not depend on the choice of $W$.

\begin{definition}\label{def:interior}\textbf{\em (Proper interior and boundary)}
Let $S$ be a $k$-solenoid. Let $W$ be a foliated manifold defining
the solenoid structure of $S$. The proper interior of $S$ is the
interior of $S$ as a subset of $W$, considered as a
$(k+l)$-dimensional manifold.

The proper boundary of $S$ is defined as the complement in $S$ of
the proper interior.
\end{definition}

Let $\hat\varphi:\hat{U} \to \RR^k\x \RR^l$ be a flow-box for $W$
such that $U=\hat{U}\cap S$ and $\varphi=\hat\varphi_{|U}: U\to D^k \x
K(U)$ is a flow-box for $S$. Then $K(U)\subset \RR^l$. The proper
interior, resp.\ the proper boundary, of $S$, intersected with $U$,
consists of the collection of leaves $\varphi^{-1}(D^k \x \{ y\})$, where
$y\in K(U)$ is in the interior, resp.\ boundary, of $K(U)\subset
\RR^l$.

Note that the proper boundary of a solenoid that is a foliation of a
manifold is empty. We prove the converse.

\begin{proposition}\label{prop:sub-solenoid}
If the proper boundary of $S$ is non-empty then it is a
sub-solenoid of $S$.
\end{proposition}

\begin{proof}
The proper boundary is a compact subset of $S$ and a union of
leaves. The result follows from proposition
\ref{prop:sub-solenoids}.
\end{proof}

\begin{theorem}\label{thm:empty-interior}
Let $S\in \cS$ be a minimal solenoid. If $S$ is not the foliation of
a manifold then $S$ has empty proper interior, i.e. $K(U)\subset
\RR^l$ has empty interior for any flow-box $(U,\varphi )$.
\end{theorem}

\begin{proof}
The proper boundary of $S$ is non-empty  because otherwise, for each
flow-box $U$, $K(U)\subset \RR^l$ is an open set. Thus $S$ would be
an open subset of $W$, where $W$ is a foliated manifold defining the
solenoid structure of $S$, and so $S$ is itself a foliated
$(k+l)$-manifold. This contradicts the assumption.

Now by minimality the proper boundary must coincide with $S$ and the
proper interior is empty.
\end{proof}

Solenoids with a one dimensional transversal play a prominent role
in the sequel. We have for these the following structure theorem.

\begin{theorem} \label{thm:1-transversal}
\textbf{\em (Minimal solenoids with a $1$-dimensional
transversal).} Let $S\in \cS$ be a minimal $k$-solenoid which
admits a $1$-dimensional transversal $T$.

Then we have two cases:
 \begin{enumerate}
 \item {}$T$ is a finite union of circles, and $S$ is a $1$-dimensional foliation of a
 connected manifold of dimension $k+1$.
 \item {}$T$ is totally disconnected, in which case we have two further possibilities:
  \begin{enumerate}
   \item {}$T$ is a finite set and $S$ is a connected manifold of dimension $k$,
   \item {}$T$ is a Cantor set.
  \end{enumerate}
 \end{enumerate}
\end{theorem}

\begin{proof}
We define the proper interior of $T$ as the intersection of the
proper interior of $S$ with $T$. Now we have two cases.

If the proper interior of $T$ is non-empty, then the proper interior
of $S$ is non-empty. Then the complement of the proper interior of
$S$, if non-empty, is a sub-solenoid of $S$, contradicting
minimality. Thus the proper interior of $S$ is all $S$, so the
proper interior of $T$ is the whole of $T$. This means that any
point $p\in T$ has a neighbourhood (in $T$) homeomorphic to an
interval. Therefore $T$ is a topological compact $1$-dimensional
manifold, thus a finite union of circles. This ends the first case.

If the proper interior of $T$ is empty, then $T$ is totally
disconnected. In this case, if $T$ has an isolated point $p$, then
$S$ has only one leaf because by minimality any other leaf must
accumulate the leaf containing $p$, and this is only possible if
it coincides with it. Then $S$ is a $k$-dimensional connected
manifold. If $T$ has no isolated points, then $T$ is non-empty,
perfect, compact and totally disconnected, i.e. it is a Cantor
set.
\end{proof}

\section{Holonomy, Poincar\'e return and suspension}\label{sec:holonomy}

\begin{definition}\label{def:holonomy}\textbf{\em (Holonomy)}
Given two points $p_1$ and $p_2$ in the same leaf, two local
transversals $T_1$ and $T_2$, at $p_1$ and $p_2$ respectively, and
a path $\gamma :[0,1]\to S$, contained in the leaf with endpoints
$\g(0)=p_1$ and $\g(1)=p_2$, we define a germ of a map (the
holonomy map)
 $$
 h_\gamma: (T_1,p_1) \to (T_2, p_2)\, ,
 $$
by lifting $\gamma$ to nearby leaves.

We denote by ${\Hol}_S (T_1, T_2)$ the set of germs of holonomy
maps from $T_1$ to $T_2$.
\end{definition}

The following result is clear.

\begin{proposition} \label{prop:holonomies} If $T_1$ and $T_2$ are global transversals then
the sets of holonomy maps from $T_1$ to $T_2$ is non-empty. In
particular, if $S$ is minimal the set of holonomy maps between two
arbitrary local transversals is non-empty.
\end{proposition}

\begin{definition}\label{def:holonomy-group} \textbf{\em (Holonomy pseudo-group)}
The holonomy pseudo-group of a local transversal $T$ is the
pseudo-group of holonomy maps from $T$ into itself. We denote it
by ${\Hol}_S(T)={\Hol}_S(T, T)$.

The holonomy pseudo-group of $S$ is the pseudo-group of all
holonomy maps. We denote it by ${\Hol}_S$,
 $$
 {\Hol}_S=\bigcup_{T_1 ,T_2} {\Hol}_S(T_1, T_2) \, .
 $$
\end{definition}

The following proposition is obvious.

\begin{proposition} \label{thm:orbit} The orbit of a point
$x\in S$ by the holonomy pseudo-group coincides with the leaf
containing $x$.

Therefore, a solenoid $S$ is minimal if and only if the action of
the holonomy pseudo-group is minimal, i.e. all orbits are dense.
\end{proposition}

The Poincar\'e return map construction reduces sometimes the
holonomy to a classical dynamical system.

\begin{definition}\label{def:Poincare}\textbf{\em (Poincar\'e return map)}
Let $S$ be an oriented minimal $1$-solenoid and $T$ be a local
transversal. Then the holonomy return map is well defined for all
points in $T$ and defines the Poincar\'e return map
 $$
 R_T:T\to T \, .
 $$
\end{definition}

The return map is well defined because in minimal solenoids ``half''
leaves are dense.

\begin{lemma} \label{lem:half-leaves}
Let $S$ be a minimal $1$-solenoid. Let $p_0\in S$ and let
$l_0\subset S$ be the leaf containing $p_0$. The point $p_0$ divides
the leaf $l_0$ into two connected components. They are both dense in
$S$.
\end{lemma}

\begin{proof}
Consider one connected component of $l_0-\{p_0\}$, and let $C$ be
its accumulation set. Then $C$ is non-empty, by compactness of $S$,
and it is compact, as a closed subset of the compact solenoid $S$.
It is also a union of leaves because if $l\subset S$ is a leaf, then
$C\cap l$ is open in $l$ as is seen in flow-boxes, and also $C\cap
l$ is closed in $l$. Therefore by connectedness of $l$, $C\cap l$ is
empty or $l\subset C$.

We conclude that $C$ is a sub-solenoid, and by minimality we have
$C=S$.
\end{proof}

When $S$ admits a global transversal (in particular when $S$ is
minimal and admits a transversal) and the Poincar\'e return map is
well defined, we have
that it is continuous (without any assumption on minimality of $S$).

\begin{proposition}\label{prop:continuity}
Let $S$ be an oriented $1$-solenoid and let $T$ be a global transversal
such that the Poincar\'e return map $R_T$ is well defined. Then the
holonomy return map $R_T$ is continuous.
If the Poincar\'e return map for the reversed orientation of $S$ is
also well defined, then $R_T$ is a homeomorphism of $T$. Moreover, if $S$ is a
solenoid of class $C^{r,s}$ then $R_T$ is a $C^s$-diffeomorphism.
\end{proposition}

\begin{proof}
The map $R_T$ is continuous because
the inverse image of an open set is clearly open.

If the Poincar{\'e} return map $R_T^-$ for the same transversal obtained for
the reverse orientation of $S$ is also well defined, then
$R_T$ is bijective because by construction its inverse is $R_T^-$.
Hence $R_T$ is a
homeomorphism of $T$. Moreover, letting $W$ be a foliated manifold defining the solenoid structure of $S$,
$T$ is a subset of an open manifold $U$ of dimension $l$, and the
map $R_T$ extends as a homeomorphism $U_1\to U_2$, where $U_1$, $U_2$ are neighborhoods of $T$
(at least locally). If the
transversal regularity of $S$ is $C^s$ then the local extension of
$R_T$ is a $C^s$-diffeomorphism.
\end{proof}

\medskip

When $T$ is only a local transversal then in general $R_T$ is not
continuous. Nevertheless the discontinuities of $R_T$ are well
controlled in practice and are innocuous when we deal with measure
theoretic properties of $R_T$, as for example in section
\ref{sec:Schwartzman-cycles}.

\medskip

The suspension construction reverses Poincar\'e construction of the
first return map.

\begin{definition}\label{def:suspension}\textbf{\em (Suspension construction)}
Let $X\subset \RR^l$ be a compact set and let $f:X\to X$ be a
homeomorphism which has a $C^s$-diffeomorphism extension to a
neighborhood of $X$ in $\RR^l$. The suspension of $f$ is the
oriented $1$-solenoid $S_f$ defined by the suspension construction
 $$
 S_f=([0,1]\times X)_{/(0,x)\sim (1,f(x))} \, .
 $$
The solenoid $S_f$ has regularity $C^{\omega ,s}$.

The transversal $T=\{0 \}\times X$ is a global transversal and the
associated Poincar\'e return map $R_T:T\to T$ is well defined and
equal to $f$.
\end{definition}

In particular, the theory of dynamical systems for $X\subset \RR^l$
and diffeomorphisms $f:X\to X$ (extending to a neighborhood of $X$)
is contained in the theory of transversal structures of solenoids.

\section{Measurable transversal structure of solenoids}\label{sec:transversal-structure}

\begin{definition} \label{def:transversal-measure} \textbf{\em (Transversal measure)}
Let $S$ be a $k$-solenoid. A transversal measure $\mu=(\mu_T)$ for
$S$ associates to any local transversal $T$ a locally finite measure
$\mu_T$ supported on $T$, which are invariant by the holonomy
pseudogroup. More precisely, if $T_1$ and $T_2$ are two transversals
and $h : V\subset T_1 \to T_2$ is a holonomy map, then
 $$
 h_* (\mu_{T_1}|_{V})= \mu_{T_2}|_{h(V)} \, .
 $$
We assume that a transversal measure $\mu$ is non-trivial, i.e. for
some $T$, $\mu_T$ is non-zero.

We denote by $S_\mu$ a $k$-solenoid $S$ endowed with a transversal
measure $\mu=(\mu_T)$. We refer to $S_\mu$ as a measured solenoid.
\end{definition}

Observe that for any transversal measure $\mu=(\mu_T)$ the scalar
multiple $c\, \mu=(c \, \mu_T)$, where $c>0$, is also a transversal
measure. Notice that there is no natural scalar normalization of
transversal measures.

\begin{definition}\label{def:support} \textbf{\em (Support of a transversal measure)}
Let $\mu=(\mu_T)$ be a transversal measure. We define the support

of $\mu$ by
 $$
 \supp \mu=\bigcup_T \supp \mu_T \, ,
 $$
where the union runs over all local transversals $T$ of $S$.
\end{definition}

\begin{proposition}\label{prop:support}
The support of a transversal measure $\mu$ is a sub-solenoid of $S$.
\end{proposition}

\begin{proof}
For any flow-box $U$, $\supp \mu\cap U$ is closed in $U$, since
$\supp \mu_{K(U)}$ is closed in $K(U)$. Hence, $\supp \mu$ is closed
in $S$. Also, locally in
flow-boxes $\supp \mu$ contains full leaves of $U$. Therefore a
leaf of $S$ is either disjoint from $\supp \mu$ or contained in
$\supp  \mu$. Also $\supp \mu $ is non-empty because $\mu$ is
non-trivial. We conclude that $\supp \mu$ is a sub-solenoid.
\end{proof}

\begin{definition} \label{def:transversal-ergodicity}\textbf{\em (Transverse ergodicity)}
A transversal measure $\mu=(\mu_T )$ on a solenoid $S$ is ergodic if for any Borel set
$A\subset T$ invariant by the pseudo-group of holonomy maps on $T$,
we have
 $$
 \mu_T(A) = 0 \ \ {\hbox{\rm{ or }}} \ \ \mu_T(A) = \mu_T(T) \, .
 $$

We say that $S_\mu$ is an ergodic solenoid.
\end{definition}

\begin{definition} \label{def:transversal-unique-ergodicity}
\textbf{\em (Transverse unique ergodicity)} Let $S$ be a
$k$-solenoid. The solenoid $S$ is transversally uniquely ergodic, or
a uniquely ergodic solenoid, if $S$ has a unique (up to scalars)
transversal measure $\mu$ and moreover $\supp \mu =S$.
\end{definition}

Observe that in order  to define these notions we only need
continuous transversals. These ergodic notions are intrinsic and
purely topological, i.e. if $S_{1}$ and $S_{2}$ are two homeomorphic solenoids by a
homeomorphism $h:S_1\to S_2$, then $S_1$ is uniquely ergodic if and only if $S_2$ is.
If $S_{1,\mu_1}$ and $S_{2,\mu_2}$ are homeomorphic and $\mu_2=h_* \mu_1$ via the
homeomorphism $h:S_1\to S_2$, then
$S_{1,\mu_1}$ is ergodic if and only if $S_{2,\mu_2}$ is.

\medskip

These notions of ergodicity generalize the classical ones and do
exactly correspond to the classical notions in the situation described
by the next theorem.

\begin{theorem}\label{thm:ergodic}
Let $S$ be an oriented $1$-solenoid. Let $T$ be a global transversal
such that the Poincar\'e return map $R_T: T\to T$ is well defined.

Then the solenoid $S$ is minimal, resp.\ ergodic, uniquely ergodic,
if and only if $R_T$ is minimal, resp.\ ergodic, uniquely ergodic.
\end{theorem}

\begin{proof}
We have by proposition \ref{prop:continuity} that $R_T$ is
continuous. A leaf of $S$ is dense if and only if its intersection
with $T$ is a dense orbit of $R_T$, hence the equivalence of
minimality.

For the ergodicity, observe that we have a correspondence between
measures on $T$ invariant by $R_T$ and transversal measures for $S$.
Each transversal measure for $S$, locally defines a measure on $T$,
hence defines a measure on $T$. Conversely, given a measure $\nu$ on
$T$, we can transport $\nu$ in order to define a measure in each
local transversal $T'$ in the following way. We can define a map
$R_{T',T}: T'\to T$ of first impact on $T$ by following leaves of
$S$ from $T'$ in the positive orientation. By the global character
of the transversal this map is well defined. By construction
$R_{T',T}$ is injective. So we can define
$\mu_{T'}=R_{T',T}^*\nu_{R_{T',T}(T')}$. Then $(\mu_{T'})$ defines a
transversal measure. The equivalence of unique ergodicity follows.
Also $\nu$ is ergodic if and only $(\mu_{T'})$ is ergodic because any
decomposition of $\nu=\nu_1+\nu_2$ induces a decomposition of
$(\mu_{T'})$ by the transversal measures corresponding to the
decomposing measures.
\end{proof}

When we have an ergodic oriented $1$-solenoid $S_\mu$ and $T$ is a
local transversal, then the Poincar\'e return map is well defined
$\mu_T$-almost everywhere and $\mu_T$ is ergodic.

\begin{proposition}\label{prop:5.7}
Let $S$ be an oriented $1$-solenoid and let $T$ be a local
transversal of $S$. Let $\mu$ be an ergodic transversal measure for
$S$. Then the Poincar\'e return map $R_T$ is well defined for
$\mu_T$-almost all points of $T$ and $\mu_T$ is an ergodic measure
of $R_T$.
\end{proposition}

\begin{proof}
Let $A_T\subset T$ be the set of wandering points of $T$, i.e. those
points whose positive half leaves through them never meet $T$ again.
Clearly $A_T$ is a Borel set. If $\mu_T (A_T)\not= 0$ we can
decompose $\mu_T$ by decomposing $\mu_{T |A_T}$ and transporting the
decomposition (back and forward) by the holonomy in order to
decompose the transversal measure. Therefore $\mu_T (A_T)=0$. As
before, a decomposition of $\mu_T$ into invariant measures by $R_T$
yields a decomposition of the transversal measure $\mu$ invariant by
holonomy.
\end{proof}

Recall that a dynamical system is minimal when all orbits are dense,
and that uniquely ergodic dynamical systems are minimal. We have the
same result for uniquely ergodic solenoids.

\begin{proposition}\label{cor:ergodic-1-solenoid}
An oriented uniquely ergodic $1$-solenoid $S$ is minimal.
\end{proposition}

\begin{proof}
If $S$ has a non-dense leaf $l\subset S$, we can consider a local
transversal $T_0$ such that $T_0\cap \bar l \not= \emptyset , T_0$.
Let $(l_n)$ be an exhaustion of $l$ by compact subsets. Let $\mu_n$
be the atomic probability measure on $T_0$ equidistributed on the
intersection of $l_n$ with $T_0$. Any limit measure $\mu_{n_k}\to
\nu$ is a probability measure on $T_0$ with $\supp \nu \subset
T_0\cap \bar l$. It follows easily that $\nu$
is invariant by the holonomy on $T_0$.
Transporting by the holonomy, $\nu$ defines a
transversal measure $\mu=(\mu_T)$ (up to normalization, in each transversal it is
also a limit of counting measures). But this contradicts unique
ergodicity since $\supp \mu \not= S$.
\end{proof}

Given a measured solenoid $S_\mu$ we can talk about ``$\mu$-almost
all leaves'' with respect to the transversal measure. More
precisely, a Borel set of leaves has zero $\mu$-measure if the
intersection of this set with any local transversal $T$ is a set of
$\mu_T$-measure zero.

Now Poincar\'e recurrence theorem for classical dynamical systems
translates as:

\begin{proposition} \label{cor:Poincare-recurrence} \textbf{\em (Poincar\'e recurrence)}
Let $S_\mu$ be an ergodic oriented $1$-solenoid with $\supp \mu =S$.
Then $\mu$-almost all leaves are dense and accumulate themselves.
\end{proposition}

\begin{proof}
For each local transversal $T$ we know by proposition \ref{prop:5.7}
that the Poincar\'e return map $R_T$ is
defined for $\mu_T$-almost every point and leaves invariant $\mu_T$.
Therefore by Poincar\'e recurrence theorem, $\mu_T$-almost every
point has a dense orbit by $R_T$ in $\supp \mu_T=T$.

Observe that $S_\mu$ ergodic implies that $S$ is connected
(otherwise we may decompose the invariant measure by restricting it to each
connected component).

By compactness we can choose a finite number of local transversals
$T_i=\varphi^{-1}(\{ 0\}\x K(U_i))$ with flow-boxes $\{U_i\}$
covering $S$. We can assume that we have that $U_i\cap U_j$ is a
flow-box if non-empty. This and connectedness of $S$ imply that any
Borel set of leaves that has either total or zero measure in a flow-box
$U_i$, has the same property in $S$.

Now, the set of leaves $S_i$ that are non-dense in a given
flow-box $U_i$ is of zero $\mu$-measure in $U_i$ (by Poincar\'e
recurrence theorem applied to $R_{T_i}$). By the above,
$S_i$ is of zero $\mu$-measure in $S$. Finally the set of
non-dense leaves in $S$ is the finite union of the $S_i$,
therefore is a set of leaves of zero $\mu$-measure.
\end{proof}

\begin{definition} \label{def:M_T}
We denote by $\cM_\cT (S)$ the space of tranversal positive measures
on the solenoid $S$ equipped with the topology generated by weak
convergence in each local tranversal.

We denote by $\overline\cM_\cT (S)$ the quotient of $\cM_\cT(S)$
by positive scalar multiplication.
\end{definition}

\begin{proposition}\label{prop:extremal-ergodic-transversal}
The space $\cM_\cT (S)$ is a cone in the vector space of all
transversal signed measures $\cV_\cT (S)$. 
Extremal measures of $\cM_\cT (S)$ correspond to ergodic tranversal measures.
\end{proposition}

\begin{proof}
Only the last part needs a proof. If $(\mu_T)$ is not ergodic, then
there exists a local tranversal $T_0$ and two disjoint Borel set
$A,B\subset T_0$ invariant by holonomy with $\mu_{T_0}(A)\not=0$,
$\mu_{T_0}(B)\not=0$ and
$\mu_{T_0}(A)+\mu_{T_0}(B)=\mu_{T_0}({T_0})$. Let $S_A\subset S$,
resp.\ $S_B\subset S$, be the union of leaves of the solenoid
intersecting $A$, resp.\ $B$. These are Borel subsets of $S$. Let
$(\mu_{T |T\cap S_A})$ and $(\mu_{T |T\cap S_B})$ be the transversal
measures conditional to $T\cap S_A$, resp.\ $T\cap S_B$. Then
 $$
 (\mu_T)=(\mu_{T |T\cap S_A})+(\mu_{T |T\cap S_B}) \, ,
 $$
and $(\mu_T)$ is not extremal.
\end{proof}

\begin{corollary}
If $\cM_\cT (S)$ is non-empty then $\cM_\cT (S)$ contains ergodic
measures.
\end{corollary}

\begin{proof} The proof follows from the application of Krein-Milman
theorem after the identification of $\overline\cM_\cT (S)$ to a
compact convex set of a locally convex topological vector space given below by
theorem \ref{thm:transverse-riemannian}.
\end{proof}

\section{Riemannian solenoids}\label{sec:riemannian}

All measures considered are Borel measures and all limits of
measures are understood in the weak-* sense. We denote by $\cM (S)$
the space of probability measures supported on $S$.

\begin{definition} \label{def:riemannian}\textbf{\em (Riemannian solenoid)}
Let $S$ be a $k$-solenoid of class $C^{r,s}$ with $r\geq 1$. A
Riemannian structure on $S$ is a Riemannian structure on the leaves
of $S$ such that there is a foliated manifold $W$ defining the
solenoid structure of $S$ and a metric $g_W$ on the leaves of $W$ of
class $C^{r-1,s}$ such that $g=g_{W|S}$.
\end{definition}

As for compact manifolds, any solenoid can be endowed with a
Riemannian structure. In the rest of this section $S$ denotes a
Riemannian solenoid.

Note that a Riemannian structure defines a $k$-volume on leaves of
$S$.

\begin{definition}\textbf{\em (Flow group)} \label{def:flow-group}
We define the flow group $G_S^0$ of a Riemannian $k$-solenoid $S$ as
the group of $k$-volume preserving homeomorphisms of $S$ isotopic to
the identity in the group of homeomorphisms of $S$. The full group
consisting of $k$-volume preserving homeomorphisms of $S$ (non
necessarily isotopic to the identity) is denoted by $G_S$.
\end{definition}

\begin{definition}\label{def:desintegrate}\textbf{\em (daval
measuers)} Let $\mu$ be a measure supported on $S$. The measure
$\mu$ is a daval measure if it desintegrates as volume along leaves
of $S$, i.e.\ for any flow-box $(U,\varphi)$ with local transversal
$T=\varphi^{-1}(\{0\}\x K(U))$, we have a measure $\mu_{U,T}$
supported on $T$ such that for any Borel set $A\subset U$
 $$
 \mu(A)=\int_{T} {\Vol}_k(A_y) \ d\mu_{U,T}(y) \, ,
 $$
where $A_y=A\cap \varphi^{-1} (D^k\times \{ y \} )\subset U$.

We denote by $\cM_\cL (S)\subset \cM (S)$ the space of probability
daval measures.
\end{definition}

It follows from this definition that the measures $(\mu_{U,T})$ do
indeed define a transversal measure as we prove in the next
proposition.

\begin{proposition} \label{prop:desintegrate} Let $\mu$ be a daval measure on $S$.
Then we have the following properties:
 \begin{enumerate}
 \item[(i)] For a local transversal $T$, the measures $\mu_{U,T}$ do not depend on
 $U$. So they define a unique measure $\mu_T$ supported on $T$
 \item[(ii)] The measures $(\mu_T)$ are uniquely determined by $\mu$.
 \item[(iii)] The measures $(\mu_T)$ are locally finite.
 \item[(iv)] The measures $(\mu_T)$ are invariant by holonomy and
 therefore define a transversal measure.
 \end{enumerate}
\end{proposition}

\begin{proof}
For (i) and (ii) notice that for any Riemannian metric $g$ we have,
denoting by $B_\epsilon^g (y)$ the Riemannian ball of radius $\epsilon$ around $y$
in its leaf,
 $$
 \lim_{\epsilon \to 0} \frac{\Vol_k(B_\epsilon^g (y))}
 {\epsilon^k}=c(k) \, ,
 $$
where $c(k)$ is a constant only depending on $k$. Therefore by
dominated convergence we have for any Borel set $C\subset T$
 $$
 \mu_{U,T}(C)=\lim_{\epsilon \to 0} \int_T
 \frac{{\Vol}_k(B_\epsilon^g (y))}{c(k) \epsilon^k} \ d\mu_{U,T}(y)
 =\lim_{\epsilon \to 0} \frac{\mu(V_{\epsilon} (C))}{c(k) \epsilon^k}
 \, ,
$$
where $V_\epsilon(C)$ denotes the $\epsilon$-neighborhood of $C$ along
leaves. The last limit is clearly independent of $U$, thus
$\mu_{U,T}$ is independent of $U$ as claimed, and $\mu_T$ is
uniquely determined by $\mu$.

\medskip

For (iii) observe that for each flow-box $U$ we have that
 $$
 y\mapsto {\Vol}_k (L_y)\, ,
 $$
$L_y=\varphi^{-1}(D^k\x\{y\})$, is continuous on $y\in T$, therefore
for any compact subset $C\subset T$ exists $\epsilon_0
>0$ such that for all $y\in C$,
 $$
 {\Vol}_k (L_y) \geq \epsilon_0 \, .
 $$
Let $V=\varphi^{-1}(D^k\x C)$, then we have
 \begin{equation}\label{eqn:formula}
 \mu (V)=\int_C {\Vol}_k (L_y) \ d\mu_T(y) \geq \epsilon_0 \
 \mu_T (C) \, ,
 \end{equation}
therefore $\mu_T(C) <+\infty $.

\medskip

Regarding (iv), consider a flow-box $(U,\varphi)$ and two local
transversals $T_1$ and $T_2$ in $U$ of the form
$T_i=\varphi^{-1}(\{x_i\} \x K(U))$, $i=1,2$, $x_i\in D^k$. These
transversals are associated to flow-boxes $(U,\varphi_i)$ with the
same domain $U$. There is a local holonomy map in $U$, $h:T_1 \to
T_2$. For any Borel set $A\subset U$, we have by definition
 $$
 \int_{T_1} {\Vol}_k(A_y) \ d\mu_{U,T_1}(y)=\mu(A)=\int_{T_2}
 {\Vol}_k(A_{y'}) \ d\mu_{U,T_2}({y'}) \, .
 $$
On the other hand, the change of variables, $y'=h(y)$, gives
 $$
 \int_{T_1} {\Vol}_k(A_y) \ d\mu_{U,T_1}(y)=\int_{T_2}
 {\Vol}_k(A_{y'}) \ dh_* \mu_{U,T_1}(y') \, .
 $$
Thus for any Borel set $A\subset U$,
 $$
 \int_{T_2} {\Vol}_k(A_{y'}) \ d\mu_{U,T_2}({y'})=\int_{T_2}
 {\Vol}_k(A_{{y'}}) \ dh_* \mu_{U,T_1}({y'}) \, .
 $$
And taking horizontal Borel sets, this implies
 $$
 \mu_{U,T_2}=h_* \mu_{U,T_1} \, .
 $$

The invariance by local holonomy implies the invariance by all
holonomies. Take two arbitrary local transversals $T_1'$ and $T_2'$,
two points $p_1\in T_1'$, $p_2\in T_2'$ in the same leaf, and a path
$\gamma$ from $p_1$ to $p_2$ inside a leaf. Then it is easy to
construct a flow-box $(U,\varphi)$ around the curve $\gamma$, so
that $T_1''=\varphi^{-1}(\{x_1\} \x K(U))\subset T_1'$ and
$T_2''=\varphi^{-1}(\{x_2\} \x K(U))\subset T_2'$ ($x_1$ and $x_2$ being
two distinct points of $D^k$) are open subsets
of the respective transversals and $\gamma$ is fully contained in a
leaf of $U$.
\end{proof}

{}From this it follows that Riemannian solenoids do not necessarily
have daval measures (i.e. $\cM_\cL(S)$ can be empty), because there
are solenoids that do not admit transversal measures (see \cite{Pl}
for interesting examples).

\begin{proposition} \label{prop:M_L-compact}
The space of probability daval measures $\cM_\cL (S)$ is a compact
convex set in the vector space of signed measures $\cV(S)$.
\end{proposition}

\begin{proof}
The convexity is clear, and by compactness of $\cM(S)$ we only need
to show that $\cM_\cL (S)$ is closed which follows from the more
precise lemma that follows.
\end{proof}

\begin{lemma}\label{lem:desintegrate}
Let $(\mu_n)$ be a sequence of measures on $S$ that desintegrate as
volume on leaves in a flow-box $U$. Then any limit $\mu$
desintegrates as volume on leaves in $U$ and the tranversal measure
is the limit of the transversal measures.
\end{lemma}

\begin{proof}
We assume that $\mu_n\to \mu$. Given the transversal $T$, the
tranversal measures $(\mu_{n,T})$ are locally finite by proposition
\ref{prop:desintegrate}. Moreover, formula (\ref{eqn:formula}) shows
that they are uniformly locally finite. Extract (with a diagonal
process) a converging subsequence $\mu_{n_k,T}$ to a locally finite
measure $\mu_T$. For any vertically compactly supported continuous
function $\varphi$ defined on $U$ and depending only on $y\in T$,
 $$
 \int_U \varphi \ d\mu_{n_k}=\int_T \varphi \Vol_k(L_y)\ d\mu_{n_k,T}(y)\, ,
 $$
with $L_y=\varphi^{-1}(D^k\x \{y\})$. Passing to the limit $k\to
+\infty$,
 $$
 \int_U \varphi \ d\mu =\int_T \varphi \Vol_k(L_y)\ d\mu_{T}(y) \, .
 $$
Therefore the limit measure $\mu$ desintegrates as volume on leaves
in $U$ with transversal measure $\mu_T$. Since $\mu_T$ is uniquely
determined by $\mu$ (by proposition \ref{prop:desintegrate}), the
only limit of the transversal measures is $\mu_T$.
\end{proof}

\begin{theorem}\label{thm:G-invariance}
A finite measure $\mu$ on $S$ is $G_S^0$-invariant if and only if
$\mu$ is a daval measure.
\end{theorem}

\begin{proof}
If $\mu$ desintegrates as volume along leaves, then it is clearly
invariant by a transformation in $G_S^0$ close to the identity as is
seen in each flow-box. Then it is $G_S^0$-invariant since a
neighborhood of the identity in $G_S^0$ generates $G_S^0$.

Conversely, assume that $\mu$ is a $G_S^0$-invariant finite measure.
We must prove that in any flow-box $(U,\varphi)$ we
have $\mu=\Vol_k \x \mu_{K(U)}$. We can find a map $h: D^k\x K(U)
\to \RR^k\x K(U)$ of class $C^{r,s}$, preserving leaves, and such
that it takes the $k$-volume for the Riemannian metric to the
Lebesgue measure on $\RR^k$. On $h(D^k\x K(U))$, we still denote by
$\mu$ the corresponding measure. We can desintegrate $\mu= \{\nu_y\}
\x \eta$, where $\eta$ is supported on $K(U)$ and $\nu_y$ is a
measure on each horizontal leaf, parametrized by $y\in K(U)$ (see
\cite{Bou,Die}), i.e.
 \begin{equation}\label{eqn:eta}
 \int_U \varphi \ d\mu=\int_{K(U)} \left ( \int_{L_y} \varphi \
 d\nu_y\right ) \ d\eta (y)\, .
 \end{equation}

The group $G_S^0$ in this chart contains all small translations.
Each translation must leave invariant $\eta$-almost all measures
$\nu_y$. Therefore a countable number of translations leave
invariant $\eta$-almost all measures $\nu_y$. Now observe that if $\tau_n$ are
translations leaving invariant $\nu_y$,  and $\tau_n\to \tau$, then $\tau$
leaves $\nu_y$ invariant. Thus taking a countable and dense set of
translations of fixed small displacement, and taking limits, it
follows that all small translations leave invariant $\nu_y$ for
$\eta$-almost all $y$. By Haar theorem these measures are
proportional to the Lebesgue measure, $\nu_y=c(y){\Vol}_k$. We have
that $c\in L^1(K(U),\eta)$ by applying (\ref{eqn:eta}) with
$\varphi$ being the characteristic function of a sub-flow-box with
horizontal leaves being balls of fixed $k$-volume. We define the
transversal measure $\mu_{K(U)}$ as
 $$
 d\mu_{K(U)} =c \ d\eta \, .
 $$
Therefore $\mu=\Vol_k \times \mu_{K(U)}$ on $U$, hence $\mu$ is a daval measure.
\end{proof}

\begin{theorem} \label{thm:transverse-riemannian}
\textbf{\em (Tranverse measures of the Riemannian solenoid)} There
is a one-to-one correspondence between transversal measures $(\mu_T)$
and finite daval measures $\mu$. Furthermore, there is an
isomorphism
  $$
  \overline{\cM}_\cT(S) \cong \cM_\cL(S)\, .
  $$
\end{theorem}

\begin{proof}
The open sets inside flow-boxes form a basis for the Borel
$\sigma$-algebra, and the formula
 $$
 \mu(A)=\int_T {\Vol}_k(A_y) \ d\mu_{T}(y) \, ,
 $$
for $A$ in a flow-box $U$ with local transversal $T$, is compatible
for different flow-boxes. So it defines a measure $\mu$. This
measure is finite because by compactness we can cover $S$ by a
finite number of flow-boxes with finite mass. By construction, $\mu$
is a daval measure. The converse was proved earlier in proposition
\ref{prop:desintegrate}.

This correspondence is clearly a topological isomorphism.
\end{proof}

\begin{definition}\label{def:volume}\textbf{\em (Volume of measured solenoids)}
For a measured Riemannian solenoid $S_\mu$ we define the volume
measure of $S$ as the unique probability measure (also denoted by
$\mu$) associated to the transversal measure $\mu=(\mu_T)$ by theorem
\ref{thm:transverse-riemannian}.
\end{definition}

For uniquely ergodic Riemannian solenoids $S$, this volume measure
is uniquely determined by the Riemannian structure (as for a compact
Riemannian manifold). We observe that, contrary to what happens with
compact manifolds, there is no canonical total mass normalization of
the volume of the solenoid depending only on the Riemannian metric.
This is the reason why we normalize $\mu$ to be a probability
measure.

\begin{definition}\label{def:controlled-growth}\textbf{\em (Controlled growth solenoids)}
Let $S$ be a Riemannian solenoid. Fix a leaf $l\subset S$ and an
exhaustion $(C_n)$ by subsets of $l$. For a flow-box $(U,\varphi)$
write
 $$
 C_n \cap U = A_n\cup B_n \, ,
 $$
where $A_n$ is composed by all full disks $L_y=\varphi^{-1}
(D^k\times \{y\})$ contained in $C_n$, and $B_n$ contains those
connected components $B$ of $C_n\cap U$ such that $B\not= L_{y}\cap
U$ for any $y$. The solenoid $S$ has controlled growth with respect
to $l$ and $(C_n)$ if for any flow-box $U$ in a finite covering of
$S$
 $$
 \lim_{n\to +\infty } \frac{{\Vol}_k(B_n) }{{\Vol}_k(A_n)}=0\, .
 $$

 The solenoid $S$ has controlled growth if $S$ contains a leaf $l$ and an
exhaustion $(C_n)$ such that $S$ has controlled growth with respect
to $l$ and $(C_n)$.
\end{definition}

For a Riemannian solenoid $S$, it is natural to consider the
exhaustion by balls $B(x_0,R_n)$ in a leaf $l$ centered at a point
$x_0\in l$ and with $R_n\to +\infty$, and test the controlled growth
condition with respect to such exhaustions.

The controlled growth condition depends a priori on the Riemannian
metric. As we see next, it guarantees the existence of daval
measures, hence the existence of transversal measures on $S$. Indeed
the measures we construct are Schwartzman measures defined as:

\begin{definition}\label{def:Schwartzman limit}\textbf{\em (Schwartzman
limits and measures)} We say that a measure $\mu$ is a Schwartzman
measure if it is obtained as the Schwartzman limit
 $$
 \mu=\lim_{n\to +\infty } \mu_{n}\, ,
 $$
where the measures $(\mu_n)$ are the normalized $k$-volume of the
exhaustion $(C_n)$ with uniformly bounded total mass. We denote by
$\cM_\cS (S)$ the space of probability Schwartzman measures.
\end{definition}

Compactness of probability measures show:

\begin{proposition}
There are always Schwartzman measures on $S$,
 $$
 \cM_\cS(S)\not= \emptyset \, .
 $$
\end{proposition}

\begin{theorem}\label{thm:desintegration-Schwartzman-measures}
If $S$ is a solenoid with controlled growth, then any Schwartzman
measure is a daval measure,
 $$
 \cM_\cS (S)\subset \cM_\cL (S) \, .
 $$
In particular, $\cM_\cL (S)\not= \emptyset$ and $S$ admits
transversal measures.
\end{theorem}

\begin{proof}
Let $\mu_{n}\to \mu$ be a Schwartzman limit as in definition
\ref{def:Schwartzman limit}. For any flow-box $U$ we prove that
$\mu$ desintegrates as volume on leaves of $U$. Since $S$ has
controlled growth, pick a leaf and an exhaustion which satisfy the
controlled growth condition. Let
 $$
 C_n\cap U = A_n\cup B_n \, ,
 $$
be the decomposition for $C_n\cap U$ described before. The set $A_n$
is composed of a finite number of horizontal disks. We define a new
measure $\nu_n$ with support in $U$ which is the restriction of
$\mu_{n}$ to $A_n$, i.e. it is proportional to the $k$-volume on
horizontal disks. The measure $\nu_n$ desintegrates as volume on
leaves in $U$. The transversal measure is a finite sum of Dirac
measures. Moreover the controlled growth condition implies that
$(\nu_n)$ and $(\mu_{n|U})$ must converge to the same limit. But we
know that $\cM_\cL(S)$ is closed, thus the limit measure $\mu_{|U}$
desintegrates on leaves in $U$. So $\mu$ is a daval measure.
\end{proof}

\begin{corollary}\label{cor:volume}
The volume $\mu$ of a uniquely ergodic solenoid with controlled
growth is the unique Schwartzman measure. Therefore there is only
one Schwartzman limit
 $$
 \mu=\lim_{n\to +\infty} \mu_n \, ,
 $$
which is independent of the leaf and the exhaustion.
\end{corollary}

\begin{proof} There are always Schwartzman limits.
Theorem \ref{thm:desintegration-Schwartzman-measures} shows that any
such limit $\mu$ desintegrates as volume on leaves. Thus the measure
$\mu$ defines the unique (up to scalars) transversal measure
$(\mu_T)$. But, conversely, the transversal measure determines the
measure $\mu$ uniquely. Therefore there is only possible limit
$\mu$, which is the volume of the uniquely ergodic solenoid.
\end{proof}

Following the proof of theorem
\ref{thm:desintegration-Schwartzman-measures} we can be more
precise. We first define irregular measures. These are measures
which have no mass that desintegrates as volume along leaves.

\begin{definition}
Let $\mu$ be a measure supported on $S$. We say that $\mu$ is
irregular if for any Borel set $A\subset S$ and for any non-zero
measure $\nu\in \cM_\cL (S)$ we do not have
 $$
 \nu_{|A} \leq \mu_{|A} \, .
 $$
\end{definition}

\begin{theorem}
Let $\mu$ be any measure supported on $S$. There is a unique
canonical decomposition of $\mu$ into a regular part
$\mu_r\in\cM_\cL (S)$ and an irregular part $\mu_i$,
 $$
 \mu=\mu_r +\mu_i \, .
 $$
We can define the regular part by
 $$
 \mu_r (A)=\sup_\nu \nu (A) \leq \mu (A) \, ,
 $$
for any Borel set $A\subset S$, where the supremum runs over all
measures $\nu\in \cM_\cL (S)$, with $\nu_{|A} \leq \mu_{|A}$ (if no
such measure exists then $\mu_r (A)=0$).
\end{theorem}

Note that this theorem corresponds to the decomposition of any
measure on a regular manifold into an absolutely continuous part
with respect to a Lebesgue measure and a singular part. Indeed, it
generalizes that decomposition to solenoids, since this theorem
reduces to the classical result when the solenoid is a manifold.

\begin{proof}
Consider all measures $\nu\in\cM_\cL (S)$ such that $\nu\leq \mu$.
We define $\mu_r=\sup \nu$. Considering a countable basis
$(A_i)$ for the Borel $\sigma$-algebra and extracting a triangular
subsequence, we can find a sequence of such measures $(\nu_n)$ such
that $\nu_n (A_i)\to \mu_r(A_i)$, for all $i$, i.e.\ $\nu_n\to
\mu_r$. Since $\cM_\cL(S)$ is closed it follows that
$\mu_r\in\cM_\cL (S)$. By construction, $\mu-\mu_r$ is a positive
measure and irregular. Moreover the decomposition
 $$
 \mu=\mu_r+\mu_i
 $$
is unique, because for another decomposition
 $$
 \mu=\nu_r+\nu_i \, ,
 $$
we have by construction of $\mu_r$,
 $$
 \nu_r\leq \mu_r \, .
 $$
Therefore
 $$
 \nu_i =(\mu_r-\nu_r)+\mu_i \, ,
 $$
and $\mu_i$ being positive this implies that
 $$
 0\leq \mu_r- \nu_r \leq \nu_i \, .
 $$
By definition of irregularity of $\nu_i$, this is only possible if
$\mu_r=\nu_r$, then also $\mu_i=\nu_i$, and the decomposition is
unique.
\end{proof}

\section{Generalized Ruelle-Sullivan currents}\label{sec:Ruelle-Sullivan}

We fix in this section a $C^\infty$ manifold $M$ of dimension $n$.

\begin{definition}\label{def:solenoid-in-manifold}
\textbf{\em (Immersion and embedding of solenoids)} Let $S$ be a
$k$-solenoid of class $C^{r,s}$ with $r\geq 1$. An immersion
 $$
 f:S\to M
 $$
is a regular map (that is, it has an extension $\hat{f}:W\to M$ of
class $C^{r,s}$, where $W$ is a foliated manifold which defines
the solenoid structure of $S$), such that the differential
restricted to the tangent spaces of leaves has rank $k$ at every
point of $S$. We denote by $(f,S)$ an immersed solenoid.


Let $r,s\geq 1$. A transversally immersed solenoid $(f,S)$ is a
regular map $f:S\to M$ such that it admits an extension $\hat{f}:W\to M$ which is
an embedding (of a $(k+l)$-dimensional manifold into an
$n$-dimensional one) of class $C^{r,s}$, such that the images of
the leaves intersect transversally in $M$.

An embedded solenoid $(f,S)$ is a transversally immersed solenoid of class $C^{r,s}$,
with $r,s\geq 1$,
with injective $f$, that is, the leaves do not intersect or
self-intersect.
\end{definition}

Note that under a transversal immersion, resp.\ an embedding,
$f:S\to M$, the images of the leaves are immersed, resp.\
injectively immersed, submanifolds.

A foliation of $M$ can be considered as a solenoid, and the
identity map is an embedding.

\begin{definition}\label{def:Ruelle-Sullivan}\textbf{\em (Generalized currents)}
Let $S$ be an oriented  $k$-solenoid of class $C^{r,s}$, $r\geq
1$, endowed with a transversal measure $(\mu_T)$. An immersion
 $$
 f:S\to M
 $$
defines a real homology class $[f,S_\mu]\in H_k(M,\RR )$ by
duality with differential forms as follows.

Let $\omega $ be an $k$-differential form in $M$. The pull-back
$f^* \omega$ defines a $k$-differential form on the leaves of $S$.
Let $S=\bigcup_i S_i$ be a measurable partition such that each
$S_i$ is contained in a flow-box $U_i$.  We define
 $$
 \la [f,S_\mu],\omega \ra=\sum_i \int_{K(U_i)} \left ( \int_{L_y\cap S_i}
 f^* \omega \right ) \ d\mu_{K(U_i)} (y) \, ,
 $$
where $L_y$ denotes the horizontal disk of the flow-box.
\end{definition}

Note that this definition does not depend on the measurable
partition (given two partitions consider the common refinement).
If the support of $f^*\omega$ is contained in a flow-box $U$ then
 $$
 \la [f,S_\mu],\omega \ra =\int_{K(U)} \left ( \int_{L_y} f^* \omega \right )
 \ d\mu_{K(U)} (y) \, .
 $$
In general, take a partition of unity $\{\rho_i\}$ subordinated to
the covering $\{U_i\}$, then
  $$
   \la [f,S_\mu],\omega \ra = \sum_i
   \int_{K(U_i)} \left( \int_{L_y} \rho_i f^* \omega \right)
   d\mu_{K(U_i)} (y) \, .
  $$
Also for any exact differential $\omega=d\a$ we have
 $$
 \begin{aligned}
 \la [f,S_\mu],d\a\ra =\, &  \sum_i
 \int_{K(U_i)} \left ( \int_{L_y} \rho_i \, f^* d\a \right )
 \ d\mu_{K(U_i)}(y) \\  =\, &  \sum_i
 \int_{K(U_i)} \left ( \int_{L_y} d (\rho_i f^* \a) \right )
 \ d\mu_{K(U_i)}(y)   \\ & \qquad  -
  \sum_i \int_{K(U_i)} \left ( \int_{L_y} d \rho_i \wedge f^* \a \right )
 \ d\mu_{K(U_i)}(y)    = 0\, . \qquad
 \end{aligned}
 $$
The first term vanishes using Stokes in each leaf (the form $\rho_i f^* \a$
is compactly
supported on $U_i$), and the second term vanishes because $\sum_i d\rho_i\equiv 0$.
Therefore $[f, S_\mu]$ is a well defined homology class of degree $k$.

In their original article \cite{RS}, Ruelle and Sullivan defined
this notion for the restricted class of solenoids embedded in $M$.

We can define for the solenoid $S$ homology groups $H_p(S)$ (whose
construction is recalled in \ref{sec:appendix2-generalities}). If
$S$ is an oriented $k$-solenoid and $\mu$ is a transversal
measure, then there is an associated fundamental class $[S_\mu]\in
H_k(S)$. The generalized current is the push-forward of this
fundamental class by $f$.

\begin{proposition}\label{prop:RS-push-forward}
  Let $S_\mu$ be an oriented measured $k$-solenoid. If $f:S\to M$
  is an immersion, we have
    $$
    f_*[S_\mu]=[f,S_\mu] \in H_k(M,\RR) \, .
    $$
\end{proposition}

\begin{proof}
The equality $\la [f, S_\mu],\omega \ra= \la [S_\mu],
f^*\omega\ra$ is clear for any $\omega\in \Omega^k(M)$ (see the
construction of the fundamental class in definition
\ref{def:B.1}). The result follows.
\end{proof}

\bigskip

{}From now on, we shall consider a $C^\infty$ compact and oriented
manifold $M$ of dimension $n$. Let $(f,S_\mu)$ be an oriented
measured $k$-solenoid immersed in $M$.
We aim to construct a $(n-k)$-form representing
  $$
  [f,S_\mu]^* \in H^{n-k}(M,\RR)\, ,
  $$
the dual of $[f,S_\mu]$ under the Poincar{\'e} duality isomorphism $H_k(M,\RR)\cong
H^{n-k}(M,\RR)$.

\medskip

Fix an accessory Riemannian metric $g$ on $M$. This endows $S$
with a solenoid Riemannian metric $f^*g$. We can define the normal
bundle
 $$
 \pi:\nu_f \to S\, ,
 $$
which is an oriented bundle of rank $n-k$, since both $S$ and $M$
are oriented. The total space $\nu_f$ is a (non-compact)
$n$-solenoid whose leaves are the preimages by $\pi$ of the leaves
of $S$.

By \ref{sec:appendix2-generalities}, there is a Thom form $\Phi\in
\Omega^{n-k}_{cv}(\nu_f)$ for the normal bundle. This is a closed
$(n-k)$-form on the total space of the bundle $\nu_f$, with
vertical compact support, and satisfying that
 $$
 \int_{\nu_{f,p}}\Phi=1\, ,
 $$
for all $p\in S$, where $\nu_{f,p}=\pi^{-1}(p)$. Denote by
$\nu_r\subset \nu$ the disc bundle formed by normal vectors of
norm at most $r$ at each point of $S$. By compactness of $S$,
there is an $r_0>0$ such that $\Phi$ has compact support on
$\nu_{r_0}$.

For any $\lambda>0$, let $T_\lambda:\nu_f\to \nu_f$ be the map
which is multiplication by $\lambda^{-1}$ in the fibers. Then set
  $$
  \Phi_r=T_{r/r_0}^*\Phi\, ,
  $$
for any $r>0$. So $\Phi_r$ is a closed $(n-k)$-form, supported in
$\nu_r$, and satisfying
 $$
 \int_{\nu_{f,p}}\Phi_r=1\, ,
 $$
for all $p\in S$. Hence it is a Thom form for the bundle $\nu_f$
as well. By \ref{sec:appendix2-generalities}, $[\Phi_r]=[\Phi]$ in
$H^{n-k}_{cv}(\nu_f)$, i.e. $\Phi_r-\Phi=d\beta$, with $\beta\in
\Omega_{cv}^{n-k-1}(\nu_f)$.

\medskip

Using the exponential map and the immersion $f$, we have a map
 $$
 j:\nu_f \to M\, ,
 $$
given as $j(p,v)=\exp_{f(p)}(v)$, which is a regular map from the
$\nu_f$ (as an $n$-solenoid) to $M$. By compactness of $S$, there are
$r_1,r_2>0$ such that for any disc $D$ of radius $r_2$ contained
in a leaf of $S$, the map $j$ restricted to $\pi^{-1}(D)\cap
\nu_{r_1}$ is a diffeomorphism onto an open subset of $M$. Let us
now define a push-forward map
  $$
  j_* : \Omega^{p}_{cv}(\nu_{r_1}) \to \Omega^p(M)\, .
  $$

Consider first a flow-box $U\cong D^k\x K(U)$ for $S$, where the
leaves of the flow-box are contained in discs of radius $r_2$.
Then
 $$
 \pi^{-1}(U) \cap \nu_{r_1} \cong D^{n-k}_{r_1} \x D^k \x K(U) \,
 ,
 $$
where $D^{n-k}_r$ denotes the disc of radius $r>0$ in $\RR^{n-k}$.
Let $\alpha\in \Omega^{p}_{cv}(\nu_{r_1})$ with support in
$\pi^{-1}(U) \cap \nu_{r_1}$. Then we define
 $$
 j_*\alpha := \int_{K(U)} \big( (j_y)_* (\a|_{D^{n-k}_{r_1}\x D^k \x
 \{y\}} ) \big) d\mu_{K(U)}(y)\, ,
 $$
where $j_y$ is the restriction of $j$ to ${D^{n-k}_{r_1}\x D^k \x
\{y\}} \subset \pi^{-1}(U)\cap \nu_{r_1}$, which is a
diffeomorphism onto its image in $M$. This is the average of the
push-forwards of $\alpha$ restricted to the leaves of $\nu_f$,
using the transversal measure.

Now in general, consider a covering $\{U_i\}$ of $S$ by flow-boxes
such that the leaves of the flow-boxes $U_i$ are contained in
discs of radius $r_2$. Then, for any form $\a\in
\Omega_{cv}^p(\nu_{r_1})$, we decompose $\a=\sum \a_i$ with $\a_i$
supported in $\pi^{-1}(U_i) \cap \nu_{r_1}$. Define
 $$
 j_*\a:=\sum j_*\a_i\, .
 $$
This does not depend on the chosen cover.

\begin{proposition} \label{prop:push-forward}
There is a well defined push-forward linear map
  $$
  j_* : \Omega^{p}_{cv}(\nu_{r_1}) \to \Omega^p(M)\, ,
  $$
such that $dj_*\a=j_*d\a$, and $j_* (\a\wedge \b)= j_*\a \wedge
j_*\beta$, for $\a,\b \in \Omega^{p}_{cv}(\nu_{r_1})$.
\end{proposition}

\begin{proof}
 $j_*d\a=dj_*\a$ holds in flow-boxes, hence it holds globally.
 The other assertion is analogous.
\end{proof}

\medskip

When $M$ is a compact and oriented $n$-manifold, the generalized
current $[f,S_\mu]\in H_k(M,\RR)$ gives an element
 $$
 [f,S_\mu]^* \in H^{n-k}(M,\RR)\, ,
 $$
under the Poincar{\'e} duality isomorphism $H_k(M,\RR)\cong
H^{n-k}(M,\RR)$.

We can construct a form representing the dual of the generalized
current.

\begin{proposition} 
\label{prop:RS-form}
  Let $M$ be a compact oriented manifold.
  Let $(f,S_\mu)$ be a oriented measured solenoid immersed in $M$.
  Let $\Phi_r$ be the Thom form of the normal bundle $\nu_f$
  supported on $\nu_r$, for $0<r<r_1$.
  Then  $j_*\Phi_r$ is a closed $(n-k)$-form representing the dual
  of the
  generalized current,
   $$
   [j_*\Phi_r]= [f,S_\mu]^*\, .
   $$
\end{proposition}

\begin{proof}
As $\Phi_r$ is a closed form, we have
 $$
 dj_*\Phi_r=j_*d\Phi_r=0\, ,
 $$
for $0<r\leq r_1$, so the class $[j_*\Phi_r]\in H^{n-k}(M,\RR)$ is
well-defined.

Now let $r,s$ such that $0<r\leq s<r_1$. Then $[\Phi_r]=[\Phi_s]$
in $H_{cv}^{n-k}(\nu_f)$, so there is a vertically compactly
supported $(n-k-1)$-form $\eta$ with
 \begin{equation}\label{eqn:BETA}
 \Phi_r-\Phi_{s}=d\eta\, .
 \end{equation}
Let $r_3>0$ be such that $\eta$ has support on $\nu_{r_3}$. We can
define a smooth map $F$ which is the identity on $\nu_s$, which
sends $\nu_{r_3}$ into $\nu_{r_1}$ and it is the identity on $\nu_f
- \nu_{2r_3}$. Pulling back (\ref{eqn:BETA}) with $F$, we get
 $$
 \Phi_r-\Phi_s=d (F^*\eta)\, ,
 $$
where $F^*\eta \in \Omega_{cv}^{n-k-1}(\nu_{r_1})$. We can apply
$j_*$ to this equality to get
 $$
 j_*\Phi_r-j_*\Phi_s=d j_*(F^*\eta)\, ,
 $$
and hence $[j_*\Phi_r]=[j_*\Phi_s]$ in $H^{n-k}(M,\RR)$.

\medskip

Now we want to prove that $[j_*\Phi_r]$ coincides with the dual of
the generalized current $[f,S_\mu]^*$. Let $\beta$ be any $k$-form
in $\Omega^k(S)$. Consider a cover $\{U_i\}$ of $S$ by flow-boxes
such that the leaves of each flow-box are contained in discs of
radius $r_2$, and let $\{\rho_i\}$ be a partition of unity
subordinated to this cover. Let $\Phi_i=\rho_i\Phi$, which is
supported on $\pi^{-1}(U_i)\cap \nu_f$, and
 $$
 \Phi_{r,i}=\rho_i\Phi_r= T^*_{r/r_0} \Phi_i
 $$
supported on $\pi^{-1}(U_i)\cap \nu_r$. For $0<\epsilon\leq r_1$,
we have
 $$
 \begin{aligned}
  \int_M j_*\Phi_{\epsilon,i} \wedge \beta &=
  \int_M  \left(  \int_{K(U_i)} (j_{i,y})_* (\Phi_{\epsilon,i|A_y^\epsilon})
  \, d\mu_{K(U_i)}(y)\right)  \wedge\beta \\
  &=   \int_{K(U_i)} \left( \int_M (j_{i,y})_* (\Phi_{\epsilon,i|A_y^\epsilon})
  \wedge\beta \right)\,
  d\mu_{K(U_i)}(y) \\
  &=  \int_{K(U_i)} \left( \int_{A_y^\epsilon} \Phi_{\epsilon,i}\wedge
  j_{i,y}^*\beta \right)\, d\mu_{K(U_i)}(y)\, ,
  \end{aligned}
 $$
where $A_y^\epsilon=D^{n-k}_\epsilon \x D^k\x \{y\}\subset
\pi^{-1}(U_i)$ for $y\in K(U_i)$, and $j_{i,y}=j_{|A_y^\epsilon}$.

In coordinates $(v_1,\ldots, v_{n-k},x_1,\ldots, x_k,y)$ for
$\pi^{-1}(U_i)\cong \RR^{n-k}\x D^k\x K(U_i)$, we can write
 $$
 \Phi=\Phi(v,x,y)= g_0 \,dv_{1}\wedge \cdots \wedge dv_{n-k} +
 \sum_{|I|>0} g_{IJ}\, dx_I\wedge dv_J\, ,
 $$
where $g_0$, $g_{IJ}$ are functions, and $I=\{i_1,\ldots,
i_a\}\subset \{1,\ldots, n-k\}$ and $J=\{j_1,\ldots, j_b\}\subset
\{1,\ldots, k\}$ multi-indices with $|I|=a$, $|J|=b$, $a+b=n-k$.
Pulling-back via $T=T_{\epsilon/r_0}$, we get
 \begin{equation}\label{eqn:Phi}
 \begin{aligned}
 \Phi_\epsilon &= \left(\frac{\epsilon}{r_0}\right)^{-(n-k)}
 \left( (g_0\circ T)
 \,dv_{1}\wedge \cdots \wedge dv_{n-k} + \sum_{|I|>0} \left(\frac{\epsilon}{r_0}\right)^{|I|} (g_{IJ}\circ T)
 \,  dx_I\wedge dv_J\right)  \\
 &=\left(\frac{\epsilon}{r_0}\right)^{-(n-k)}
 \left( (g_0\circ T) \,dv_{1}\wedge \cdots \wedge dv_{n-k} + O(\epsilon)\right) \,
 ,
 \end{aligned}
  \end{equation}
since $|g_{IJ}\circ T|$ are uniformly bounded. Note that the support
of $\Phi_{\epsilon|\RR^{n-k}\x D^k\x \{y\}}$ is inside
$D^{n-k}_\epsilon \x D^k\x \{y\}$.

Also write
 $$
 j_{i,y}^*\beta (v,x) =
 h_0(x,y)\, dx_1\wedge \ldots \wedge dx_k + \sum_{|J|>0}
 h_{IJ}(x,y)\, dx_I\wedge dv_J + O(|v|) \, ,
 $$
and note that $f^*\beta_{|D^k\x \{y\}} = h_0(x,y)\, dx_1\wedge
\ldots \wedge dx_k$.

So
 $$
  \begin{aligned}
  \int_{A_y^\epsilon} & \Phi_{\epsilon,i}\wedge j_{i,y}^*\beta =
  \int_{\RR^{n-k}\x D^k} \rho_i \Phi_{\epsilon}\wedge j_{i,y}^*\beta \\ &=
  \left(\frac{\epsilon}{r_0}\right)^{-(n-k)} \left( \int_{\RR^{n-k} \x D^{k}}
  \rho_i (g_0\circ T) \,dv_{1}\wedge \cdots \wedge dv_{n-k} \wedge
  j_{i,y}^*\beta + O(\epsilon^{n-k+1} ) \right) \\ &=
  \left(\frac{\epsilon}{r_0}\right)^{-(n-k)}  \bigg( \int_{\RR^{n-k}\x D^{k}}
  \rho_i (g_0\circ T) \,dv_{1} \wedge \cdots \wedge dv_{n-k} \wedge
  (h_0\, dx_1\wedge \ldots \wedge dx_k + O(|v|)) \bigg) +  O(\epsilon)
  \\ &=
  \left(\frac{\epsilon}{r_0}\right)^{-(n-k)}  \left( \int_{\RR^{n-k}\x D^{k}}
  \rho_i h_0 (g_0\circ T) \,dv_{1} \wedge \cdots \wedge dv_{n-k} \wedge
  dx_1\wedge \ldots \wedge dx_k  \right)+  O(\epsilon) \\ &=
  \int_{D^{k}}
  \rho_i h_0 \,dx_1\wedge \ldots \wedge dx_k   +  O(\epsilon) \\ &=
  \int_{D^k\x \{y\}} \rho_i  \,
  f^*\beta_{|D^k\x \{y\}}  + O(\epsilon)\, .
  \end{aligned}
 $$
The second equality holds since $|\rho_i|\leq 1$, $|j_{i,y}^*\beta|$
is uniformly bounded, and the support of $\rho_i (g_{IJ}\circ T) \,
dx_I\wedge dv_J \wedge j_{i,y}^*\beta$ is contained inside
$D^{n-k}_\epsilon \x D^k$, which has volume $O(\epsilon^{n-k})$. In
the fourth line we use that $|v|\leq \epsilon$ and
 $$
 \left(\frac{\epsilon}{r_0}\right)^{-(n-k)}  \int_{\RR^{n-k}}
 (g_0\circ T) \,dv_{1} \wedge \cdots \wedge dv_{n-k} =
 \int_{\nu_{f,p}} T^* \Phi = 1.
 $$
The same equality is used in the fifth line.

Adding over all $i$, we get
  $$
  \begin{aligned}
  \la [j_* \Phi_\epsilon], [\b]\ra =
  \int_M j_*\Phi_\epsilon\wedge \b 
  &= \sum_i
  \int_M j_*\Phi_{\epsilon,i} \wedge \beta \\ &=
  \sum_i  \int_{K(U_i)} \left( \int_{A_y^\epsilon} \Phi_{\epsilon,i}\wedge
  j_{i,y}^*\beta \right)\, d\mu_{K(U_i)}(y)
  \\ &= \sum_i \int_{K(U_i)} \left(  \int_{D^k\x \{y\}} \rho_i  \,
  f^*\beta|_{D^k\x \{y\}}  + O(\epsilon) \right)
  d\mu_{K(U_i)}(y) \\ &= \la [f, S_\mu], \b \ra + O(\epsilon)\, .
  \end{aligned}
  $$
Taking $\epsilon\to 0$, we get that
 $$
 [j_*\Phi_r] = \lim_{\epsilon\to 0}\  [j_* \Phi_\epsilon]=[f,S_\mu]^*\,
 ,
 $$
for all $0<r<r_1$.
\end{proof}

\medskip

For $M$ non-compact, we have the isomorphism $H_k(M,\RR)\cong
H^{n-k}_c(M,\RR)$, where $H_c^*(M,\RR)$ denotes compactly
supported cohomology of $M$. Then the generalized current
$[f,S_\mu]$ of an immersed oriented measured
solenoid $(f,S_\mu)$ gives an element
  $$
  [f,S_\mu]^* \in H^{n-k}_c(M,\RR)\, .
  $$
The construction of the proof of proposition \ref{prop:RS-form}
gives a smooth compactly supported form $j_*\Phi_r$ on $M$, for
$r$ small enough, with
  $$
  [j_*\Phi_r]=[f,S_\mu]^* \in H^{n-k}_c(M,\RR)\, .
  $$

\medskip

For $M$ non-oriented, let $\fro$ be the local system defining the
orientation of $M$. Let $(f,S_\mu)$ be an immersed oriented
measured solenoid. Then both $[f,S_\mu]$ and $[j_*\Phi_r]$ are classes
which correspond under the isomorphism
 $$
 H_k(M,\RR) \cong H^{n-k}_c(M,\fro)\, .
 $$
The same proof shows that they are equal.

\begin{theorem}\label{thm:self-intersection} \textbf{ \em
(Self-intersection of embedded solenoids)} Let $M$ be a compact,
oriented, smooth manifold. Let $(f,S_\mu)$ be an embedded oriented
solenoid, such that the transversal measures
$(\mu_T)$ have no atoms. Then we have
 $$
 [f,S_\mu]^* \cup [f,S_\mu]^*=0
 $$
in $H^{2(n-k)}(M,\RR)$.
\end{theorem}

\begin{proof}
If $n-k>k$ then $2(n-k)>n$, therefore the self-intersection is $0$ by degree reasons.
So we may assume $n-k\leq k$.

Let $\b$ be any closed $(n-2(n-k))$-form on $M$. We must
prove that
 $$
 \la [f,S_\mu]^*\cup [f,S_\mu]^*\cup [\b],[M]\ra =0\,,
 $$
where $[M]$ is the fundamental class of $M$.
By proposition \ref{prop:RS-form},
 $$
 \la [f,S_\mu]^*\cup [f,S_\mu]^*\cup [\b],[M]\ra =
 \la [f,S_\mu], j_*\Phi_\epsilon\wedge \b\ra\, ,
 $$
for $\epsilon>0$ small enough.

Consider a covering of $f(S)\subset M$ by open sets $\hat U_i\subset
M$ and another covering of $f(S)$ by open sets $\hat V_i\subset M$
such that the closure of $\hat{V}_i$ is contained in $\hat U_i$. We
may assume that the covering is chosen so that $\{V_i=f^{-1}(\hat
V_i) \}$ satisfies the properties needed for computing
$j_*\Phi_\epsilon$ locally (the auxiliary Riemannian structure is
used). Let $\{\rho_i\}$ be a partition of unity of $S$ subordinated to
$\{V_i \}$ and decompose $\Phi_\epsilon=\sum \Phi_{\epsilon,i}$ with
$\Phi_{\epsilon,i}=\rho_i \, \Phi_\epsilon$. We take $\epsilon>0$
small enough so that $j(\supp \Phi_{\epsilon,i})\subset \hat{U}_i$.
Then
 $$
  \la [f,S_\mu]^*\cup [f,S_\mu]^*\cup [\b],[M]\ra=
 \la [f,S_\mu], j_*\Phi_\epsilon\wedge \b\ra\
  =  \sum_i\la [f,S_\mu], j_*\Phi_{\epsilon,i} \wedge \b\ra\,
 .
 $$

As $f$ is an embedding, we may suppose the open sets
$U_i=f^{-1}(\hat U_i)$ are flow-boxes of $S$. Therefore
 $$
 \la [f,S_\mu], j_*\Phi_{\epsilon,i} \wedge \b \ra =\int_{K(U_i)} \left(\int_{L_y}
 f^*( j_* \Phi_{\epsilon,i}\wedge \beta) \right) d\mu_{K(U_i)}(y)\,
 .
 $$
We may compute
 $$
 \begin{aligned}
 \int_{L_y} f^*(j_* \Phi_{\epsilon,i}\wedge \beta )&= \int_{f(L_y)}
 \left(\int_{K(V_i)}(j_{i,z})_*\Phi_{\epsilon,i}
 \wedge \beta \, d\mu_{K(V_i)}(z) \right) \\ &= \int_{K(V_i)} \left( \int_{f(L_y)}
 (j_{i,z})_*\Phi_{\epsilon,i} \wedge \beta \right) d\mu_{K(V_i)}(z)\, .
 \end{aligned}
 $$
Note that $(j_{i,z})_* \Phi_{i,\epsilon} |_{f(L_y)}$ consists of
restricting the form $\Phi_{i,\epsilon}$ to $\pi^{-1}(L_z)$, the
normal bundle over the leaf $L_z$, then sending it to $M$ via $j$,
and finally restricting to the leaf $f(L_y)$.

Since $f$ is an embedding, we may suppose that in a local chart
$f:U_i=D^k\x K(U_i) \to \hat U_i\subset M$ is the restriction of a
map (that we denote with the same letter) $f: D^k\x B \to \hat U_i$,
where $B\subset \RR^l$ is open and $K(U_i)\subset B$, which in suitable
coordinates for $M$ is written as $f(x,y)=(x,y,0)$. The map $j$
extends to  a map from the normal bundle to the horizontal foliation
of $D^k\x B$, as $j: D^{n-k}_\epsilon \x D^k \x B \to M$,
 $$
 j(v,x,z)= (x_1,\ldots, x_k, z_1+v_1,\ldots, z_l+v_l, v_{l+1},
 \ldots, v_{n-k}) +O(|v|^2)\, .
 $$

Using the formula of $\Phi_{\epsilon}$ given in \eqref{eqn:Phi}, we
have
 $$
  (j_{i,z})_* \Phi_\epsilon (x,y) = \sum_{|I|+|J|=n-k}
  \left(\frac{\epsilon}{r_0}\right)^{|I|-(n-k)}
  (g_{IJ}\circ T)(x, y-z) \,  dx_I\wedge dy_J + O(|y-z|)\, .
 $$
We restrict to $L_y$, and multiply by $\beta$, to get
 $$
 ((j_{i,z})_* \Phi_{\epsilon,i}\wedge \beta)|_{L_y} =
 \sum_{|I|=n-k} (\rho_i \cdot (g_{I0}\circ T))(x,y-z)
 \,  dx_I \wedge \beta
 + O(|y-z|) \, ,
 $$
which is bounded by a universal constant.

Hence
 $$
 |\la [f, S_\mu], j_*\Phi_{\epsilon,i}\wedge\b\ra| \leq
 C_0\ \mu_{K(U_i)}(K(U_i))\ \mu_{K(V_i)}(K(V_i))\leq C_0 \
 \mu_{K(U_i)}(K(U_i))^2 \,,
 $$
where $C_0$ is a constant that is valid for any refinement of the
covering $\{U_i\}$. So
 $$
 |\la [f, S_\mu], j_*\Phi_\epsilon \wedge\b \ra| \leq
 C_0 \sum_i \ \mu_{K(U_i)}(K(U_i))^2\,.
 $$

Observe that $\mu_{K(U_i)}(K(U_i))\leq C_1 \, \mu (U_i)$ and that
$\sum_i \mu (U_i)\leq C_2$ for some positive constants $C_1$ and
$C_2$ independent of the refinements of the covering. Therefore,
 $$
 \begin{aligned}
 |\la [f, S_\mu], j_*\Phi_\epsilon \wedge\b \ra| &\,\leq
 C_0  (\max_i \mu_{K(U_i)} (K(U_i)) )\ C_1 \sum_i \ \mu (U_i) \\
 &\,\leq C_0 C_1 C_2 \max_i \mu_{K(U_i)} (K(U_i))\,.
  \end{aligned}
 $$
When we refine the covering, if the transversal measures have no
atoms, we get that $\max_i \mu_{K(U_i)} (K(U_i)) \to 0$ and then
$$
 \la [f, S_\mu], j_*\Phi_\epsilon \wedge\b \ra =0\,,
 $$
as required.
\end{proof}

Note that for a compact solenoid, atoms of transversal measures must
give compact leaves (contained in the support of the atomic part),
since otherwise at the accumulation set of the leaf we would have a
transversal $T$ with $\mu_T$ not locally finite. In particular if
$S$ is a minimal solenoid which is not a $k$-manifold, then all
transversal measures have no atoms. Therefore, the existence
of transversal measures with atomic part is equivalent to the
existence of compact leaves. This observation gives the following
corollary.

\begin{corollary}\label{cor:self-intersection-no-compact-leaves}
Let $M$ be a compact, oriented, smooth manifold. Let $(f,S)$ be an
embedded oriented solenoid, such that $S$
has no compact leaves. Then for any tranversal measure $\mu$, we
have
 $$
 [f,S_\mu]^* \cup [f,S_\mu]^*=0
 $$
in $H^{2(n-k)}(M,\RR)$.
\end{corollary}

\bigskip

We conclude this section observing that if we want to represent a
homology class $a\in H_k(M,\RR)$ by an immersed solenoid in an
$n$-dimensional manifold $M$ and $a\cup a\neq 0$, then the solenoid
cannot be embedded. Note that when $n-k$ is odd, there is no
obstruction. We shall see later that if 
$n-k$ is odd then we can always obtain a transversally immersed solenoid representing $a$,
and that if $n-k$ is even then we can obtain an immersed solenoid. 

\section{Schwartzman clusters and asymptotic cycles} \label{sec:1-schwartzman}

Let $M$ be a compact  $C^\infty$ Riemannian manifold. Observe that
since $H_1(M,\RR)$ is a finite dimensional real vector space, it
comes equipped with a unique topological vector space structure.

The map $\gamma \mapsto [\gamma ]$ that associates to each loop its
homology class in $H_1(M,\ZZ) \subset H_1(M,\RR)$ is continuous when
the space of loops is endowed with the Hausdorff topology.
Therefore, by compactness, oriented rectifiable loops in $M$ of
uniformly bounded length define a bounded set in $H_1(M,\RR)$.

We have a more precise quantitative version of this result.

\begin{lemma} \label{lem:9.1}
Let $(\g_n)$ be a sequence of oriented rectifiable loops in $M$, and
$(t_n)$ be a sequence with $t_n>0$ and $t_n\to +\infty$. If
 $$
 \lim_{n\to +\infty } \frac{l(\g_n )}{t_n}=0 \, ,
 $$
then in $H_1(M,\RR)$ we have
 $$
 \lim_{n\to +\infty} \frac{[\g_n]}{t_n} =0 \, .
 $$
\end{lemma}

\begin{proof}
Via the map
 $$
 \omega \mapsto \int_\g \omega \, ,
 $$
each loop $\g$ defines a linear map $L_\g$ on $H^1(M,\RR )$ that
only depends on the homology class of $\g$. We can extend this map
to $\RR \otimes H_1(M,\ZZ)$ by
 $$
 c\otimes \g \mapsto c \cdot L_\g \, .
 $$

We have the isomorphism
 $$
 H_1(M,\RR)=\RR \otimes H_1(M,\ZZ) \cong \left ( H^1(M,\RR ) \right )^*
 \, .
 $$
The Riemannian metric gives a $C^0$-norm on forms. We consider the
norm in $H^1(M,\RR )$ given as
 $$
 || [\omega]||_{C^0} =\min_{\omega \in [\omega]} ||\omega|| \, ,
 $$
and the associated operator norm in $H_1(M,\RR )\cong \left (
H^1(M,\RR ) \right )^*$.

We have
 $$
 |L_\gamma([\omega])|= \left | \int_\g \omega \right |\leq l(\g ) ||\o ||_{C^0}
 \leq l(\g ) ||[\o ] ||_{C^0} \, ,
 $$
so
 $$
 ||L_\g || \leq l(\g ) \, .
 $$
Hence $l(\g_{n})/t_n\to 0$ implies $L_{\g_{n}}/t_n \to 0$ which is
equivalent to $[\g_n]/t_n\to 0$.
\end{proof}

\begin{definition}\label{def:asymptotic-cycle} \textbf{\em
(Schwartzman asymptotic $1$-cycles)} \label{def:9.2} Let $c$ be a
parametrized continuous curve $c:{\RR}\to M$ defining an immersion
of $\RR$. For $s, t\in \RR$, $s<t$, we choose a rectifiable
oriented curve $\gamma_{s,t}$ joining $c(s)$ to $c(t)$ such that
 $$
 \lim_{t\to +\infty \atop s\to -\infty}  \frac{l(\g_{s, t} )}{t-s}=0 \
 .
 $$

The parametrized curve $c$ is a Schwartzman asymptotic $1$-cycle
if the juxtaposition of $c|_{[s,t]}$ and $\gamma_{s,t}$, denoted
$c_{s ,t}$ (which is a $1$-cycle), defines a homology class
$[c_{s,t}]\in H_1(M,{\ZZ})$ such that the limit
 \begin{equation}\label{eqn:X}
 \lim_{t\to +\infty \atop s \to -\infty} \frac{[c_{s,t}]}{t-s} \in H_1(M,{\RR})
 \end{equation}
exists.

We define the Schwartzman asymptotic homology class as
 $$
 [c]:=\lim_{t\to +\infty \atop s \to -\infty} \frac{[c_{s,t}]}{t-s} \, .
 $$
\end{definition}

Thanks to lemma \ref{lem:9.1} this definition does not depend on
the choice of the closing curves $(\g_{s,t})$. If we take another
choice $(\g'_{s,t})$, then as homology classes,
 $$
 [c_{s,t}]=[c'_{s,t}]+[\g'_{s,t}-\g_{s,t}]\, ,
 $$
and
 $$
 \frac{l(\g_{s,t}'-\g_{s,t})}{t-s}=
 \frac{l(\g_{s,t}')}{t-s}+\frac{l(\g_{s,t})}{t-s}\to 0\, ,
 $$
as $t\to \infty$, $s\to -\infty$. By lemma \ref{lem:9.1},
 $$
 \lim_{t\to +\infty \atop s \to -\infty}
 \frac{[\g_{s,t}-\g'_{s,t}]}{t-s}=0 \, ,
 $$
thus
 $$
 [c]=\lim_{t\to +\infty \atop s \to -\infty}
 \frac{[c_{s,t}]}{t-s}=\lim_{t\to +\infty \atop s \to -\infty}
 \frac{[c'_{s,t}]}{t-s} \, .
 $$

Note that we do not assume that $c(\RR )$ is an embedding of
$\RR$, i.e. $c(\RR )$ could be a loop. In that case, the
Schwartzman asymptotic homology class coincides with a scalar
multiple (the scalar depending on the parametrization) of the
integer homology class $[c(\RR)]$. This shows that the Schwartzman
homology class is a generalization to the case of immersions
$c:\RR\to M$. More precisely we have:

\begin{proposition} \label{prop:9.3}
If $c: \RR \to M$ is a loop then it is a Schwartzman asymptotic
$1$-cycle and the Schwartzman asymptotic homology class is a
scalar multiple of the homology class of the loop  $[c(\RR )]\in
H_1(M,\ZZ)$.

If $c: \RR \to M$ is a rectifiable loop with its arc-length
parametrization, and $l(c)$ is the length of the loop $c$, then
 $$
 [c]=\frac{1}{l(c)}\, [c(\RR)] \, .
 $$
\end{proposition}

\begin{proof}
Let $t_0>0$ be the minimal period of the map $c:\RR \to M$. Then
 $$
 [c_{s,t}]= \left [\frac{t-s}{t_0} \right ] [c(\RR )] + O (1) \, .
 $$
Then
 $$
 \lim_{t\to +\infty \atop s\to -\infty} \frac{[c_{s,t}]}{t-s}
 =\frac{1}{t_0} [c(\RR )] \, .
 $$
When  $c: \RR \to M$ is the arc-length parametrization of a
rectifiable loop, the period $t_0$ coincides with the length of
the loop.
\end{proof}

We will assume also in the definition of Schwartzman asymptotic
$1$-cycle that we choose $(\g_{s,t})$ such that
$l(\g_{s,t})/(t-s)\to 0$ uniformly and separately on $s$ and $t$
when $t\to +\infty$ and $s\to -\infty$. For simplicity we can
decide to choose always $\g_{s,t}$ with uniformly bounded length,
and even with $\{\g_{s,t} ; s<t\}$ contained in a compact subset
of the space of continua of $M$. Then the uniform boundedness will
hold for any Riemannian metric and the notions defined will not
depend on the Riemannian structure.

\begin{definition}\textbf{\em (Positive and negative asymptotic cycles)}
Under the assumptions of definition \ref{def:9.2}, if the limit
 \begin{equation}\label{eqn:Y}
 \lim_{t\to +\infty} \frac{[c_{s,t}]}{t-s} \in H_1(M,{\RR})
 \end{equation}
exists then it does not depend on $s$, and we say that the
parametrized curve $c$ defines a positive asympotic cycle. The
positive Schwartzman homology class is defined as
 $$
 [c_+]=\lim_{t\to +\infty} \frac{[c_{s,t}]}{t-s} \, .
 $$

The definition of negative asymptotic cycle and negative
Schwartzman homology class is the same but taking $s\to -\infty$,
 $$
 [c_-]=\lim_{s\to -\infty} \frac{[c_{s,t}]}{t-s} \, .
 $$
\end{definition}

The independence of the limit (\ref{eqn:Y}) on $s$ follows from
 $$
  \lim_{t\to +\infty} \frac{[c_{s',t}]}{t-s'} =
  \lim_{t\to +\infty} \frac{[c_{s,t}] + [c_{s',s}] +O(1)}{t-s} \cdot \frac{t-s}{t-s'} =
  \lim_{t\to +\infty} \frac{[c_{s,t}]}{t-s}
 \, .
 $$

\begin{proposition}\label{prop:criterium-Schwartzman-cycle}
A parametrized curve $c$ is a Schwartzman asymptotic $1$-cycle if
and only if it is both a positive and a negative asymptotic cycle
and
 $$
 [c_+]=[c_-] \, .
 $$
In that case we have
 $$
 [c]=[c_+]=[c_-] \, .
 $$
\end{proposition}

\begin{proof}
If $c$ is a Schwartzman asymptotic $1$-cycle, then for $t\to +\infty$ take
$s\to -\infty$ very slowly, say satisfying the relation $t=s^2 l(c_{|[s,0]})$,
which defines $s=s(t)<0$ uniquely as a function of $t>0$. Then
 $$
 \begin{aligned} \
 [c]= & \lim_{t\to \infty \atop s =s(t) \to -\infty} \frac{[c_{s,t}]}{t-s} =
  \lim_{t\to +\infty} \frac{[c_{s,0}] +[c_{0,t}]+O(1)}{t-s} \\
  =&
  \lim_{t\to +\infty} \left(\frac{[c_{s,0}]+O(1)}{t}
   + \frac{[c_{0,t}]}{t} \right)\frac{t}{t-s} = \lim_{t\to +\infty}
 \frac{[c_{0,t}]}{t} \, ,
 \end{aligned}
 $$
since $\frac{t}{t-s} \to 1$ because $\frac{s}{t}\to 0$, and $\frac{[c_{s,0}]}{t}\to
0$ by lemma \ref{lem:9.1}.
So $c$ is a positive asymptotic cycle and $[c]=[c_+]$.
Analogously, $c$ is a negative asymptotic cycle and $[c]=[c_-]$.

Conversely, assume that $c$ is a positive and negative asymptotic
cycle with $[c_+]=[c_-]$. For $t$ large we have
 $$
 \frac{[c_{0,t}]}{ t}=[c_+]+o(1) \, .
 $$
For $-s$ large we have
 $$
 \frac{[c_{s,0}]}{-s}=[c_-]+o(1) \, .
 $$
Now
  $$
  \frac{[c_{s,t}]}{t-s}= \frac{-s}{t-s} \cdot\frac{[c_{s,0}]}{-s} +\frac{t}{t-s} \cdot\frac{[c_{0,t}]}{t} +
  \frac{O(1)}{t-s} = \frac{-s}{t-s} [c_+] + \frac{t}{t-s} [c_-]
  +o(1)\, .
  $$
As $[c_+]=[c_-]$, we get that this limit exists and equals
$[c]=[c_+]=[c_-]$.
\end{proof}

\begin{definition}\textbf{\em (Schwartzman clusters)}
Under the assumptions of definition \ref{def:9.2}, we can
consider, regardless of whether (\ref{eqn:X}) exists or not, all
possible limits
 \begin{equation}\label{eqn:XX}
 \lim_{n\to +\infty} \frac{[c_{s_n,t_n}]}{t_n-s_n} \in H_1(M,\RR ) \, ,
 \end{equation}
with $t_n\to +\infty$ and $s_n\to -\infty$, that is, the derived
set of $([c_{s,t}]/(t-s))_{t\to \infty,s\to -\infty}$. The limits
(\ref{eqn:XX}) are called Schwartzman asymptotic homology classes
of $c$, and they form the Schwartzman cluster of $c$,
 $$
 \cC (c)\subset H_1(M,\RR) \, .
 $$
A Schwartzman asymptotic homology class (\ref{eqn:XX}) is balanced
when the two limits
 $$
 \lim_{n\to +\infty} \frac{[c_{0,t_n}]}{t_n} \in H_1(M,\RR ) \, ,
 $$
and
 $$
 \lim_{n\to +\infty} \frac{[c_{s_n,0}]}{ -s_n} \in H_1(M,\RR ) \, ,
 $$
do exist in $H_1(M,\RR)$. We denote by $\cC_b (c)\subset \cC
(c)\subset H_1(M,\RR)$ the set of those balanced Schwartzman
asymptotic homology classes. The set $\cC_b (c)$ is named the
balanced Schwartzman cluster.

We define also the positive and negative Schwartzman clusters,
$\cC_+ (c)$ and $\cC_- (c)$, by taking only limits $t_n\to
+\infty$ and $s_n\to -\infty$ respectively.
\end{definition}

\begin{proposition}\label{prop:closed-clusters}
The Schwartzman clusters $\cC (c)$, $\cC_+ (c)$ and $\cC_- (c)$
are closed  subsets of $H_1(M,\RR)$.

If $\{[c_{s,t}]/(t-s) ; s<t\}$ is bounded in $H_1(M,\RR )$, then
the Schwartzman clusters $\cC (c)$, $\cC_+ (c)$ and $\cC_- (c)$
are non-empty, compact and connected  subsets of $H_1(M,\RR)$.
\end{proposition}

\begin{proof}
The Schwartzman cluster $\cC (c)$ is the derived set of
 $$
 ( [c_{s,t}]/(t-s) )_{t\to \infty, s\to - \infty} \, ,
 $$
in $H_1(M,\RR)$, hence closed.

Under the boundedness assumption, non-emptiness and compactness
follow. Also the oscillation of $([c_{s,t}])_{s,t}$ is bounded by
the size of $[\g_{s,t}]$. Therefore the magnitude of the
oscillation of $([c_{s,t}]/(t-s))_{s,t}$ tends to $0$ as $t\to
\infty$, $s\to -\infty$. This forces the derived set to be
connected under the boundedness assumption, since it is
$\epsilon$-connected for each $\epsilon >0$. (A compact metric
space is $\epsilon$-connected for all $\epsilon >0$ if and only if
it is connected.)

Also $\cC_+(c)$, resp.\ $\cC_-(c)$, is closed because it is the
derived set of
 $$
 ( [c_{0,t}]/t)_{t\to \infty} \, ,
 $$
resp.\  $$
 ( [c_{s,0}]/(-s))_{s\to - \infty}\, ,
 $$
in $H_1(M,\RR)$. Non-emptiness, compactness and connectedness
under the boundedness assumption follow for the cluster sets
$\cC_\pm (c)$ in the same way as for $\cC (c)$.
\end{proof}

Note that all these cluster sets may be empty if the parametrization
is too fast.

The balanced Schwartzman cluster $\cC_b (c)$ does not need to be
closed, as shown in the following counter-example.

\begin{counterexample} \label{counter}
We consider the torus $M=\TT^2$. We identify $H_1(M,\RR) \cong
\RR^2$, with $H_1(M,\ZZ )$ corresponding to the lattice $\ZZ^2
\subset \RR^2$. Consider a line $l$ in $H_1(M, \RR^2)$ of
irrational slope passing through the origin, $y=\sqrt{2}\  x$ for
example. We can find a sequence of pairs of points $(a_n,b_n)\in
\ZZ^2 \times \ZZ^2$ in the open lower half plane $H_l$ determined
by the line $l$, such that the sequence of segments $[a_n, b_n]$
do converge to the line $l$, and the middle point $(a_n+b_n)/2 \to
0$ (this is an easy exercise in diophantine approximation). We
assume that the first coordinate of $b_n$ tends to $+\infty$, and
the first coordinate of $a_n$ tends to $-\infty$. Now we can
construct a parametrized curve $c$ on $\TT^2$ such that for all
$n\geq 1$ there are an infinite number of times $t_{n,i}\to
+\infty$ with $[c_{0,t_{n,i}}]/t_{n,i}=b_n$, and for an infinite
number of times $s_{n,i}\to -\infty$,
$[c_{s_{n,i},0}]/(-s_{n,i})=a_n$. Thus in homology the curve $c$
oscillates wildly.  We can adjust the velocity of the
parametrization so that $-s_{n,i}=t_{n,i}$. Hence for these times
 $$
 \frac{[c_{s_{n,i},t_{n,i}}]}{t_{n,i}-s_{n,i}} =
 \frac{a_n (-s_{n,i}) + b_n(t_{n,i}) +O(1)}{t_{n,i}-s_{n,i}}
 \to \frac{a_n+b_n}2\, ,
 $$
when $i\to +\infty$, and the two ends balance each other. We have
great freedom in constructing $c$, so that we may arrange to have
always $[c_{s,t}]\subset H_l$. Then we get that $0\in \cC (c)$ and
all $(a_n+b_n)/2\in \cC_b (c)$ but $0\notin \cC_b (c)$.
\end{counterexample}

\medskip

We have that  $c$ is a Schwartzman asymptotic $1$-cycle (resp.\
positive, negative) if and only if $\cC (c)$ (resp.\ $\cC_+(c)$,
$\cC_-(c)$) is reduced to one point. In that case the Schwartzman
asymptotic $1$-cycle is balanced. The next result generalizes
proposition \ref{prop:criterium-Schwartzman-cycle}. We need first
a definition.

\begin{definition}
Let $A,B\subset V$ be subsets of a real vector space $V$. For
$a,b\in V$ the segment $[a,b]\subset V$ is the convex hull of $\{ a,
b\}$ in $V$. The additive hull of $A$ and $B$ is
 $$
 A\widehat + B=\bigcup_{a\in A \atop b\in B } [a,b] \, .
 $$
\end{definition}

\begin{proposition}\label{prop:balanced-cluster}
The Schwartzman balanced cluster $\cC_b (c)$ is contained in the
additive hull of $\cC_+(c)$ and $\cC_-(c)$
 $$
 \cC_b (c)\subset \cC_+(c)\widehat + \, \cC_-(c)\, .
 $$
Moreover, for each $a\in \cC_+(c)$ and $b\in \cC_-(c)$, we have
 $$
 \cC_b (c) \cap [a,b] \not= \emptyset \, .
 $$
\end{proposition}

\begin{proof}
Let $x\in \cC_b (c)$,
 $$
 x=\lim_{n\to +\infty  } \frac{[c_{s_n,t_n}]}{t_n-s_n} \, .
 $$
We write
 $$
 \frac{[c_{s_n,t_n}]}{t_n-s_n}=\frac{[c_{s_n,0}]}
 {-s_n}\cdot\frac{-s_n}{t_n-s_n}
 +\frac{[c_{0,t_n}]}{t_n}\cdot\frac{t_n}{t_n-s_n} +o(1) \, ,
 $$
and the first statement follows.

For the second, consider
 $$
 a=\lim_{n\to +\infty} \frac{[c_{0,t_n}]}{t_n} \in \cC_+(c) \, ,
 $$
and
 $$
 b=\lim_{n\to +\infty} \frac{[c_{s_n,0}]}{-s_n} \in \cC_-(c) \, .
 $$
Then taking any accumulation point $\tau\in [0,1]$ of the sequence
$(t_n/(t_n-s_n))_n \subset [0,1]$ and taking subsequences in the
above formulas, we get a balanced Schwartzman homology class
 $$
 c=\tau a+(1-\tau )b \in \cC_b(c) \, .
 $$
\end{proof}

\begin{corollary} \label{cor:8.11}
If $\cC_+(c)$ and $\cC_-(c)$ are non-empty, then $\cC_b(c)$ is
non-empty, and therefore $\cC (c)$ is also non-empty.
\end{corollary}

Note that we can have $\cC_+(c)=\cC_-(c)=\emptyset$ (then $\cC_b
(c)=\emptyset$) but $\cC (c)\not= \emptyset$ (modify appropriately
counter-example \ref{counter}).

\medskip

There is one situation where we can assert that the balanced
Schwartzman cluster set is closed.

\begin{proposition} \label{prop:9.11}
If $B=\{[c_{s,t}]/(t-s); s<t\} \subset H_1(M,\RR )$ is a bounded
set, then $\cC (c)$, $\cC_+(c)$, $\cC_-(c)$ and $\cC_b (c)$ are
all compact sets. More precisely, they are all contained in the
convex hull of ${\overline B}$.
\end{proposition}

\begin{proof}
Obviously $\cC (c)$, $\cC_+(c)$ and $\cC_-(c)$ are bounded as
cluster sets of bounded sets, hence compact by proposition
\ref{prop:closed-clusters}.

In order to prove that $\cC_b (c)$ is bounded, we observe that the
additive hull of bounded sets is bounded, therefore boundedness
follows from proposition \ref{prop:balanced-cluster}. We show that
$\cC_b (c)$ is closed. Since $\cC_b(c)\subset \cC (c)$ and
$\cC(c)$ is closed, any accumulation point $x$ of $\cC_b(c)$ is in
$\cC (c)$. Let
 $$
 x=\lim_{n\to +\infty } \frac{[c_{s_n,t_n}]}{t_n-s_n} \, ,
 $$
and write as before
 $$
 \frac{[c_{s_n,t_n}] }{t_n-s_n}=
 \frac{[c_{s_n,0}]}{-s_n}\cdot\frac{-s_n}{t_n-s_n} +\frac{[c_{0,t_n}] }{t_n}
 \cdot\frac{t_n}{t_n-s_n} +o(1) \, .
 $$
Note that $([c_{s_n,0}] /(-s_n) )_n$ and $([c_{0,t_n}] /t_n )_n$
stay bounded. Therefore we can extract converging subsequences and
also for the sequence $(t_n/(t_n-s_n))_n\subset [0,1]$. The limit
along these subsequences $t_{n_k}\to +\infty$ and $s_{n_k}\to
-\infty$ give the same Schwartzman homology class $x$ which turns
out to be balanced.

The final statement follows from the above proofs.
\end{proof}

The situation described in proposition \ref{prop:9.11} is indeed
quite natural. It arises each time that $M$ is a Riemannian manifold
and $c$ is an arc-length parametrization of a rectifiable curve. In
the following proposition we make use of the natural norm
$||\cdot||$ in the homology of a Riemannian manifold defined in
\ref{sec:appendix1-norm}.

\begin{proposition}\label{prop:Riemannian-clusters}
Let $M$ be a Riemannian manifold and denote by $||\cdot||$ the
norm in homology. If $c$ is a rectifiable curve parametrized by
arc-length then the cluster sets $\cC (c)$, $\cC_+(c)$, $\cC_-(c)$
and $\cC_b (c)$ are compact subsets of $\bar B(0,1)$, the closed
ball of radius $1$ for the norm in homology.

So $\cC (c)$ and $\cC_\pm(c)$ are non-empty, compact and connected,
and $\cC_b (c)$ is non-empty and compact.
\end{proposition}

\begin{proof}
Observe that we have
 $$
 l(c_{s,t})=l(c|_{[s,t]})+l(\g_{s,t})=t-s+l(\g_{s,t}) \, .
 $$
Thus
 $$
 l([c_{s,t}])\leq t-s+ l(\g_{s,t}) \, .
 $$
By theorem \ref{thm:A.4},
 $$
 || [c_{s,t}]||\leq t-s + l(\g_{s,t}) \, ,
 $$
and
 $$
 \left\| \frac{[c_{s,t}]}{t-s} \right\| \leq 1+\frac{l(\g_{s,t} )}{t-s} \, .
 $$
Since $\frac{l(\g_{s,t} )}{t-s}\to 0$ uniformly, we get that
$B=\{[c_{s,t}]/(t-s); s<t\} \subset H_1(M,\RR )$ is a bounded set.

By proposition \ref{prop:closed-clusters}, $\cC (c)$ and
$\cC_\pm(c)$ are non-empty, compact and connected.
By corollary \ref{cor:8.11}, $\cC_b (c)$ is non-empty and by
proposition \ref{prop:9.11}, it is compact.
\end{proof}

Obviously the previous notions depend heavily on the
parametrization. For a non-parametrized curve we can also define
Schwartzman cluster sets.

\begin{definition}\label{def:8.14.si}
For a non-parametrized oriented curve $c\subset M$, we define the
Schwartzman cluster $\cC (c)$ as the union of the Schwartzman
clusters for all orientation preserving parametrizations of $c$.
We define the positive $\cC_+(c)$, resp.\ negative $\cC_-(c)$,
Schwartzman cluster set as the union of all positive, resp.\
negative, Schwartzman cluster sets for all orientation preserving
parametrizations.
\end{definition}

\begin{proposition} \label{prop:8.15.si}
For an oriented curve $c \subset M$ the Schwartzman clusters $\cC
(c)$, $\cC_+(c)$ and $\cC_-(c)$ are non-empty closed cones of
$H_1(M,\RR)$. These cones are degenerate (i.e. reduced to $\{
0\}$) if and only if $\{ [c_{s,t}] ; s<t \}$ is a bounded subset
of $H_1(M, \ZZ )$.
\end{proposition}

\begin{proof}
We can choose the closing curves $\g_{s,t}$ only depending on
${c}(s)$ and ${c}(t)$ and not on the parameter values $s$ and $t$,
nor on the parametrization. Then the integer homology class
$[{c}_{s,t}]$ only depends on the points ${c}(s)$ and ${c}(t)$ and
not on the parametrization. Therefore, we can adjust the speed of
the parametrization so that $[{c}_{s,t}]/(t-s)$ remains in a ball
centered at $0$. This shows that $\cC (c)$ is not empty. Adjusting
the speed of the parametrization we equally get that it contains
elements that are not $0$, provided that the set $\{[c_{s,t}] ;
s<t\}$ is not bounded in $H_1(M, \ZZ )$. Certainly, if
$\{[c_{s,t}] ; s<t\}$ is bounded, all the cluster sets are reduced
to $\{ 0\}$. Observe also that if $a\in \cC (c)$ then any multiple
$\lambda a$, $\lambda>0$, belongs to $\cC (c)$, by considering the new
parametrization with velocity multiplied by $\lambda$. So $\cC
(c)$ is a cone in $H_1(M ,\RR )$.

Now we prove that $\cC (c) $ is closed. Let $a_n\in \cC (c)$ with
$a_n\to a \in H_1(M,\RR )$. For each $n$ we can choose a
parametrization of $c$, say $c^{(n)}=\tilde{c}\circ \psi_n$ (here
$\tilde{c}$ is a fixed parametrization and $\psi_n$ is an orientation
preserving homeomorphism of $\RR$), and parameters $s_n$ and $t_n$
such that $||[c_{s_n,t_n}^{(n)}]-a|| \leq 1/n$ (considering any
fixed norm in $H_1(M,\RR )$). For each $n$ we can choose $t_n$ as
large as we like, and $s_n$ negative as we like. Choose them
inductively such that $(t_n)$ and $(\psi_n(t_n))$ are both increasing
sequences converging to $+\infty$, and $(s_n)$ and $(\psi_n(s_n))$
are both decreasing sequences converging to $-\infty$. Construct a
homeomorphism $\psi$ of $\RR$ with $\psi(t_n)=\psi_n(t_n)$ and
$\psi(s_n)=\psi_n(s_n)$. It is clear that $a$ is obtained as
Schwartzman limit for the parametrization $\tilde{c}\circ \psi$ at
parameters $s_n,t_n$.

The proofs for $\cC_+(c)$ and $\cC_-(c)$ are similar.
\end{proof}

\begin{remark} \label{rem:8.16}
The image of these cluster sets in the projective space $\PP
H_1(M, \RR)$ is not necessarily connected: On the torus
$M=\TT^2=\RR^2/\ZZ^2$, choose a curve in $\RR^2$ that oscillates
between the half $y$-axis $\{ y>0\}$ and the half $x$-axis $\{
x>0\}$, remaining in a small neighborhood of these axes and being
unbounded for $t\to +\infty$,  and  being bounded when $s\to -\infty$.
Then its Schwartzman cluster consists is two lines through $0$ in
$H_1(\TT^2, \RR)\cong \RR^2$, and its projection in the projective
space consists of two distinct points.
\end{remark}

\begin{remark}\label{rem:8.16bis}
  Let $c$ be a parametrized Schwartzman asymptotic $1$-cycle, and consider the
  unparametrized oriented curve defined by $c$, denoted by $\bar c$. Assume
  that the asymptotic Schwartzman homology class is $a=[c]\neq 0$. Then
   $$
   \cC_\pm(\bar c)=\cC(\bar c)=\RR_{\geq 0}\cdot a\,,
   $$
  as a subset of $H_1(M,\RR)$. This follows since any
  parametrization of $\bar c$ is of the form $c'=c\circ \psi$, where $\psi:\RR\to \RR$
  is a positively oriented homeomorphism of $\RR$. Then
   \begin{equation} \label{eqn:esa}
   \frac{c'_{s,t}}{t-s} =\frac{c_{\psi(s),\psi(t)}}{\psi(t)-\psi(s)} \cdot
   \frac{\psi(t)-\psi(s)}{t-s}\, .
   \end{equation}
   The first term in the right hand side tends to $a$ when $t\to +\infty$, $s\to -\infty$. If the
   left hand side is to converge, then the second term in the right hand side
   stays bounded. After extracting a subsequence, it converges to some $\lambda\geq 0$.
   Hence (\ref{eqn:esa}) converges to $\lambda\, a$.
\end{remark}

\medskip

We define now the notion of asymptotically homotopic curves.

\begin{definition}\textbf{\em (Asymptotic homotopy)} \label{def:8.18}
Let $c_0, c_1: \RR \to M$ be two parametrized curves. They are
asymptotically homotopic if there exists a continuous family
$c_u$, $u\in[0,1]$, interpolating between $c_0$ and $c_1$, such
that
 $$
 c: \RR \times [0,1] \to M \, , \ c(t,u)=c_u(t)\, ,
 $$
satisfies that $\delta_t(u)=c(t,u)$, $u\in [0,1]$  is rectifiable
with
 \begin{equation} \label{eqn:star}
 l(\delta_t)=o(|t|) \, .
  \end{equation}

Two oriented curves are asymptotically homotopic if they have
orientation preserving parametrizations that are asymptotically
homotopic.
\end{definition}

\begin{proposition}\label{prop:8.19}
If $c_0$ and $c_1$ are asymptotically homotopic parametrized curves
then their cluster sets coincide:
 \begin{align*}
 \cC_\pm (c_0)&=\cC_\pm (c_1)\, , \\
 \cC_b (c_0)&=\cC_b (c_1)\, , \\
 \cC (c_0)&=\cC (c_1)\, .
 \end{align*}

If $c_0$ and $c_1$ are asymptotically homotopic oriented curves
then their cluters sets coincide:
\begin{align*}
\cC_\pm (c_0)&=\cC_\pm (c_1)\, , \\
\cC (c_0)&=\cC (c_1)\, .
\end{align*}
\end{proposition}


\begin{proof}
For parametrized curves  we have
$$
[c_{0,s,t}]=[c_{1,s,t}]+[\d_s-\g_{1,s,t}-\d_t+\g_{0,s,t}] \, .
$$
The length of the displacement by the
homotopy is bounded by (\ref{eqn:star}), so
$$
l(\d_s-\g_{1,s,t}-\d_t+\g_{0,s,t})=l(\g_{1,s,t})+l(\g_{0,s,t})+o(|t|+|s|)
\, ,
$$
thus
 $$
 \frac{[c_{0,s,t}]}{t-s}=\frac{[c_{1,s,t}]}{t-s} +o(1) \, .
 $$

For non-parametrized curves, the homotopy between two particular
parametrizations yields a one-to-one correspondence between
points in the curves
 $$
 c_0(t)\mapsto c_1(t) \, .
 $$
Using this correspondence, we have a correspondence of pairs of
points $(a,b)=(c_0(s),c_0(t))$ with pairs of points
$(a',b')=(c_1(s),c_1(t))$. Thus if the sequence of pairs of points
$(a_n,b_n)$ gives a cluster value for $c_0$, then the
corresponding sequence $(a'_n,b'_n)$ gives a proportional cluster
value, since (with obvious notation)
 $$
 [c_{0,a_n,b_n}]=[c_{1,a'_n,b'_n}] +O(1) \, .
 $$
So we can always normalize the speed of the parametrization of
$c_1$ in order to assure that the limit value is the same. This
proves that the clusters sets coincide.
\end{proof}

\section{Calibrating functions}\label{sec:calibrating}

Let $M$ be a $C^\infty$ smooth compact manifold. We define now the
notion of calibrating function.

Let $\pi:\tilde{M}\to M$ be the universal cover of $M$ and let
$\Gamma$ be the group of deck transformations of the cover.

Fix a point $\tilde{x}_0\in \tilde{M}$ and $x_0=\pi(\tilde{x}_0)$.
There is a faithful and transitive action of $\Gamma$ in the fiber
$\pi^{-1}(x_0)$ induced by the action of $\Gamma$ in $\tilde {M}$,
and we have a group isomorphism $\Gamma \cong \pi_1(M,x_0)$. Thus
from the group homomorphism
 $$
 \pi_1(M,x_0)\to H_1(M,{\ZZ}) \, ,
 $$
we get a group homomorphism
 $$
 \rho:\Gamma \to H_1(M,{\ZZ}) \, .
 $$

\begin{definition} \textbf{\em (Calibrating function)}
\label{def:calibrating} A map $\Phi: \tilde{M} \to H_1(M,{\RR})$
is a calibrating function if the diagram
 $$
 \begin{array}{ccc}
  \Gamma\cong\pi_1(M,x_0 )  & \hookrightarrow  & \tilde{M}  \quad \\
   \rho \downarrow  & & \downarrow \Phi\\
   H_1(M,{\ZZ}) & \to & H_1(M,{\RR})
 \end{array}
 $$
is commutative and $\Phi$ is equivariant for the action of
$\Gamma$ on $\tilde {M}$, i.e. for any $g\in \Gamma$ and $\tilde
x\in \tilde {M}$,
 $$
 \Phi (g\cdot \tilde x)=\Phi (\tilde x )+\rho (g) \, .
 $$

 If $\tilde x_0 \in \tilde M$ we say that the calibrating function
 $\Phi$ is associated to $\tilde x_0$ if $\Phi (\tilde x_0 )=0$.
\end{definition}

\begin{proposition}\label{prop:calibrating-existence}
There are smooth calibrating functions associated to any point
$\tilde x_0 \in\tilde M$.
\end{proposition}

\begin{proof}
Fix a smooth  non-negative function $\varphi:\tilde{M}\to {\RR}$
with compact support $K={\overline U}$ with $U=\{\varphi >0\}$
such that $\pi(U)=M$. Moreover, we can request that $U\cap
\pi^{-1}(x_0)=\{\tilde{x}_0\}$.

For any $g_0\in\Gamma$, define $\varphi_{g_0}(\tilde{x})=
\varphi(g_0^{-1} \cdot \tilde{x})$. The support of $\varphi_{g_0}$
is $g_0\, K$, and $(g_0\, K)_{g_0\in \Gamma}$  is a locally finite
covering of $\tilde{M}$, as follows from the compactness of $K$.
Set
 $$
 \psi_{g_0}(\tilde{x}) := \frac{\varphi_{g_0}(\tilde{x})}{\sum_{g\in
 \Gamma} \varphi_{g}(\tilde{x})} \ \, .
 $$
Then $\psi_{g_0}(\tilde{x})=\psi_e(g_0^{-1}\cdot\tilde{x})$ and
 $$
 \sum_{g\in \Gamma} \psi_g \equiv 1 \, .
 $$
Also $\psi_{g_0}$ has compact support $g_0\, K$, and it is a smooth
function since the denominator is strictly positive (because $\pi
(U)=M$) and it is at each point a finite sum of smooth
functions.

We define the map
 $$
 \Phi: \tilde{M} \to H_1(M,{\RR}) \ \, ,
 $$
by
 $$
 \Phi(\tilde{x})= \sum_{g\in \Gamma} \psi_g(\tilde{x})\, \rho(g) \, .
 $$
We check that $\Phi$ is a calibrating function:
 $$
 \begin{aligned}
 \Phi(g\cdot \tilde{x}) &= \sum_{h\in \Gamma}\psi_h(g\cdot\tilde{x})\,\rho(h)\\
    &= \sum_{h\in \Gamma}\psi_{g^{-1}h}(\tilde{x})\,(\rho(g)+\rho(g^{-1}h)) \\
    &= \sum_{h'\in \Gamma}\psi_{h'}(\tilde{x})\,\rho(g) \  + \
       \sum_{h'\in \Gamma}\psi_{h'}(\tilde{x})\, \rho(h') \\
    &= \rho(g)+ \Phi(\tilde{x})\, .
 \end{aligned}
 $$
Notice that by construction $\Phi (\tilde x_0)=0$.
\end{proof}

We note also that choosing a function $\phi$ of rapid decay, we
may do a similar construction, as long as $\sum_{g\in\Gamma}
\phi_g$ is summable (we may need to add a translation to $\Phi$ in
order to ensure $\Phi (\tilde x_0)=0$).

Observe that the calibrating property implies that for a curve
$\gamma:[a,b]\to M$, the quantity
$\Phi(\tilde{\gamma}(b))-\Phi(\tilde{\gamma}(a))$ does not depend
on the lift $\tilde \gamma$ of $\gamma$, because for another
choice $\tilde \gamma '$, we would have for some $g\in \Gamma$,
 $$
 \tilde \gamma'(a) =g\cdot\tilde \gamma (a) \, ,
 $$
and
 $$
 \tilde \gamma'(b) =g\cdot\tilde \gamma (b) \, .
 $$
Therefore
 $$
 \Phi(\tilde{\gamma}'(b))-\Phi(\tilde{\gamma}'(a))
 =\Phi(g\cdot\tilde{\gamma}(b))-\Phi(g\cdot\tilde{\gamma}(a))
 =\Phi(\tilde{\gamma}(b))-\Phi(\tilde{\gamma}(a)) \, .
 $$

This justifies the next definition.

\begin{definition}
Given a calibrating function $\Phi$, for any curve $\gamma:[a,b]\to
M$, we define
$\Phi(\gamma):=\Phi(\tilde{\gamma}(b))-\Phi(\tilde{\gamma}(a))$ for
any lift $\tilde \gamma$ of $\gamma$.
\end{definition}

\begin{proposition} \label{def:8.14}
For any loop $\gamma \subset M$ we have
 $$
 \Phi(\gamma)=[\gamma]\in H_1(M,{\ZZ})\, .
 $$
\end{proposition}

\begin{proof}
Modifying $\gamma$, but without changing its endpoints nor $\Phi
(\gamma)$ nor $[\gamma ]$, we can assume that $x_0\in \gamma$.
Since $\Gamma \cong \pi_1(M,x_0)$, let $h_0\in\Gamma$ be
the element corresponding to $\gamma$. Then $\gamma$ lifts to a curve
joining $\tilde x_0$ to $h_0\cdot \tilde x_0$, and
 $$
 \Phi(\gamma)=\Phi(h_0 \cdot \tilde{x}_0)-
 \Phi (\tilde x_0)
 =\rho(h_0)=[\gamma]\in H_1(M,{\ZZ})\, .
 $$
\end{proof}

\begin{proposition} \label{prop:8.15}
We assume that $M$ is endowed with a Riemannian metric and that the
calibrating function $\Phi$ is smooth. Then for any rectifiable
curve $\gamma$ we have
$$
|\Phi(\gamma)|\leq C \cdot l(\gamma) \, ,
$$
where $l(\gamma )$ is the length of $\gamma$, and $C>0$ is a
positive constant depending only on the metric.
\end{proposition}

\begin{proof}
The calibrating function $\Phi$ is a smooth function on $\tilde{M}$
and $\Gamma$-equivariant, hence it is bounded as well as its
derivatives. The result follows.
\end{proof}

\medskip

\begin{example}
For $M=\TT$, $\tilde{M}={\RR}$, $H_1(M,\ZZ)=\ZZ\subset \RR
=H_1(M,\RR)$, $\Gamma =\ZZ$ and $\rho : \Gamma \to H_1(M,\ZZ)$ is
given (with these identifications) by $\rho (n)=n$. We can take
$\varphi(x)= |1-x|$, for $x\in [-1,1]$, and $\varphi (x)=0$
elsewhere. Then
 $$
 \sum_{n=-\infty}^\infty \varphi(x-n)=1 \, ,
 $$
and
 $$
 \psi_n(x)=\varphi_n(x)=\varphi (x-n) \, .
 $$
Therefore we get the calibrating function
 $$
 \Phi (x)= \sum_{n=-\infty}^\infty \varphi(x-n) \, n= x \,  .
 $$
It is a smooth calibrating function (despite that $\varphi$ is not).

A similar construction works for higher dimensional tori.
\end{example}

\begin{proposition}\label{prop:asymptotic}
Let $c:{\RR}\to M$ be a $C^1$ curve. Consider two sequences $(s_n)$
and $(t_n)$ such that $s_n<t_n$, $s_n \to -\infty$, and $t_n\to
+\infty$.

Then the following conditions are equivalent:
\begin{enumerate}
 \item The limit
 $$
[c]=\lim_{n\to +\infty} \frac{[c_{s_n , t_n}]}{t_n-s_n} \in
H_1(M,{\RR})
 $$
 exists.

\medskip

 \item The limit
  $$
  [c]_{\Phi}=\lim_{n\to \infty} \frac{\Phi(c_{|[s_n,t_n]}) }{t_n-s_n}   \in H_1(M,{\RR})
  $$
  exists.

\medskip

 \item For any closed $1$-form $\alpha \in \Omega^1(M)$, the
  limit
  $$
  [c](\alpha )=\lim_{n\to \infty} \frac{1}{t_n-s_n}\int_{c([s_n , t_n])} \alpha
   $$
  exists.

\medskip

\item For any cohomology class $[\alpha ]\in H^1(M,\RR )$, the limit
  $$
  [c][\alpha ]=\lim_{n\to \infty} \frac{1}{t_n-s_n}\int_{c([s_n , t_n])} \alpha
  $$
exists, and does not depend on the closed  $1$-form $\alpha \in
\Omega^1(M)$ representing the cohomology class.

\medskip

\item For any continuous map $f:M\to \TT$, let $\widetilde{f\circ
c}:{\RR} \to {\RR}$ be a lift of $f\circ c$,  the limit
 $$
 \rho (f)=\lim_{n\to +\infty} \frac{\widetilde{f\circ c}(t_n) -
 \widetilde{f\circ c}(s_n)}{t_n-s_n}
 $$
exists.

\medskip

\item For any (two-sided, embedded, transversally oriented)
hypersurface $H\subset M$ such that all intersections
$c({\RR})\cap H$ are transverse, the limit
  $$
  [c]\cdot [H]=\lim_{n\to \infty}  \frac{\# \{u \in [s_n,t_n]\ | \
  c(u)\in H\}}{t_n-s_n}
  $$
exists. The notation $\#$ means a signed count of intersection
points.

\end{enumerate}

When these conditions hold, we have $[c]=[c]_{\Phi}$ for any
calibrating function $\Phi$. If $\alpha \in \Omega^1 (M)$ is a
closed form, then $[c](\alpha)=[c][\alpha]=\la [c],[\alpha] \ra$.
If $f:M\to \TT$ is a continuous map and $a=f^*[dx]\in H^1(M,\ZZ)$
is the pull-back of the generator $[dx]\in H^1(\TT,\ZZ)$, and $H$
is a hypersurface such that $[H]$ is the Poincar{\'e} dual of $a$,
then $\la [c],a\ra= \rho (f)=[c]\cdot [H]$.
\end{proposition}

\begin{proof}
The equivalence of (1) and (2) follows from the properties of
$\Phi$. Let $c:{\RR}\to M$ be a curve. Then
 $$
 \Phi(c_{|[s_n,t_n]}) = \Phi([c_{s_n,t_n}])- \Phi(\gamma_{s_n,t_n}) =[c_{s_n,t_n}]+O(l(\gamma_{s_n,t_n})) \, .
 $$
Dividing by $t_n-s_n$ and passing to the limit the equivalence of
(1) and (2) follows.

\medskip

We prove that (1) is equivalent to (3). First note that
 $$
  \left| \int_{\gamma_{s_n,t_n}} \alpha  \right| \leq C
   \, l(\g_{s_n,t_n})\,||\alpha||_{C^0}\, .
 $$
We have when $t_n-s_n\to +\infty$,
$$
  \frac{1}{t_n-s_n} \int_{c([s_n,t_n])} \alpha  = \frac{1}{t_n-s_n}\int_{c_{s_n,t_n}} \alpha
  + O\left(\frac{l(\g_{s_n,t_n})}{t_n-s_n}\right)=\frac{[c_{s_n,t_n}](\alpha )}{t_n-s_n}
  +o(1) \, .
$$
and the equivalence of (1) and (3) results.

\medskip

The equivalence of (3) and (4) results from the fact that the limit
  $$
  [c](\alpha)=\lim_{n\to \infty} \frac{1}{t_n-s_n}\int_{c([s_n,t_n])}
  \alpha
  $$
does not depend on the representative of the cohomology class
$a=[\alpha]$. If $\beta=\alpha+d\phi$, with $\phi:M\to {\RR}$
smooth, then $[c](\alpha)=[c](\beta)$ since
  $$
  [c](d\phi)=\lim_{n\to \infty} \frac{1}{t_n-s_n} \int_{c([s_n,t_n])} d\phi
  =\lim_{n\to \infty}  \frac{\phi(c(t_n))-\phi(c(s_n))}{t_n-s_n} \to 0,
  $$
since $\phi$ is bounded. Also $[c][\alpha]=[c](\alpha)$.

\medskip

We turn now to (4) implies (5). First note that there is an
identification $H^1(M,{\ZZ}) \cong [M,K({\ZZ},1)]=[M, \TT]$, where any
cohomology class $[\alpha]\in H^1(M, {\ZZ})$ is associated to a
(homotopy class of a) map $f:M\to \TT$ such that
$[\alpha]=f^*[\TT]$, where $[\TT]\in H^1(\TT,\ZZ)$ is the
fundamental class. To prove (5), assume first that $f$ is smooth.
With the identification $\TT={\RR}/{\ZZ}$, the class $f^*(d
x)=df\in \Omega^1(M)$ represents $[\alpha]$. Therefore
 \begin{equation} \label{eqn:masstar}
 \begin{aligned}
 &\frac{\widetilde{f\circ c}(t_n) -
 \widetilde{f\circ c}(s_n)}{t_n-s_n}= \frac{1}{t_n-s_n} \int_{[s_n,t_n]} d(f\circ
 c)=\\
 & = \frac{1}{t_n-s_n} \int_{[s_n,t_n]} (df)(c') =
 \frac{1}{t_n-s_n} \int_{c([s_n,t_n])} df \, ,
 \end{aligned}
  \end{equation}
and from the existence of the limit in (4) we get the limit in (5)
that we identify as
 $$
 \rho (f)= [c][df] \, .
 $$
If $f$ is only continuous, we approximate it by a smooth function,
which does not change the limit in (5).

Conversely, if (5) holds, then any integer cohomology class admits
a representative of the form $\alpha =df$, where $f: M\to \TT$ is
a smooth map. Then using (\ref{eqn:masstar}) we have
 $$
 \frac{1}{t_n-s_n} \int_{c([s_n,t_n])} \alpha \to \rho(f)\, .
 $$
So the limit in (4) exists for $\alpha=df$.
This implies that the limit in (4) exists for any closed $\alpha \in
\Omega^1 (M)$, since $H^1(M,\ZZ )$ spans $H^1(M,\RR )$.

\medskip

We check the equivalence of (5) and (6). First, let us see that
(6) implies (5). As before, it is enough to prove (5) for a smooth
map $f:M\to \TT$. Let $x_0\in \TT$ be a regular value of $f$, so
that $H=f^{-1}(x_0)\subset M$ is a smooth (two-sided)
hypersurface. Then $[H]$ 
represents the
Poincar{\'e} dual of $[df]\in H^1(M,{\ZZ})$. Choose $x_0$ such that it
is also a regular value of $f\circ c$, so all the
intersections of $c(\RR)$ with $H$ are transverse. Now for any $s<t$, 
 $$
 [c_{s,t}]\cdot [H]= \# c([s,t]) \cap H + \# \gamma_{s,t} \cap
 H\, ,
 $$
where $\#$ denotes signed count of intersection points (we may
assume that all intersections of $\gamma_{s,t}$ and $H$ are
transverse, by a small perturbation of $\g_{s,t}$; also we do not
count the extremes of $\gamma_{s,t}$ in $\# \gamma_{s,t} \cap H$
in case that either $c(s)\in H$ or $c(t)\in H$).

Now
 $$
 \# c ([s,t])\cap H =[\widetilde{f\circ c}(t)] + [ -
 \widetilde{f\circ c}(s)]=\widetilde{f\circ c}(t) -
 \widetilde{f\circ c}(s) + O(1),
 $$
where $[\cdot]$ denotes the integer part, and $| \# \gamma_{s,t}
\cap H|$ is bounded by the total variation of $\widetilde{f\circ
\gamma_{s,t}}$, which is bounded by the maximum of $df$ times the
total length of $\gamma_{s,t}$, which is $o (t-s)$ by assumption.
Hence
 $$
 \lim_{n\to +\infty} \frac{\widetilde{f\circ c}(t_n) -
 \widetilde{f\circ c}(s_n)}{t_n-s_n} =
 \lim_{n\to +\infty}  \frac{\# c([s_n,t_n])\cap H}{t_n-s_n}
 $$
exists.

Conversely, if (5) holds, consider a two-sided embedded
topological hypersurface $H\subset M$. Then there is a collar
$[0,1]\times H$ embedded in $M$ such that $H$ is identified with
$\{\frac12 \}\x H$. There exists a continuous map $f:M\to \TT$
such that $H=f^{-1}(x_0)$ for $x_0=\frac12\in \TT$, constructed by
sending $[0,1]\times H \to [0,1]\to \TT$  and collapsing the
complement of $[0,1]\times H$ to $0$.

Now if all intersections of $c(\RR)$ and $H$ are transverse, that
means that for any $t\in\RR$ such that $c(t)\in H$, we have that
$c(t-\epsilon)$ and $c(t+\epsilon)$ are at opposite sides of the
collar, for $\epsilon>0$ small (the sign of the intersection point
is given by the direction of the crossing). So $f(c(s))$ crosses
$x_0$ increasingly or decreasingly (according to the sign of the
intersection). Hence
 $$
 \frac{\# \{u \in [s_n,t_n] \ | \ c(u)\in H\}}{t_n-s_n} =
 \frac{\widetilde{f\circ c}(t_n) - \widetilde{f\circ
 c}(s_n)}{t_n-s_n}+o(1).
 $$
The required limit exists.
\end{proof}

\medskip

\begin{remark}
Proposition \ref{prop:asymptotic} holds if we only assume the
curve $c$ to be rectifiable.
\end{remark}

\medskip

\begin{corollary}\label{prop:asymptotic2}
Let $c:{\RR}\to M$ be a $C^1$  curve. The following conditions are
equivalent:

\begin{enumerate}
 \item The curve $c$ is a Schwartzman asymptotic cycle.

\medskip

 \item The limit
  $$
  \lim_{t\to +\infty \atop s\to -\infty } \frac{\Phi(c_{|[s,t]})}{t-s}   \in H_1(M,{\RR})
  $$
  exists.

\medskip

 \item For any closed $1$-form $\alpha \in \Omega^1(M)$, the
  limit
  $$
  \lim_{t\to +\infty \atop s\to -\infty } \frac{1}{t-s}\int_{c([s,t])} \alpha
   $$
  exists.

\item For any cohomology class $[\alpha ]\in H^1(M,\RR )$, the limit
  $$
  [c][\alpha ]=\lim_{t\to +\infty \atop s\to -\infty} \frac{1}{t-s}\int_{c([s , t])} \alpha
   $$
  exists, and does not depend on the closed  $1$-form $\alpha \in
  \Omega^1(M)$ representing the cohomology class.

\medskip

 \item For any continuous map $f:M\to \TT$, let $\widetilde{f\circ c}:{\RR}
 \to {\RR}$ be a lift of $f\circ c$, we have that the limit
 $$
 \lim_{t\to +\infty \atop s\to -\infty } \frac{\widetilde{f\circ c}(t) -
 \widetilde{f\circ c}(s)}{t-s}
 $$
 exists.

\medskip

 \item For a (two-sided, embedded, transversally oriented) hypersurface $H\subset M$ such that all
intersections $c({\RR})\cap H$ are transverse, the limit
  $$
  \lim_{t\to +\infty \atop s\to -\infty }  \frac{\# \{u \in [s,t] | c(u)\in H\}}{t-s}
  $$
  exists.
\end{enumerate}

When $c$ is a Schwartzman asymptotic cycle, we have
$[c]=[c]_{\Phi}$ for any calibrating function $\Phi$. If $\alpha
\in \Omega^1 (M)$ is a closed form then
 $$
 [c](\alpha)=[c][\alpha]=\la [\alpha],[c]\ra\,.
 $$
If $f:M\to \TT$ and
$a=f^*[dx]\in H^1(M,\ZZ)$, where $[dx]\in H^1(\TT,\ZZ)$ is the
generator, and $H \subset M$ is a hypersurface such that $[H]$ is
the Poincar\'e dual of $a$, then we have
 $$
 \la [c],[\alpha] \ra =\rho (f)=[c]\cdot [H] \, .
 $$
\end{corollary}

\section{Schwartzman $1$-dimensional cycles}\label{sec:Schwartzman-cycles}

We assume that $M$ is a compact 
$C^\infty$ Riemannian manifold, with Riemannian metric $g$.

\begin{definition} \textbf{\em (Schwartzman representation of homology classes)}
\label{def:representation-Schwartzman-1-solenoid} Let $(f,S)$ be
an immersion in $M$ of an oriented $1$-solenoid $S$. Then $S$ is a
Riemannian solenoid with the pull-back metric $f^* g$.

\begin{enumerate}
\item If $S$ is endowed with a transversal measure
 $\mu=(\mu_T)\in \cM_\cT (S)$, the immersed solenoid $(f,S_\mu)$
 represents an homology class $a\in H_1(M, \RR)$ if for
 $(\mu_T)$-almost all leaves $c:{\RR} \to S$, parametrized
 positively and by arc-length, we have that $f\circ c$ is a
 Schwartzman asymptotic $1$-cycle with $[f\circ c]=a$.

\item The immersed solenoid $(f,S)$ fully represents an homology
class $a\in
 H_1(M, \RR)$ if for all leaves $c:{\RR} \to S$,
 parametrized positively and by arc-length, we have that $f\circ c$ is a
 Schwartzman asymptotic $1$-cycle with $[f\circ c]=a$.
\end{enumerate}

\end{definition}

Note that if $(f,S)$ fully represents an homology class $a\in
H_1(M, \RR)$, then for all oriented leaves $c\subset  S$, we have
that $f \circ c$ is a Schwartzman asymptotic cycle and
 $$
 \cC_+(f \circ c)=\cC_-(f \circ c)
 =\cC(f \circ c)=\RR_{\geq 0}\cdot a \subset H_1(M,\RR) \, ,
 $$
by remark \ref{rem:8.16bis}.

Observe that contrary to what happens with Ruelle-Sullivan cycles,
we can have an immersed solenoid fully representing an homology
class without the need of a transversal measure on $S$.

\begin{definition} \textbf{\em (Cluster of an immersed solenoid)}
\label{def:clusters-Schwartzman-1-solenoid}
Let $(f,S)$ be an immersion in $M$ of an oriented $1$-solenoid
$S$. The homology cluster of $(f,S)$, denoted by $\cC
({{f,S}})\subset H_1(M,\RR)$, is defined as the derived set of 
$([(f\circ c)_{s,t}]/(t-s))_{c,t\to \infty, s\to -\infty}$, taken over all
images of orientation preserving parametrizations $c$ of all leaves of
$S$, and $t\to +\infty$ and $s\to -\infty$. Analogously, we define
the corresponding positive and negative clusters.

The Riemannian cluster of $(f,S)$, denoted by  $\cC^g ({{f,S}})$,
is defined in a similar way, using arc-length orientation
preserving parametrizations. Analogously, we define the positive,
negative and balanced Riemannian clusters.
\end{definition}

\medskip

As in section \ref{sec:1-schwartzman}, we can prove with arguments
analogous to those of propositions \ref{prop:Riemannian-clusters}
and \ref{prop:8.15.si}\ :

\begin{proposition} \label{prop:10.3}
The homology clusters $\cC({{f,S}})$, $\cC_\pm ({{f,S}})$ are
non-empty, closed cones of $H_1(M,\RR)$. If these cones are
non-degenerate, their images in $\PP H_1(M,\RR)$ are non-empty and
compact sets.

The Riemannian homology clusters $\cC^g ({{f,S}})$, $\cC_\pm^g
({{f,S}})$ are non-empty, compact and connected subsets of
$H_1(M,\RR)$.
\end{proposition}

The following proposition is clear, and gives the relationship with
the clusters of the images by $f$ of the leaves of $S$.

\begin{proposition}
Let $(f,S)$ be an immersion in $M$ of an oriented $1$-solenoid
$S$. We have
 $$
 \bigcup_{c\subset S} \cC (f\circ c) \subset \cC ({{f,S}}) \, ,
 $$
where the union runs over all parametrizations of leaves of $S$.
We also have
 $$
 \bigcup_{c\subset S} \cC_\pm (f\circ c) \subset \cC_\pm ({{f,S}}) \, ,
 $$
and
 $$
 \bigcup_{c\subset S} \cC_b (f\circ c) \subset \cC_b ({{f,S}}) \, .
 $$
And similarly for all Riemanniann clusters with $\cC_* (f\circ
c)$ denoting the Schwartzman clusters for the arc-length
parametrization.
\end{proposition}

We recall that given an immersion $(f,S)$ of an oriented
$1$-solenoid, $S$ becomes a
Riemannian solenoid and theorem \ref{thm:transverse-riemannian}
gives a one-to-one correspondence
between the space of transversal measures (up to scalar
normalization) 
and the space of daval measures,
 $$
 \cM_\cT (S) \cong \overline\cM_\cL (S) \, .
 $$

Moreover, in the case of $1$-solenoids that we consider here, 
they do satisfy the controlled growth condition of definition
\ref{def:controlled-growth}. Therefore all Schwartzman measures
desintegrate as length on leaves by theorem
\ref{thm:desintegration-Schwartzman-measures}.

Giving any transversal measure $\mu$ we can consider the associated
generalized current $[f,S_\mu]$.

\begin{definition} We define the Ruelle-Sullivan map
 $$
 \Psi : \cM_\cT (S) \to H_1(M,\RR)
 $$
by
 $$
 \mu \mapsto \Psi (\mu )=[f,S_\mu] \, .
 $$

The Ruelle-Sullivan cluster cone of $(f,S)$ is the image of $\Psi$
 $$
 \cC_{RS} ({{f,S}}) =\Psi (\cM_\cT (S))=\left \{ [f,S_\mu] ; \mu
 \in \cM_\cT (S)\right \} \subset H_1(M,\RR ) \, .
 $$

The Ruelle-Sullivan cluster set is
 $$
 \PP\cC_{RS} ({{f,S}}) \cong \left \{[f,S_\mu] ; \mu \in \cM_\cL (S)\right \} \subset
 H_1(M,\RR ) \, ,
 $$
i.e. using transversal measures which are normalized (using the Riemannian
metric of $M$).
\end{definition}

\medskip

\begin{proposition} \label{prop:11.6}
Let $\cV_\cT(S)$ be the set of all signed measures, with finite
absolute measure and invariant by holonomy, on the solenoid $S$.
The Ruelle-Sullivan map $\Psi$ extended by linearity to
$\cV_\cT(S)$ is a linear continuous operator,
 $$
 \Psi : \cV_\cT (S) \to H_1(M,\RR) \, .
 $$
\end{proposition}

\begin{proof}
Coming back to the definition of generalized current, it is
clear that $\mu\mapsto [f,S_\mu]$ is linear in flow-boxes,
therefore globally. It is also continuous because if $\mu_n\to
\mu$, then $[f,S_{\mu_n}]\to [f,S_\mu ]$ as can be seen in a fixed
flow-box covering of $S$.
\end{proof}

\begin{corollary}
The Ruelle-Sullivan cluster $\cC_{RS} ({{f,S}})$ is a non-empty,
convex, compact cone of $H_1(M,\RR)$. Extremal points of the
convex set $\cC_{RS} ({{f,S}})$ come from the generalized currents
of ergodic measures in $\cM_\cL (S)$.
\end{corollary}

\begin{proof}
Since $\cM_\cL (S)$ is non-empty, convex and compact set, its
image by the continuous linear map $\Psi$ is also a non-empty,
convex and compact set. Any extremal point of $\cC_{RS} ({{f,S}})$
must have an extremal point of $\cM_\cL (S)$ in its pre-image, and
these are the ergodic measures in $\cM_\cL (S)$ (according to the
identification of $\overline\cM_\cL (S)$ to $\cM_\cT (S)$ and by
proposition \ref{prop:extremal-ergodic-transversal}).
\end{proof}

It is natural to investigate the relation between the Schwartzman
cluster and the Ruelle-Sullivan cluster.

\begin{theorem} \label{thm:11.8}
Let $S$ be a $1$-solenoid. For any immersion $f:S\to M$ we have
 $$
 \bigcup_{c\subset S} \cC (f\circ c) \subset \cC_{RS} ({{f,S}}) \, .
 $$
\end{theorem}

\begin{proof}
It is enough to prove the theorem for minimal solenoids, since
each leaf $c\subset S$ is contained in a minimal solenoid
$S_0\subset S$, and
 $$
 \cC (f\circ c) \subset \cC_{RS} (f,S_0)\subset
 \cC_{RS} ({{f,S}}) \, .
 $$
The last inclusion holds because if $\mu$ is a transversal measure
for $S_0$, then it defines a transversal measure $\mu'$ for $S$,
which is clearly invariant by holonomy. Now the generalized
currents coincide, $[f,S_{\mu'}]=[f,S_{0,\mu}]$, as can be seen by
in a fixed flow-box covering of $S$.

The statement for minimal solenoids follows from theorem
\ref{thm:11.9} below.
\end{proof}

\begin{theorem} \label{thm:11.9}
Let $S$ be a minimal $1$-solenoid. For any immersion $f:S\to M$ we have
 $$
 \cC({{f,S}})\subset \cC_{RS} ({{f,S}}) \, .
 $$
\end{theorem}


\begin{proof}
Consider an element $a\in \cC({{f,S}})$ obtained as limit of a
sequence $([(f\circ c_{n})_{s_n,t_n}])$, where $c_n$ is an positively
oriented parametrized leaf of $S$ and $s_n<t_n$, $s_n\to -\infty$,
$t_n\to \infty$. The points $(c_n(t_n))$ must accumulate a
point $x\in S$, and taking a subsequence, we can
assume they converge to it. Choose a
small local transversal $T$ of $S$ at this point, such that
$f(T)\subset B$ where $B\subset M$ is a contractible ball in $M$.
By definition \ref{def:Poincare},
the return map $R_T: T\to T$ is well defined.

Note that we may assume that $\bar T\subset T'$, where $T'$ is also a local
transversal. By compactness of $\bar T$, the return time for $R_{T'}:T'\to T'$
of any leaf, measured with the arc-length parametrization, for any
$x\in \bar T$, is universally bounded. Therefore we can adjust
the sequences $(s_n)$ and $(t_n)$ such that $c_n(s_n)\in \bar T$ and
$c_n(t_n) \in \bar T$, by changing each term by an amount $O(1)$.
Now, after further taking a subsequence, we can arrange that
$c_n(s_n),c_n(t_n) \in T$.

Taking again a subsequence if necessary we can assume that we have a
Schwartzman limit of the measures $\mu_n$ which correspond to the
arc-length on $c_n([s_n,t_n])$ normalized with total mass $1$. The
limit measure $\mu$ desintegrates on leaves because of theorem
\ref{thm:desintegration-Schwartzman-measures}, so it defines a
trasnversal measure $\mu$.

The transversal measures corresponding to $\mu_n$ are atomic,
supported on $T\cap c_n([s_n, t_n))$, assigning the weight
$l([x,R_T(x)])$ to each point in $T\cap c_n([s_n, t_n))$. The
transversal measure corresponding to $\mu$ is its normalized
limit. For each $1$-cohomology class, we may choose a closed
$1$-form $\omega$ representing it and vanishing on $B$ (this is
so because $H^1(M,B)=H^1(M)$, since $B$ is contractible). Assume
that we have constructed $[(f\circ c_{n})_{s_n,t_n}]$ by using
$\gamma_{n,s_n,t_n}$ inside $B$. So
 $$
 \la [f,S_{\mu_n}],\omega \ra = \int_S f^*\omega \, d\mu_n
  =\int_{f\circ c_n([s_n,t_n])} \omega = \la [(f\circ c_{n})_{s_n,t_n}], [\omega]\ra
 \, ,
 $$
thus
 \begin{equation*}
 \la [f,S_\mu],[\omega] \ra = \lim_{n\to \infty}
 \frac{1}{t_n-s_n} \la [f,S_{\mu_n}],\omega \ra =
 \lim_{n\to \infty} \la \frac{[(f\circ c_{n})_{s_n,t_n}]}{t_n-s_n}, [\omega]\ra =\la
 a,[\omega]\ra\, .
 \end{equation*}
Thus the generalized current of the limit measure coincides with
the Schwartzman limit.
\end{proof}

We use the notation $\partial^* C$ for the extremal points of a
compact convex set $C$. For the converse result, we have:

\begin{theorem} Let $S$ be a minimal solenoid and an immersion $f:S\to M$. We have
 $$
 \partial^* \cC_{RS} ({{f,S}}) \subset \bigcup_{c\subset S} \cC (f\circ
 c)\subset \cC({{f,S}}) \, .
 $$
\end{theorem}

\begin{proof}
We have seen that the points in $\partial^* \cC_{RS} ({{f,S}})$
come from ergodic measures in $\cM_\cL (S)$ by the Ruelle-Sullivan
map. Therefore it is enough to prove the following theorem that
shows that the Schwartzman cluster of almost all leaves is reduced
to the generalized current for an ergodic $1$-solenoid.
\end{proof}

\begin{theorem}\label{thm:Ruelle-1-solenoid}
Let $S$ be a minimal $1$-solenoid endowed with an ergodic measure
$\mu \in \cM_\cL (S)$. Consider an immersion
$f:S\to M$. Then for $\mu$-almost all leaves $c\subset S$ we have that
$f\circ c$ is a Schwartzman asymptotic $1$-cycle and
 $$
 [f\circ c]=[f,S_\mu] \in H_1(M, \RR)\, .
 $$
Therefore the immersion  $(f, S_\mu)$ represents its generalized
current.

In particular, this homology class is independent of the metric $g$
on $M$ up to a scalar factor.
\end{theorem}

\begin{proof}
The proof is an application of Birkhoff's ergodic theorem. Choose a
small local transversal $T$ such that $f(T)\subset B$, where
$B\subset M$ is a small contractible ball. Consider the associated
Poincar\'e first return map $R_T: T\to T$. Denote by $\mu_T$ the
transversal measure supported on $T$.

For each $x\in T$ we consider $\varphi_T (x)$ to be the homology
class in $M$ of the loop image by $f$  of the leaf $[x,R_T(x)]$
closed by a segment in $B$ joining $x$ with $R_T(x)$. In this way
we have defined a measurable map 
 $$
 \varphi_T : T \to H_1(M, \ZZ) \, .
 $$
Also for $x\in S$, we denote by $l_T(x)$ the length of the leaf
joining $x$ with its first impact on $T$ (which is  $R_T(x)$ for
$x\in T$).
We have then an upper semi-continuous map 
 $$
 l_T: S \to \RR_+ \, .
 $$
Therefore $l_T$ is bounded by compactness of $S$. In particular,
$l_T$ is bounded on $T$ and thefore in $L^1(T,\mu_T)$. The
boundedness of $l_T$ implies also the boundedness of $\varphi_T$ by
lemma \ref{lem:9.1}.

Consider $x_0\in T$ and its return points $x_i=R_T^i(x_0)$. Let
$0<t_1 < t_2 < t_3 < \ldots $ be the times of return for the
positive arc-length parametrization. We have
 $$
 t_{i+1}-t_i =l_T(x_{i}) \, .
 $$
Therefore
 $$
 t_n=\sum_{i=0}^{n-1} (t_{i+1}-t_i) =\sum_{i=0}^{n-1} l_T\circ R_T^i
 (x_0) \, ,
 $$
and by Birkhoff's ergodic theorem
 $$
 \lim_{n\to +\infty } \frac{1}{n} t_n =\int_T l_T(x) \ d\mu_T (x)
 =\mu (S) =1\, .
 $$

Now observe that, by contracting $B$, we have
 \begin{align*}
 [f\circ c_{0,t_n}] &= [f\circ c_{0,t_1}]+ [f\circ c_{t_1,t_2}]+\ldots
  +[f\circ c_{t_{n-1},t_n}] \\
  &=\varphi_T(x_0) +\varphi_T\circ R_T
 (x_0)+\ldots +\varphi_T \circ R_T^{n-1}(x_0) \, .
 \end{align*}
We recognize a Birkhoff's sum and  by Birkhoff`s ergodic theorem we
get the limit
 $$
 \lim_{n\to +\infty } \frac{1}{n} [f\circ c_{0,t_n}] =\int_T \varphi_T (x) \
 d\mu_T (x) \in H_1(M,\RR)\, .
 $$
Finally, putting these results together,
 $$
 \lim_{n\to +\infty } \frac{1}{t_n} [f\circ c_{0,t_n}]=\lim_{n\to +\infty }
 \frac{[f\circ c_{0,t_n}]/n}{t_n/n} = \frac{ \int_T \varphi_T (x) \ d\mu_T
 (x)}{\int_T l_T(x) \ d\mu_T (x) } = 
 \int_T \varphi_T (x) \ d\mu_T
 (x)   \, .
 $$

Let us see that this equals the generalized current. Take a
closed $1$-form $\omega\in \Omega^1(M)$, which we can assume to
vanish on $B$. Then
 $$
 \la [f,S_\mu],\omega\ra = \int_T \left(\int_{[x,R_T(x)]} f^*\omega
 \right) d\mu_T(x) = \int_T \la \varphi_T(x),\omega\ra d\mu_T(x)\,
 ,
 $$
and so
 $$
 [f,S_\mu] = \int_T \varphi_T(x)\, d\mu_T(x)\, .
 $$

Observe that so far we have only proved that $\cC_+^g(f\circ
c)=\{[f,S_\mu]\}$ for almost all leaves $c\subset S$. Considering
the reverse orientation, the result follows for the negative
clusters, and finally for the whole cluster of almost all leaves.

The last statement follows since $[f,S_\mu]$ only depends on
$\mu\in\cM_\cT(S)$, which is independent of the metric up to scalar
factor, thanks to the isomorphism of theorem \ref{thm:transverse-riemannian}.
\end{proof}

Therefore for a minimal oriented ergodic $1$-solenoid, the generalized
current coincides with the Schwartzman asymptotic homology
class of almost all leaves. It
is natural to ask when this holds for all leaves, i.e. when the
solenoid fully represents the generalized current. This
indeed happens when the solenoid $S$ is uniquely ergodic (recall
that unique ergodicity implies that all orbits are dense and
therefore minimality).

\begin{theorem} \label{thm:11.12} Let $S$ be a uniquely
ergodic oriented $1$-solenoid, and let $\cM_\cL (S)=\{\mu \}$. Let
$f:S\to M$ be an immersion. Then for each leaf $c\subset S$ we
have that $f\circ c$ is a Schwartzman asymptotic cycle with
 $$
 [f\circ c]=[f,S_\mu] \in H_1(M,\RR) \, ,
 $$
and we have
 $$
 \cC^g(f\circ c)=\cC^g ({{f,S}})=\PP\cC_{RS}({{f,S}})=\{ [f,S_\mu]\} \subset
 H_1(M,\RR ) \, .
 $$

Therefore $(f,S)$ fully represents its generalized current
$[f, S_\mu]$.
\end{theorem}

%
%
%
%
%

\section{Schwartzman $k$-dimensional cycles} \label{sec:k-schwartzman}

We study in this section how to extend Schwartzman theory to
$k$-dimensional submanifolds of $M$. We assume that
$M$ is a compact $C^\infty$ Riemannian manifold.

Given an immersion
$c: N \to M$ from an oriented smooth manifold $N$ of dimension
$k\geq 1$, it is natural to consider exhaustions $(U_n)$ of $N$ with
$U_n\subset N$ being $k$-dimensional compact submanifolds with
boundary $\bd U_n$. We close $U_n$ with a $k$-dimensional oriented
manifold $\G_n$ with boundary $\bd \G_n=-\bd U_n$ (that is, $\bd
U_n$ with opposite orientation, so that $N_n=U_n\cup \G_n$ is a
$k$-dimensional compact oriented manifold without boundary), in such
a way that $c_{|U_n}$ extends to a piecewise 
smooth map $c_n:N_n \to M$. We may consider the associated homology
class $[c_n(N_n)]\in H_k(M,\ZZ)$. By analogy with section
\ref{sec:1-schwartzman}, we consider
  \begin{equation}\label{eqn:cn(Nn)}
  \frac{1}{t_n} [c_n(N_n) ] \in H_k(M,\RR ) \, ,
  \end{equation}
for increasing sequences $(t_n)$, $t_n>0$, and $t_n \to +\infty$,
and look for sufficient conditions for (\ref{eqn:cn(Nn)}) to have
limits in $H_k(M,\RR )$. Lemma \ref{lem:9.1} extends to higher
dimension to show that, as long as we keep control of the $k$-volume
of $c_n(\G_n)$, the limit is independent of the closing procedure.

\begin{lemma} \label{lem:closing-k-dim}
Let $(\G_n)$ be a sequence of closed (i.e. compact without boundary)
oriented $k$-dimensional manifolds with piecewise smooth maps $c_n :
\G_n \to M$, and let $(t_n)$ be a sequence with $t_n>0$ and $t_n\to
+\infty$. If
 $$
 \lim_{n\to +\infty } \frac{{\Vol}_k (c_n(\G_n) )}{t_n}=0 \, ,
 $$
then in $H_k(M,\RR)$ we have
 $$
 \lim_{n\to +\infty} \frac{[c_n(\G_n)]}{t_n} =0 \, .
 $$
\end{lemma}

The proof follows the same lines as the proof of lemma
\ref{lem:9.1}. We define now $k$-dimensional Schwartzman asymptotic
cycles.

\begin{definition} \textbf{\em (Schwartzman asymptotic $k$-cycles and clusters)}
\label{def:schwartzman-cluster-k-dim} Let $c:N\to M$ be an immersion
from a $k$-dimensional oriented manifold $N$ into $M$. For all
increasing sequences $(t_n)$, $t_n\to +\infty$, and exhaustions
$(U_n)$ of $N$ by $k$-dimensional compact submanifolds with
boundary,
we consider all possible Schwartzman limits
 $$
 \lim_{n\to +\infty } \frac{[c_n(N_n)]}{t_n} \in H_k(M,\RR) \, ,
 $$
where $N_n=U_n\cup \G_n$ is a closed oriented manifold with
 \begin{equation}\label{eqn:condition}
 \frac{\Vol_k (c_n(\G_n) )}{t_n } \to 0 \, .
 \end{equation}
Each such limit is called a Schwartzman asymptotic $k$-cycle. These
limits form the Schwartzman cluster $\cC ({c,N})\subset H_k(M,\RR)$
of $N$.
\end{definition}

Observe that a Schwartzman limit does not depend on the choice of
the sequence $(\G_n)$, as long as it satisfies
(\ref{eqn:condition}). Note that this condition is independent of
the particular Riemannian metric chosen for $M$.

As in dimension $1$ we have

\begin{proposition}\label{prop:schwartzman-cluster-k-dim}
The Schwartzman cluster $\cC (c,N)$ is a closed cone of
$H_k(M,\RR)$.
\end{proposition}

The Riemannian structure on $M$ induces a Riemannian structure on
$N$ by pulling back by $c$. We define the Riemannian exhaustions
$(U_n)$ of $N$ as exhaustions of the form
  $$
  U_n=\bar B(x_0 , R_n) \, ,
  $$
i.e. the $U_n$ are Riemannian (closed) balls in $N$ centered at a
base point $x_0\in N$ and $R_n\to +\infty$. If the $R_n$ are
generic, then the boundary of $U_n$ is smooth

We define the Riemannian Schwartzman cluster of $N$ as follows. It
plays the role of the balanced Riemannian cluster of section
\ref{sec:1-schwartzman} for dimension $1$.

\begin{definition}\label{def:schwartzman-cluster-k-dim2}
The Riemann-Schwartzman cluster of $(c,N)$, $\cC^g({c,N})$, is the
set of all limits, for all Riemannian exhaustions $(U_n)$,
 $$
 \lim_{n\to +\infty } \frac{1}{\Vol_k (c_n(N_n) )} [c_n(N_n)] \in
 H_k(M,\RR) \, ,
 $$
 such that $N_n=U_n\cup \G_n$ and
 \begin{equation} \label{eqn:condition2}
{\Vol_k (c_n(\Gamma_n) ) \over \Vol_k (c_n(N_n) )}\to 0\, .
  \end{equation}
All such limits are called Riemann-Schwartzman asymptotic
$k$-cycles.
\end{definition}

\begin{definition}\label{def:regular-asymptotic-k-cycles}
The immersed manifold $(c,N)$ represents an homology class $a\in
H_k(M,\RR )$ if the Riemann-Schwartzman cluster $\cC^g({c,N})$
contains only $a$,
 $$
 \cC^g({c,N})=\{ a\} \, .
 $$
We denote $[{c,N}]=a$, and call it the Schwartzman homology class of $(c,N)$.
\end{definition}

Now we can define the notion of representation of homology classes
by immersed solenoids extending definition
\ref{def:representation-Schwartzman-1-solenoid} to higher dimension.

\begin{definition} \textbf{\em (Schwartzman representation of homology classes)}
\label{def:representation-Schwartzman-k-solenoid} Let $(f,S)$ be an
immersion in $M$ of an oriented $k$-solenoid $S$. Then $S$ is a
Riemannian solenoid with the pull-back metric $f^* g$.
\begin{enumerate}
\item If $S$ is endowed with a transversal measure
 $\mu=(\mu_T)\in \cM_\cT (S)$, the immersed solenoid $(f,S_\mu)$
 represents an homology class $a\in H_1(M, \RR)$ if for
 $(\mu_T)$-almost all leaves $l\subset S$, we have that $(f,l)$ is a
 Riemann-Schwartzman asymptotic $k$-cycle with $[{f,l}]=a$.

\item The immersed solenoid $(f,S)$ fully represents an homology
 class $a\in H_1(M, \RR)$ if for all leaves $l\subset S$,
 we have that $(f,l)$ is a
 Riemann-Schwartzman asymptotic $k$-cycle with $[{f,l}]=a$.
\end{enumerate}
\end{definition}

\begin{definition}\textbf{\em (Equivalent exhaustions)}
Two exhaustions $(U_n)$ and $(\hat U_n)$ are equivalent if
 $$
 \frac{\Vol_k (U_n-\hat U_n) +\Vol_k(\hat U_n -U_n )}{\Vol_k (U_n)}\to 0 \, .
 $$
\end{definition}

Note that if two exhaustions $(U_n)$ and $(\hat U_n)$ are equivalent, then
 $$
 \frac{\Vol_k (\hat U_n)}{\Vol_k (U_n)}\to 1 \, .
 $$
Moreover, if $N_n=U_n\cup \Gamma_n$ are closings satisfying
(\ref{eqn:condition2}), then we may close $\hat U_n$ as follows:
after slightly modifying $\hat U_n$ so that $U_n$ and $\hat U_n$
have boundaries intersecting transversally, we glue $F_1=U_n-\hat
U_n$ to $\hat U_n$ along $F_1\cap \bd \hat U_n$, then we glue a copy
of $F_2=\hat U_n- U_n$ (with reversed orientation) to $\hat U_n$
along $F_2\cap \bd \hat U_n$. The boundary of $\hat U_n\cup F_1\cup
F_2$ is homeomorphic to $\bd U_n$, so we may glue $\G_n$ to it, to
get $\hat N_n=\hat U_n\cup F_1\cup F_2 \cup \Gamma_n$. Note that
 $$
 \Vol_k(\hat N_n) =\Vol_k(N_n) + 2\Vol_k (\hat U_n- U_n) \approx \Vol_k(N_n)\,.
 $$
Define $\hat{c}_n$ by $\hat{c}_{n|F_1}= c_{|(U_n-\hat U_n)}$,
$\hat{c}_{n|F_2}= c_{|(\hat U_n- U_n)}$ and $\hat{c}_{n|\G_n}=
c_{n|\G_n}$. Then
 $$
 [c_n(N_n)] = [\hat{c}_n(\hat N_n)]\,,
 $$
so both exhaustions define the same Schwartzman asymptotic
$k$-cycles.

\begin{definition}\textbf{\em (Controlled solenoid)} \label{def:controlled}
Let $V\subset S$ be an open subset of a solenoid $S$.
We say that $S$ is controlled by $V$ if
for any Riemann exhaustion $(U_n)$ of any leaf of $S$ there is an
equivalent exhaustion $(\hat U_n)$ such that for all $n$ we have
$\partial \hat U_n \subset V$. 
\end{definition}

\begin{definition}\textbf{\em (Trapping region)} \label{def:trapping}
An open subset $W\subset S$ of a solenoid $S$ is a trapping region
if there exists a continuous map $\pi : S \to \TT$ such that
\begin{enumerate}
\item[(1)] For some $0<\epsilon_0<1/2$, $W=\pi^{-1} ((-\epsilon_0 ,\epsilon_0))$.

\medskip

\item[(2)]  There is a global transversal $T\subset \pi^{-1} (\{ 0\} )$.

\medskip

\item[(3)] Each connected component of $\pi^{-1} (\{ 0\} )$ intersects $T$ in exactly
one point.

\medskip

\item[(4)] $0$ is a regular value for $\pi$, that is, $\pi$ is smooth
in a neighborhood of $\pi^{-1} (\{ 0\} )$ and it $d\pi$ is surjective
at each point of $\pi^{-1} (\{ 0\} )$ (the differential $d\pi$ is
understood leaf-wise).


\medskip

\item[(5)] For each connected component $L$ of $\pi^{-1} (\TT-\{0\} )$ we
have ${\overline L}\cap T=\{ x,y\}$, where $\{x\} \in {\overline L}\cap T\cap
\pi^{-1} ((-\epsilon_0 ,0])$ and $\{y\} \in {\overline L}\cap T\cap \pi^{-1}
([0,\epsilon_0 ))$. We define $R_T : T\to T$ by $R_T(x)=y$.
\end{enumerate}

\end{definition}

Let $C_x$ be the (unique) component of $\pi^{-1} (\{ 0\} )$ through
$x\in T$. By (4), $C_x$ is a smooth $(k-1)$-dimensional manifold. By
(5), there is no holonomy in $\pi^{-1} ((-\epsilon_0 ,\epsilon_0
))$, so $C_x$ is a compact submanifold. Let $L_x$ be the connected
component of $\pi^{-1} (\TT-\{0\} )$ with ${\overline L_x}\cap T=\{
x,y\}$. This is a compact manifold with boundary
 \begin{equation}\label{eqn:Lx}
 \bd \overline L_x=C_x\cup C_y =C_x\cup C_{R_T(x)} \, .
  \end{equation}

\begin{proposition}\label{prop:trapping}
If $S$ has a trapping region $W$ with global transversal $T$, then
holonomy group of $T$ is generated by the map $R_T$.
\end{proposition}

\begin{proof}
If $\g$ is a path with endpoints in $T$, we may homotop it so that
each time it traverses $\pi^{-1} (\{ 0\} )$, it does it through $T$.
Then we may split $\g$ into sub-paths such that each path has
endpoints in $T$ and no other points in $\pi^{-1} (\{ 0\} )$. Each
of this sub-paths therefore lies in some $\overline L_x$ and has
holonomy $R_T$, $R_T^{-1}$ or the identity. The result follows.
\end{proof}

\begin{theorem}\label{thm:trapping}
A solenoid $S$ with a trapping region $W$ is controlled by $W$.
\end{theorem}

\begin{proof}
Fix a base point $y_0\in S$ and a exhaustion $(U_n)$ of the leaf $l$
through $y_0$ of the form $U_n=\bar B(y_0, R_n)$, $R_n\to +\infty$.
Consider $x_0\in T$ so that $y_0\in \overline
L_{x_0}$. The leaf $l$ is the infinite union
 $$
 l=\bigcup_{n\in \ZZ} \overline L_{R_T^n(x_0)} \, .
 $$
If $R_T^n(x_0)=x_0$ for some $n\geq 1$ then $l$ is a compact manifold.
Then for some $N$, we have $U_N=l$, so the controlled condition of definition
\ref{def:controlled} is satisfied for $l$.

Assume that $R_T(x_0) \not= x_0$. Then $l$ is a non-compact manifold.
For integers $a < b$, denote
  \begin{equation}\label{eqn:Uab}
  \hat{U}_{a,b}:= \bigcup_{k=a}^{b-1} \overline L_{R_T^k(x_0)} \, .
  \end{equation}
This is a manifold with boundary
 $$
 \bd \hat{U}_{a,b}= C_{R_T^a(x_0)}\cup C_{R_T^b(x_0)}\, .
 $$

Given $U_n$, pick the maximum $b\geq 1$ and minimum $a\leq 0$ such
that $\hat{U}_{a,b}\subset U_n$, and denote $\hat U_n=\hat{U}_{a,b}$
for such $a$ and $b$. Clearly $\bd \hat U_n \subset W$. Let us see
that $(U_n)$ and $(\hat U_n)$ are equivalent exhaustions, i.e. that
 $$
 \frac{\Vol_k(U_n-\hat U_n)}{\Vol_k(U_n)} \to 0\,.
 $$

Let $b'\geq 1$ the minimum and $a'\leq 0$ the maximum such that
$U_n\subset \hat{U}_{a',b'}$. Let us prove that
 $$
 \Vol_k(\hat{U}_{a',b'}-\hat{U}_{a,b})
 $$
is bounded. This clearly implies the result.

Take $y\in \overline L_{R^{b'-1}_T(x_0)} \cap U_n$. Then $d(y_0,y)\leq R_n$.
By compactness of $T$, there is a lower bound $c_0>0$ for the distance from
$C_x$ to $C_{R_T(x)}$ in $L_x$, for all $x\in T$.
Taking the geodesic path from $y_0$ to $y$, we see that there are points
in $y_i\in \overline L_{R^{b'-i}_T(x_0)}$ with $d(y_0,y_i)\leq R_n - (i-2) \, c_0$,
for $2\leq i\leq b'$.

As $\overline L_{R^{b}_T(x_0)}$ is not totally contained in $U_n$, we may take
$z\in \overline L_{R^{b}_T(x_0)} -U_n$, so $d(y_0,z)>R_n$.
Both $z$ and $y_{b'-b}$ are on the same leaf $\overline L_{R^{b}_T(x_0)}'$.
By compactness of $T$, the diameter for a leaf $\overline L_x$ is bounded above by
some $c_1>0$, for all $x\in T$. So
 $$
 R_n-(b'-b-2) \, c_0 \geq d(y_0,y_{b'-b}) \geq d(y_0,z) - d(y_{b'-b},z) >R_n - c_1\, ,
 $$
hence
 $$
 b'-b < \frac{c_1}{c_0} +2\, .
 $$
Analogously,
 $$
 a-a' < \frac{c_1}{c_0} +2\, .
 $$

Again by compactness of $T$, the $k$-volumes of $\overline L_x$ are
uniformly bounded by some $c_2>0$, for all $x\in T$. So
 $$
 \Vol_k(\hat{U}_{a',b'}-\hat{U}_{a,b}) \leq (b'-b + a-a') c_2 < 
 2\left(\frac{c_1}{c_0} +2\right)\, c_2 \,,
 $$
concluding the proof.
\end{proof}

\begin{theorem} \label{thm:11.12bis}
Let $S$ be a minimal solenoid endowed with a transversal ergodic
measure $\mu\in\cM_\cL(S)$ and with a trapping region $W\subset S$.
Consider an immersion $f:S\to M$ such that $f(W)$ is contained in a
contractible ball in $M$. Then $(f,S_\mu )$ represents its
generalized current $[f,S_\mu ]$, i.e. for $\mu_T$-almost all
leaves $l\subset S$,
 $$
 [{f,l}]=[f,S_\mu] \in H_k(M,\RR) \, .
 $$

If $S_\mu$ is uniquely ergodic, then $(f,S_\mu)$ fully represents
its generalized current.

In particular, this homology class is independent of the metric $g$
on $M$ up to a scalar factor.
\end{theorem}

\begin{proof}
We define a map $\varphi_T : T\to H_k(M,\ZZ)$ as follows: given
$x\in T$, consider $f(\overline L_x)$. Since $\bd f(\overline L_x)$
is contained in a contractible ball $B$ of $M$, we can close
$f(L_x)$ locally as $N_x =f(\overline L_x)\cup \G_x$ and define an
homology class $\varphi_T (x)=[N_x]\in H_k(M,\ZZ )$. This is
independent of the choice of the closing. This map $\varphi_T$ is
measurable and bounded in $H_k(M,\ZZ)$ since the $k$-volume of
$\G_x$ may be chosen uniformly bounded.  Also we can define a map
$l_T: T\to \RR_+$ by $l_T(x)=\Vol_k(\overline L_x)$. It is also a
measurable and bounded map.

We have seen that every Riemann exhaustion $(U_n)$ is equivalent to
an exhaustion $(\hat U_n)$ with $\partial \hat U_n \subset W$. Note
also that we can saturate the exhaustion $(\hat U_n)$ into $(\hat
U_{n,m})_{n\leq 0\leq m}$, with $\hat U_{n,m}$ defined in
(\ref{eqn:Uab}), where $\partial \hat U_{n,m} = C_{R_T^n(x_0)}\cup
C_{R_T^m(x_0)}$, and $x_0\in T$ is a base point. Since $f(W)$ is
contained in a contractible ball $B$ of $M$, we can always close
$f(\hat U_{n,m})$, with a closing inside $B$, to get $N_{n,m}$
defining an homology class $[N_{n,m}]\in H_k(M,\ZZ )$. Moreover we
have
$$
[N_{n,m}]=\sum_{i=n}^{m-1} \varphi_T(R_T^i(x_0)) \, .
$$
Thus by ergodicity of $\mu$ and Birkhoff's ergodic theorem, we have that for
$\mu_T$-almost all $x_0\in T$,
$$
\frac{1}{m-n} [N_{n,m}]\to \int_T \varphi_T \ d\mu_T \, .
$$

Also
$$
\Vol_k (\hat U_{n,m})=\sum_{i=n}^{m-1} l_T(R_T^i(x_0)) \, ,
$$
where $\Vol_k(N_{n,m})$ differs from $\Vol_k (\hat U_{n,m})$ by a
bounded quantity due to the closings.
By Birkhoff's ergodic theorem, for $\mu_T$-almost all $x_0\in
T$,
$$
\frac{1}{m-n} \Vol_k f(\hat U_{n,m}) \to \int_T l_T \ d\mu_T = \mu(S)=1\, .
$$
Thus we conclude that for $\mu_T$-almost $x_0\in T$,
$$
\frac{1}{\Vol_k (N_{n,m})} [N_{n,m}] \to \int_T \varphi_T \ d\mu_T \, ,
$$
It is easy to see as in theorem \ref{thm:Ruelle-1-solenoid}
that $\int_T \varphi_T \ d\mu_T$ is the generalized current
$[f,S_\mu]$.
\end{proof}

\section{Realization of $H_1(M,\RR)$ } \label{sec:1-solenoids}

Let $M$ be a $C^\infty$ smooth compact Riemannian manifold. Given a
real $1$-homology class $a\in H_1(M,\RR)$, we want to construct an
immersion $(f,S)$ in $M$ of a uniquely ergodic solenoid $S_\mu$
fully representing $a$ up to scalar factor (see definition
\ref{def:representation-Schwartzman-1-solenoid}), and with
generalized current $[f,S_\mu]=a$. By theorem \ref{thm:11.12}, it is
enough to construct an immersed, oriented, uniquely ergodic
$1$-solenoid $(f,S_\mu)$ with $[f,S_\mu]$ equals to a positive
multiple $\lambda \, a$ of $a$, since in this case
$[f,S_{\lambda^{-1}\mu}]=a$. If $\mu$ is the normalized ergodic
measure, theorem \ref{thm:11.12} implies that $(f,S)$ fully
represents $\lambda \, a$. Moreover, if we change the Riemannian
metric of $M$, the property of fully representing $a$ up to positive
scalar factor is preserved (although the scalar factor may vary).

In some situations (depending on the dimension) we will achieve an
embedding. Actually the $1$-solenoid $S$ that we will construct is
independent of $a$ and of $M$, and moreover it has a $1$-dimensional
transversal structure.

Let $h:\TT\to \TT$ be a diffeomorphism with an irrational rotation
number (and therefore uniquely ergodic), which is a Denjoy
counter-example, i.e. has the unique invariant probability measure
supported in a Cantor set $K\subset \TT$. Let $\mu_K$ denote the
invariant probability measure. For the original construction of
Denjoy counter-examples see  \cite{Denjoy}. Actually $h$ can be
taken to be of class $C^{2-\epsilon}$ with $\epsilon>0$ (see
\cite{Herman}).

The suspension of $h$,
  $$
  S_h= ([0,1]\x \TT )_{/ (0,x) \sim (1,h(x))}
  $$
is $C^{2-\epsilon}$-diffeomorphic to the $2$-torus $T^2$. More
explicitly, the diffeomorphism is as follows: take $c>0$ small, let
$h_t$, $t\in [0,c]$, be a (smooth) isotopy from $\id$ to $h$,
then we define the diffeomorphism $H:T^2\to X$ by
 $$
   H(x,t)=\left\{ \begin{array}{ll} (t,h^{-1}(h_t(x))),
   \qquad &\text{for $t\in [0,c]$\, ,} \\
   (t,x), & \text{for $t\in [c,1]$\, .} \end{array} \right.
 $$
Note that $S_h$ is foliated by the horizontal leaves, so $T^2$ is
foliated accordingly. It can be considered also as a $1$-solenoid
of class $C^{\omega,2-\epsilon}$.

The sub-solenoid
  $$
  S= ( [0,1]\x K)_ {/ \sim} \, \subset S_h
  $$
is an oriented $1$-solenoid of class $C^{\omega,2-\epsilon}$, with
transversal $T= (\{0\}\x \TT ) \cap S= \{0\} \x  K$. The holonomy is
given by the map $h$, which is uniquely ergodic. Moreover, the
associated transversal measure is $\mu_K$ on the transversal $K\cong
\{0\}\x K$. So $S$ is an oriented, uniquely ergodic $1$-solenoid.

Using the diffeomorphism $H$, we may see the solenoid $S$ inside the
$2$-torus, $S\subset T^2$, consisting of the paths $(t,x)$, $x\in
K$, $t\in [c,1]$, together with the paths $(t,h_t(x))$, $x\in K$,
$t\in [0,c]$. The embedding $S \inc T^2$ is of class
$C^{\infty,2-\epsilon}$, so we shall think of $S$ as an oriented
$1$-solenoid of regularity $C^{\infty,2-\epsilon}$.

\begin{figure}[h]\label{figure1}
\centering
\resizebox{6cm}{!}{\includegraphics{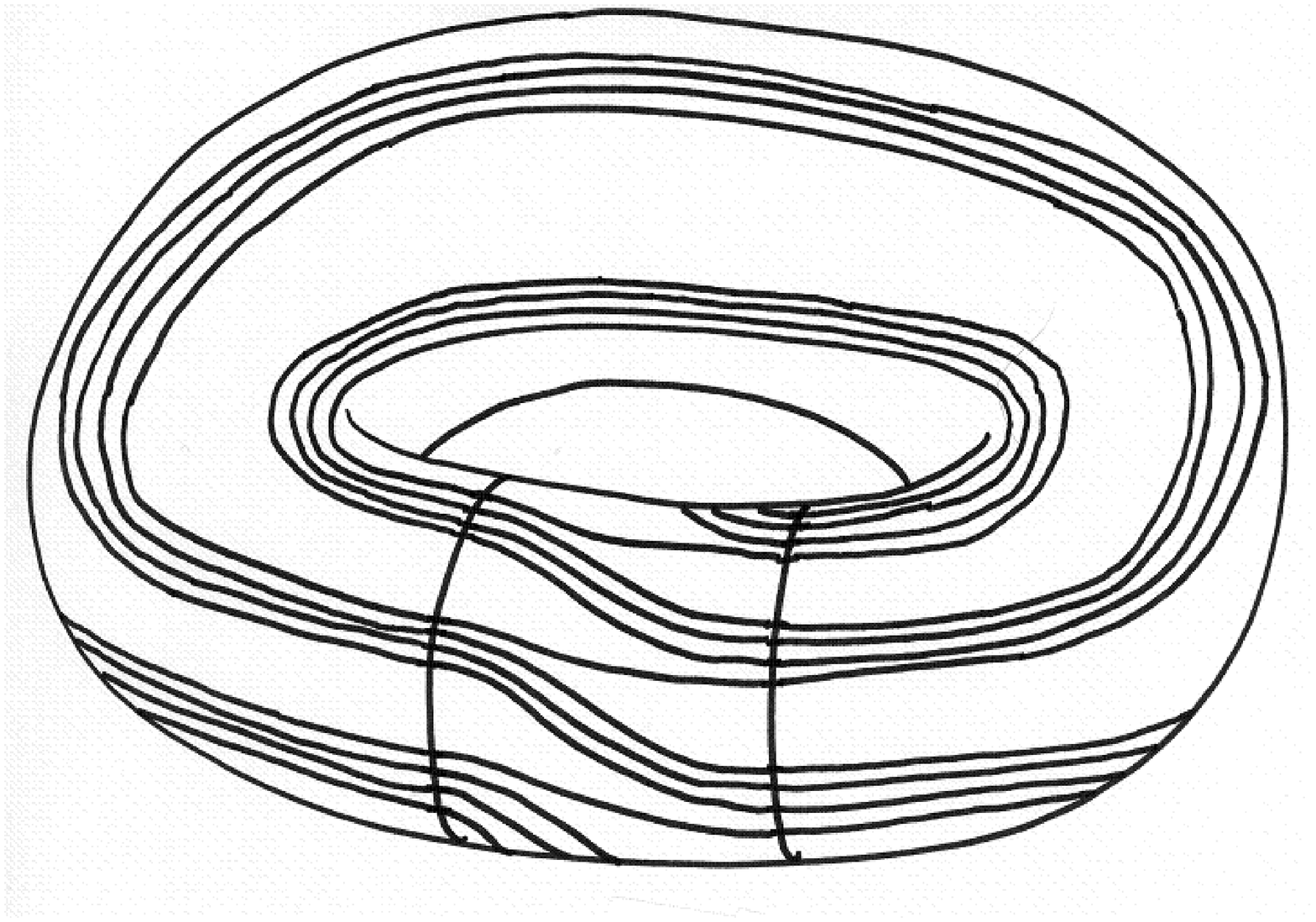}}    
\caption{The $1$-solenoid $S$.}
\end{figure}

\begin{theorem}\label{thm:construct-1-solenoid}
Let $M$ be a compact smooth manifold, and let $a\in H_1(M,{\RR})$ be
a non-zero $1$-homology class. If $\dim M\geq 3$ then (a positive
multiple of) $a$ can be fully represented by an embedding (of class
$C^{\infty,2-\epsilon}$) of the (oriented, uniquely ergodic)
$1$-solenoid $S$ into $M$. If $\dim M=2$ then (a positive multiple
of) $a$ can be fully represented by a transversal immersion of $S$
into $M$.
\end{theorem}

\begin{proof}
Let $C_1,\ldots, C_{b_1}$ be (integral) $1$-cycles which form a basis of the
(real) $1$-homology of $M$. Switch orientations and reorder the
cycles if necessary so that there are real numbers
$\lambda_1,\ldots, \lambda_r>0$ such that
 $$
 a=\lambda_1 C_1 + \cdots+ \lambda_r C_r.
 $$
By dividing by $\sum \lambda_i$ if necessary, we can assume that
$\sum \lambda_i=1$.

Consider the solenoid $S$ constructed above and partition the
cantor set $K$ into $r$ disjoint compact subsets $K_1,\ldots, K_r$
in cyclic order, each of which with
 $$
 \mu_K(K_i)=\lambda_i\, .
 $$
Consider the transversal $T=\{0\}\x\TT$ in $S_h$. We consider angles
$\t_1,\t_2,\ldots, \t_n\in \TT$ in the same cyclic order as the
$K_i$, such that $K_i$ is contained in the open subset $U_i\subset
T$ with boundary points $\t_i$ and $\t_{i+1}$ (denoting
$\t_{n+1}=\t_1$). We may assume that $\t_1=0$. Remove the segments
$[c,1]\x \{\t_i\}$ from $S_h$ to get the open $2$-manifold
 $$
 U=S_h - \cup_i ([c,1]\x \{\t_i\})\, .
 $$

\begin{figure}[h] \label{figure2}
\centering
\resizebox{6cm}{!}{\includegraphics{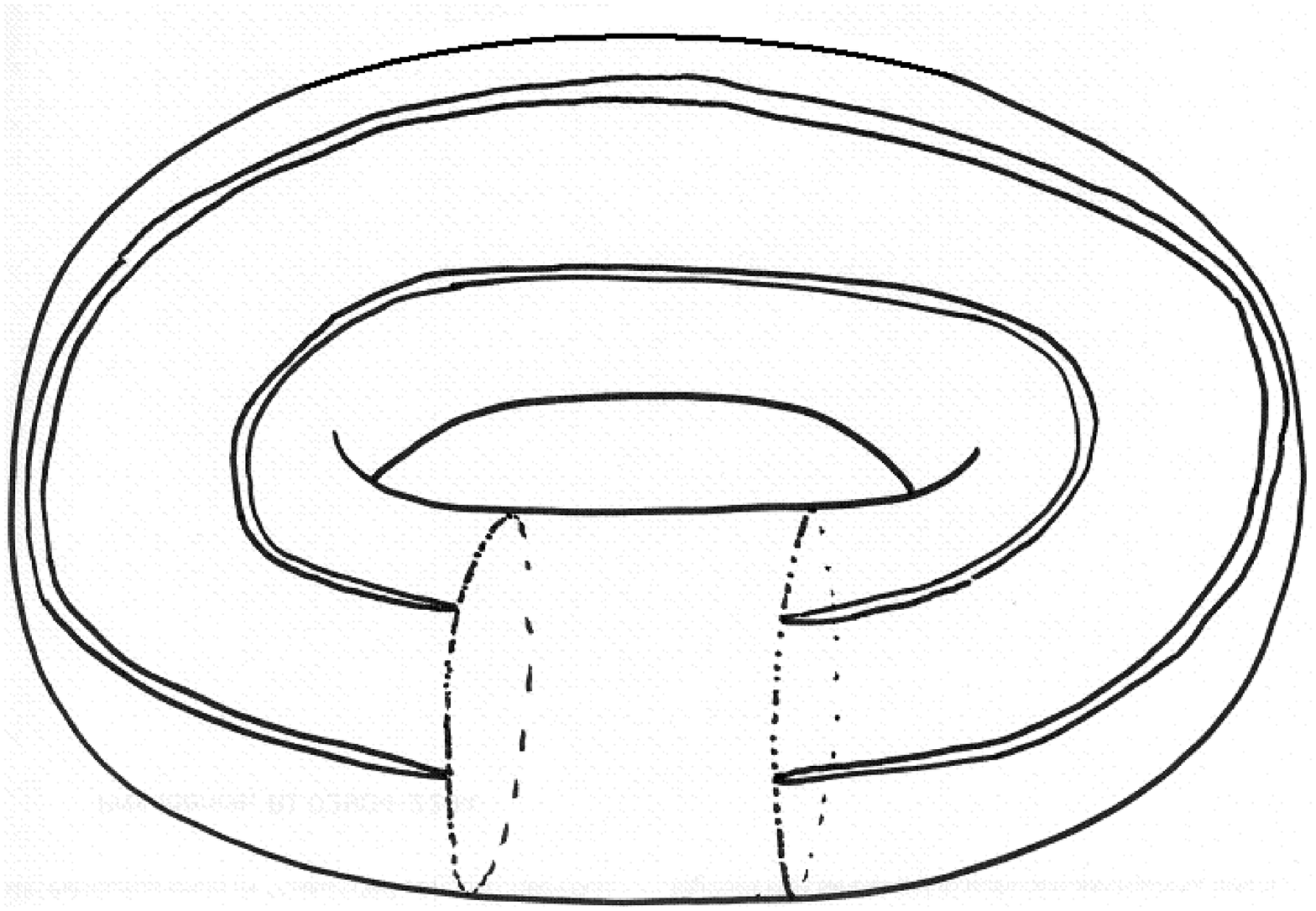}}    
\caption{The open manifold $U$.}
\end{figure}

By construction, our solenoid $S\subset U$.

Suppose that $\dim M\geq 3$. Then we can $C^\infty$-smoothly embed
$F:U\to M$ as follows: suppose that all cycles $C_i$ share a common
base-point $p_0\in M$ (and are otherwise disjoint to each other).
Then embed the central part $(0,c)\x \TT\subset U$ in a small ball
$B$ around $p_0$ and embed each of the $[c,1]\x U_i$ in $M-B$ in
such a way that if we contract $B$ to $p_0$ then the images of
$[c,1]\x \{t\}$, $t\in U_i$, represent cycles homologous to $C_i$.

The embedding of $f$ of $S$ into $M$ is defined as the composition
$S\inc U \stackrel{F}{\to} M$. This is an embedding according to
definition \ref{def:solenoid-in-manifold}. By theorem
\ref{thm:11.12}, as $S$ is uniquely ergodic, to prove that $(f,S)$
fully represents $a$, it is enough to see that $[f,S_{\mu}]=a$.

Let $\alpha$ be any closed $1$-form on $M$. Since $H^1(M)=H^1(M,B)$,
we may assume that $\alpha$ vanishes on $B$. We cover the solenoid
$S$ by the flow-boxes $((0,c)\x \TT)\cap S$ and $[c,1]\x K_i$,
$i=1,\ldots, r$. As $f^*\alpha$ vanishes in the first flow-box, we
have
  \begin{align*}
  \la [f,S_\mu],[\alpha]\ra &= \sum_{i=1}^r \int_{K_i} \left( \int_{[c,1]}
  f^*\alpha \right) d\mu_{K_i}(y) = \sum_{i=1}^r \int_{K_i} \la
  C_i ,[\alpha] \ra d\mu_{K_i}(y) \\ &= \sum_{i=1}^r \la
  C_i ,[\alpha] \ra \mu (K_i) = \sum_{i=1}^r \lambda_i \la
  C_i ,[\alpha] \ra =\la a,[\alpha]\ra\, ,
  \end{align*}
proving that $[f,S_\mu]=a$. Now the result follows from theorem
\ref{thm:11.12}.

\medskip

Now suppose that $\dim M=2$. Let us do the appropriate modifications
to the previous construction. Choose cycles $C_i$ sharing a common
base-point $p_0\in M$, and such that their intersections (and
self-intersections) away from $p_0$ are transversal. Changing $C_i$
by $2C_i$ if necessary, we suppose that going around $C_i$ does not
change the orientation (that is, the normal bundle to $C_i$ is
oriented). From the manifold $U$ in Figure 2, remove $[0,c]\x
\{\t_1\}$ to get the open $2$-manifold
 $$
 V=\big( (0,c)\x (0,1) \big) \bigcup \cup_i \big( [c,1]\x U_i \big) \, .
 $$

\begin{figure}[h] \label{figure3}
\centering
\resizebox{5cm}{!}{\includegraphics{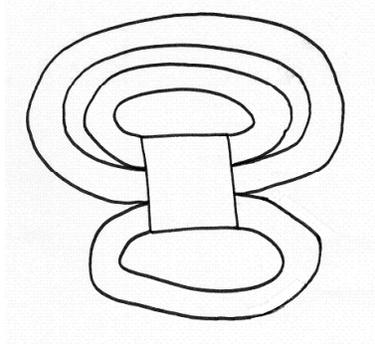}}    
\caption{The open manifold $V$}
\end{figure}

The manifold $V$ can be immersed into the surface $M$, $F:V\to M$,
in such a way that $(0,c)\x (0,1)$ is sent to a ball $B$ around
$p_0$, $[c,1]\x U_i$ are sent to $M-B$, the images of $[c,1]\x
\{t\}$, $t\in U_i$, represent cycles homologous to $C_i$ if we
contract $B$ to a point, and the intersections and
self-intersections of horizontal leaves are always transverse.

Note that the solenoid $S$ is not contained in $V$, since we have
removed $[0,c]\x\t_1$ from $U$. So we cannot define an immersion
$f:S\to M$ by restricting that of $F$. To define $f$ in $S\cap
((0,c)\x \TT)$, we need to explicit out our isotopy $h_t$.
Consider $h:\TT\to\TT$ and lift it to $\tilde{h}:\RR \to \RR$ with
$r:=\tilde{h}(0)\in (0,1)$. Consider a smooth function
$\rho:\RR\to [0,1]$, with $\rho(t)=1$ for $t\leq 0$, $\rho(t)=0$
for $t\geq c$, and $\rho'(t)<0$ for $t\in (0,c)$. Then we can
define
 $$
 h_t(x)= \tilde{h} ( \tilde{h}^{-1}(x) \rho(t) + x (1-\rho(t))) \mod \ZZ\, .
 $$

Define the immersion $f:S\to M$ as follows: $f$ equals $F$ for
$(t,x)\in [c,1]\x K \subset V$. For $(t,h^{-1}(h_t(x)))\in S\cap
([0,c]\x \TT)$, we set
  $$
  f(t,h^{-1}(h_t(x)))=\left\{\begin{array}{ll}
   F(t, (\tilde{h}^{-1}(x)+1) \rho(t) + x (1-\rho(t))), \qquad &
   x \in K \cap (0,r)\, ,\\
  F(t, \tilde{h}^{-1}(x) \rho(t) + x (1-\rho(t))), \qquad &
  x \in K \cap (r,1)\, .
   \end{array}\right.
  $$
It is easily checked that $f$ sends $S\cap ([0,c]\x \TT)$ into the
ball $B$ and the intersections of the leaves in this portion of the
solenoid are transverse.

The proof that $(f,S)$ fully represents a positive multiple of $a$
and that the generalized current $[f,S_\mu]=a$ goes as before.
\end{proof}

\begin{remark}\label{rem:no-precisa-compactness}
We  do not need $M$ compact for the above construction to work. If
$M$ is non-compact, take integer $1$-cycles $C_1,C_2,\ldots$
(possibly infinitely many) which form a basis of $H_1(M,\RR)$. Then
for any $a\in H_1(M,\RR)$ there exist an integer $r\geq 1$ and
$\lambda_1,\ldots,\lambda_r\in \RR$ with $a=\sum \lambda_i C_i$. The
construction of theorem \ref{thm:construct-1-solenoid} works.

The solenoid $S$ is oriented, regardless of $M$ being oriented or
not.
\end{remark}

\section{Realization of $H_k(M,\RR)$}\label{sec:k-solenoids}

Let $M$ be a smooth compact oriented Riemannian $C^\infty$ manifold
and let $a\in H_k(M,\RR)$ be a non-zero real $k$-homology class. We
are going to generalize the construction of section
\ref{sec:1-solenoids} to obtain uniquely ergodic $k$-solenoids
$(f,S)$ with a $1$-dimensional transversal structure, immersed in
$M$ and fully representing a positive multiple of $a$.

By theorem \ref{thm:11.12bis}, it is enough to produce a uniquely
ergodic oriented $k$-solenoid $S$, with a trapping region $W\subset
S$ and an immersion $f:S\to M$ such that $f(W)$ is sent to a
contractible ball in $M$, and with associated generalized current
$[f,S_\mu]=a$. Then, for the normalized (uniquely ergodic) measure
$\nu=\mu/\mu(S)$, it is $[f,S_\nu]=a/\mu(S)$, and so $(f,S)$ fully
represents  (by theorem \ref{thm:11.12bis}) the class $a/\mu(S)$.

Note also that if we change the Riemannian metric of $M$, then the
solenoid $(f,S)$ will still fully represent a positive multiple of
$a$ (although the scalar factor may change).

\medskip

To start with, fix a collection of compact $k$-dimensional smooth
oriented manifolds $S_1,\ldots, S_r$ and positive numbers
$\lambda_1,\ldots, \lambda_r>0$ such that $\sum \lambda_i=1$. Let
$h:\TT\to \TT$ be a diffeomorphism of the circle which is a Denjoy
counter-example with an irrational rotation number and of class
$C^{2-\epsilon}$, for some $\epsilon>0$. Hence $h$ is uniquely
ergodic. Let $\mu_K$ be the unique invariant probability measure,
which is supported in a Cantor set $K\subset \TT$. Partition the
Cantor set $K$ into $r$ disjoint compact subsets $K_1,\ldots, K_r$
in cyclic order, each of which with $\mu_K(K_i)=\lambda_i$.

We fix two points on each manifold $S_i$, and remove two small
balls, $D_i^+$ and $D_i^-$, around them. Denote
 $$
 S_i'=S_i- (D_i^+ \cup D_i^-)\, ,
 $$
so that $S_i'$ is a manifold with oriented boundary $\bd S_i'=\bd
D_i^+\sqcup \bd D_i^-$. Fix two diffeomorphisms: $\bd D_i^+ \cong
S^{k-1}$, orientation preserving, and $\bd D_i^- \cong S^{k-1}$,
orientation reversing. There are inclusions
 $$
  A_\pm :=\bigsqcup (\bd D_i^\pm \x K_i)
  \stackrel{i_\pm}{\inc} S^{k-1}\x S^1\, ,
 $$
with image $S^{k-1}\x K\subset S^{k-1}\x S^1$. Define
  $$
   S= \bigsqcup (S_i' \times K_i)_{/ x \sim i_+^{-1}(\id\x h)
   i_-(x),\,
   x\in A_- } \, .
  $$

This is an oriented $k$-solenoid of class $C^{\infty,2-\epsilon}$,
with $1$-dimensional transversal dimension. As $S^{k-1}\x K\subset
S$ in an obvious way, fixing a point $p\in S^{k-1}$, we have a
global transversal $T=\{p\}\x K \subset S^{k-1}\x K\subset S$.
Identifying $T\cong K$, the holonomy pseudo-group is generated by
$h:K\to K$. Hence $S$ is uniquely ergodic. Let $\mu$ denote the
tranversal measure corresponding to $\mu_K$.

We want to give an alternative description of $S$. Fix an isotopy
$h_t$, $t\in [0,1]$, from $\id$ to $h$. Define
 $$
  W' : = \{ (t,x,h^{-1}(h_t(y)))\, ;\, t\in [0,1], x\in S^{k-1}, y\in K
   \}\subset [0,1]\x S^{k-1} \times S^1\, .
 $$
Then
  $$
   S= \left( \bigsqcup (S_i' \times K_i) \sqcup W'
   \right)_{/ x \sim (0,i_-(x)),\, x\in \bd D_i^- \x K_i
  \atop x\sim (1,i_+(x)), \, x\in \bd D_i^+ \x K_i} \, .
  $$
Strictly speaking, we should say that they are diffeomorphic, but we
shall fix an identification. We define a map $\pi:S\to \TT$ by
  $$
  \begin{aligned}
  & \pi(t,x,h^{-1}(h_t(y))) =t-\frac12, \qquad  (t,x,h^{-1}(h_t(y)))\in W' \\
  & \pi(p)= \frac12, \qquad p\in S-W'
  \end{aligned}
  $$
Then $W=\Int(W')=\pi^{-1}(-\frac12,\frac12)$ is a trapping region as
in definition \ref{def:trapping}.

\medskip

Consider angles $\t_1,\t_2,\ldots, \t_n\in \TT$ in the same cyclic
order as the $K_i$, such that $K_i$ is contained in the open subset
$U_i\subset T$ with boundary points $\t_i$ and $\t_{i+1}$ (denoting
$\t_{n+1}=\t_1$). We may assume that $\t_1=0$. Then solenoid $S$
sits inside the $(k+1)$-dimensional open manifold
  $$
  X= \bigsqcup (S_i'\times U_i) \sqcup ([0,1]\x S^{k-1} \times S^1
  )_{/ x \sim (0,i_-(x)),\, x\in \bd D_i^- \x U_i
  \atop x\sim (1,i_+(x)), \, x\in \bd D_i^+ \x U_i} \, ,
  $$
as the collection of points $(x,y)$, $x\in S_i'$, $y\in K_i$,
together with the points $(t,x,h^{-1}(h_t(y)))$, $x\in S^{k-1}$,
$y\in K$, $t\in [0,1]$.

\begin{figure}[h]\label{figure4}
\centering
\resizebox{10cm}{!}{\includegraphics{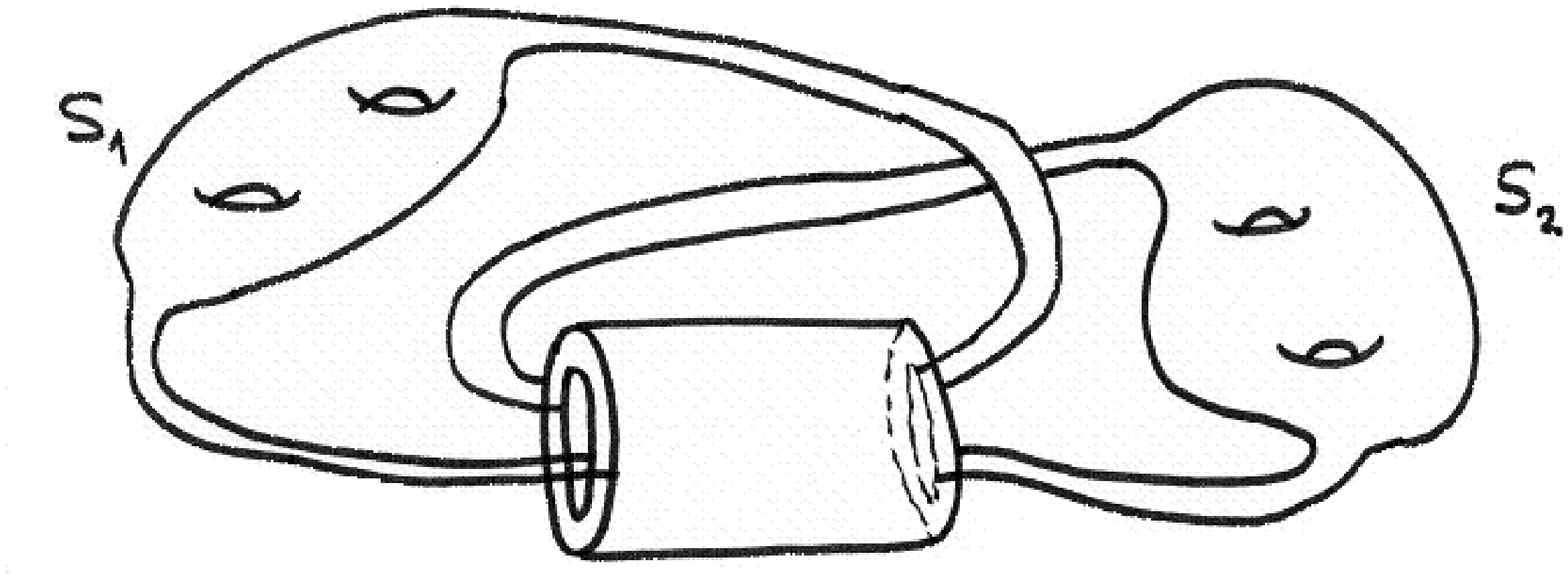}}    
\caption{The manifold $X$.}
\end{figure}

\begin{remark}
The $1$-solenoid constructed in section \ref{sec:1-solenoids}
corresponds to the case $S_i=S^1$, $i=1,\ldots, r$.
\end{remark}

\begin{theorem}\label{thm:construct-k-solenoid}
Let $M$ be a compact oriented smooth Riemannian manifold of
dimension $n$, and let $a\in H_k(M,{\RR})$ be a non-zero real
$k$-homology class.
 \begin{enumerate}
 \item If $n\geq 2k+1$ then (a positive multiple of) $a$ can be fully represented by an
embedding of a uniquely ergodic oriented $k$-solenoid $S$ into $M$.
 \item If $n-k$ is odd then (a positive multiple of) $a$ can be fully represented by a
transversally immersion of a uniquely ergodic oriented $k$-solenoid
$S$ into $M$.
 \item In other cases (a positive multiple of) $a$ can be fully represented by an immersed
uniquely ergodic oriented $k$-solenoid $S$ in $M$.
\end{enumerate}
\end{theorem}

\begin{proof}
By a theorem of Thom \cite{Thom}, if $a\in H_k(M,\ZZ)$ then there
exists $N>>0$ such that $N\, a$ is represented by a smooth
submanifold of $M$. This submanifold is oriented because it
represents a non-zero homology class. Moreover, if $n\geq 2k+1$ or
$n-k$ is odd then it can be arranged that the normal bundle to the
submanifold is trivial \cite{Thom}.

Take a collection $C_1,\ldots, C_{b_k}\in H_k(M,{\ZZ})$ which are a
basis of $H_k(M,{\QQ})$ and such that $C_i$ is represented by a
smooth submanifold $S_i\subset M$ (with trivial normal bundle if
$n\geq 2k+1$ or $n-k$ is odd). They can be assumed to be in general
position. After switching the orientations of $C_i$ if necessary,
reordering the cycles and multiplying $a$ by a suitable positive
real number, we may suppose that
 $$
 a=\lambda_1 C_1 +\ldots +\lambda_r C_r,
 $$
for some $r\geq 1$, $\lambda_i>0$, $1\leq i\leq r$, and $\sum
\lambda_i=1$. We construct the solenoid $S$ with the procedure above
starting with the manifolds $S_i$ and coefficients $\lambda_i$. This
is a uniquely ergodic $k$-solenoid with a $1$-dimensional
transversal structure, and a trapping region $W\subset S$.

Now we want to define an immersion $f:S\to M$, and to prove that
$(f,S)$ fully represents $a$. We have the following cases:

\begin{enumerate}
\item $n\geq 2k+1$. The general position property on the $S_i$
implies that all $S_i$ are disjoint submanifolds of $M$. As the
normal bundle to $S_i$ is trivial and $U_i$ is an interval, we can
embedded $S_i\x U_i$ in a small neighbourhood of $S_i$.

Fix a base point $p_0\in M$ off all $S_i$. Take a small box
$B\subset M$ around $p_0$ of the form $B=[0,1]\x D^{n-1}$, where
$D^{n-1}$ is the open $(n-1)$-dimensional ball. Consider a circle
$S^1\subset D^{k+1}\subset D^{n-1}$ and let $D^k\x S^1 \subset
D^{k+1}\subset D^{n-1}$ be a tubular neighbourhood of it, with
boundary $S^{k-1}\x S^1$.

For each $i=1,\ldots, r$, fix $x_i\in U_i$, and consider two paths
in $M-\Int(B)$, $\g_i^\pm$, where $\g_i^-$ goes from the point
$(0,x_i)\in \{0\}\x U_i \subset \{0\}\x S^1\subset \{0\}\x
D^{n-1}\subset B$ to the point $(p_i^-,x_i) \in S_i\x U_i$, and
$\g_i^+$ goes from $(1,x_i)\in \{1\}\x U_i \subset \{1\}\x
S^1\subset \{1\}\x D^{n-1}\subset B$ to $(p_i^+,x_i) \in S_i\x U_i$.
We arrange that $\g_i^\pm$ are transverse to $S_i\x U_i$ at
$(p_i^\pm,x_i)$ and are disjoint from all $S_j$ otherwise.

We thicken $\g_i^\pm$ to immersions $\g_i^\pm \x D^{k}\x U_i$ into
$M-\Int(B)$ such that one extreme goes to $D_i^\pm\x U_i$ and the
other goes to either $D^{k} \x U_i \x\{0\} \subset D^{k}\x S^1
\x\{0\} \subset D^{n-1}\x\{0\} \subset B$ for $\g_i^-$, or $D^{k} \x
U_i \x\{1\} \subset D^{k}\x S^1 \x\{1\} \subset D^{n-1}\x\{1\}
\subset B$ for $\g_i^-$. It is possible to do this in such a way
that the $U_i$ directions match, since $n\geq k+2$.

Recall that $S_i'=S_i - (D_i^+\cup D_i^-)$, and set
 \begin{equation*} 
  S_i''= S_i' \cup \g_i^\pm \x S^{k-1} \,,
 \end{equation*}
which is diffeomorphic to $S_i'$ (to be rigorous, we should smooth
out corners). Then
 $$
 U:= \bigcup ((S_i'\x U_i) \cup ( \g_i^\pm \x S^{k-1} \x U_i) \cup
 ([0,1]\x S^{k-1} \x S^1)
 $$
is a $(k+1)$-dimensional open manifold embedded in $M$. The manifold
$U$ is foliated as follows: $S_i''\x U_i$ is foliated by
$S_i''\x\{y\}$, for $y\in U_i$, and $[0,1]\x S^{k-1} \x S^1$ is
foliated by
 $$
 L_y=\{(t, x, h^{-1}(h_t(y))) \ ; \ t\in [0,1], x\in S^{k-1} \} \, ,
 $$
for $y \in S^1$. Clearly the solenoid $S$ is a subsolenoid of $U$,
$S\subset U$. Restricting the embedding $F:U\to M$ to $S$ we get an
embedding $f:S\to M$.

By construction $f(W)\subset \Int(B)$, i.e. the image of the
trapping region is contained in a contractible ball.

 \item $n-k>1$ and odd.
 As the codimension is odd, the general position property
 implies that the submanifolds $S_i$ and $S_j$,
 $i\neq j$, intersect transversally. This is argued as follows: represent $S_i$ as
 the preimage of a regular value $v_i$ of a smooth application $f_i:M\to S^{n-k}$.
 Then consider $f_i\x f_j: M\to S^{n-k}\x S^{n-k}$ and take a regular value
 $(v_i',v_j')$ of this map near $(v_i,v_j)$.

 The same construction as in (1) works now,
 with the modification that we have to allow
 intersections of different leaves, but we may take them to be always transversal.
 So we get a transversal immersion $f:S\to M$.

 \item $n-k$ even. The embedded submanifolds $S_i\subset M$ may have
 a non-trivial normal bundle. Moreover, they may intersect each other.
 The intersections are generically transversal along $S_i\cap S_j$, but not
 everywhere transversal.

 Take a generic section $s_i$ of the normal bundle to $S_i$, so that the image of
 the section intersects transversally the zero section.
 We have maps
 $f_i:S_i\x (-\epsilon,\epsilon) \to M$ defined by $f_i(x,t)=\exp_x (t s_i(x))$,
 where  $\exp$ is the exponential mapping (with respect to a Riemannian metric
 on $M$), which is an immersion for each
 $t$ fixed, for $\epsilon>0$ small enough. With these maps, the
 previous construction gives a solenoid immersion $f:S\to M$.

 \item $n-k=1$. The submanifolds $S_i$ have trivial normal bundle and
 they intersect each other transversally. We cannot avoid
 that the paths $\g_i^\pm$ intersect other $S_j$, but we arrange these
 intersections to be transverse. This produces a transversal
 immersion $f$ of the region $S-W$ of the solenoid into $M-\Int(B)$.

 We have to modify the previous construction of the immersion of $W$ into $B$,
 as the codimension one does not leave enough room for it to work. Consider the box $B=[0,1]\x
 D^{n-1}$ and remove the axis $A=[0,1]\x \{0\}$. Use polar
 coordinates to identify $B-A= [0,1] \x S^{k-1}\x (0,1)$, where the
 third coordinate corresponds to the radius.
 By construction, $W'\subset S$ embeds into $C=[0,l]\x S^{k-1} \x
 S^1$, as the set of points $(t,x,h^{-1}(h_t(y)))$, $t\in [0,1]$, $x\in S^{k-1}$
 and $y\in K$. We remove $D=[0,1]\x S^{k-1}\x \tau_1$ from $C$, so that $C-D=[0,1] \x S^{k-1}\x
(0,1)$. Then $W'$ immerses into $C-D$, by using the process at the
end of the proof of theorem \ref{thm:construct-1-solenoid} (now
there is an extra factor $S^{k-1}$ which plays no role). This is a
transversal immersion.

There is one extra detail that we should be careful about. When
connecting $p_i^\pm$ with the two faces of $B$, the orientations of
the $U_i$ should match. This happens because the normal bundle to
$S_i$ is trivial, and in this case $S_i\x U_i$ is (diffeomorphic to)
the normal bundle to $S_i$.
\end{enumerate}

To prove that $(f,S)$ fully represents (a positive multiple of) $a$,
we use theorem \ref{thm:11.12bis}. The solenoid $S$ has a trapping
region $W$, and $f(W)\subset \Int(B)$, a contractible ball in $M$.
So we only need to see that $[f,S_\mu]=a$. Recall that the
associated transversal measure is $\mu_K$ on the transversal $K$.
Let $\alpha$ be any closed $1$-form on $M$. Since $H^1(M)=H^1(M,B)$,
we may assume that $\alpha$ vanishes on $B$. We cover the solenoid
$S$ by the flow-boxes $S_i''\x K_i$, $i=1,\ldots, r$, and $W'$
(where the form $\alpha$ vanishes). Thus
  \begin{align*}
  \la [f,S_\mu],[\alpha]\ra &= \sum_{i=1}^r \int_{K_i} \left( \int_{S_i''}
  f^*\alpha \right) d\mu_{K_i}(y) = \sum_{i=1}^r \int_{K_i} \la
  C_i ,[\alpha] \ra d\mu_{K_i}(y) \\ &= \sum_{i=1}^r \la
  C_i ,[\alpha] \ra \mu (K_i) = \sum_{i=1}^r \lambda_i \la
  C_i ,[\alpha] \ra =\la a,[\alpha]\ra\, ,
  \end{align*}
proving that $[f,S_\mu]=a$.
\end{proof}

\begin{remark}
A similar comment to that of remark \ref{rem:no-precisa-compactness}
applies to the present situation, that is, the compactness of $M$ is
not necessary.

If $M$ is non-orientable, we may consider its oriented double cover
$\pi:\tilde{M}\to M$. For any non-zero $a\in H_k(M,\RR)$, there
exists $\tilde a \in H_k(\tilde M, \RR)$ with $\pi_*(\tilde a)=a$.
We may construct an oriented uniquely-ergodic $k$-solenoid $(f,S)$
immersed in $\tilde M$ fully representing $\lambda\, \tilde a$, for
some $\lambda>0$. Then $(\pi\circ f,S)$ is immersed in $M$ and fully
represents $\lambda\,a$.

\end{remark}

\section{Homotopy of solenoids}\label{sec:homotopy}

\begin{definition}\textbf{\em (Solenoid with boundary)}
  Let $0\leq r,s\leq \omega$, and let $k,l\geq 0$ be two integers.
  A foliated manifold with boundary (of dimension $k+l$, with $k$-dimensional
  leaves, of class $C^{r,s}$) is a smooth
  manifold $W$ with boundary, of dimension $k+l$, endowed with an
  atlas
  $\{ (U_i,\varphi_i)\}$ of charts
  $$
  \varphi_i: U_i \to \varphi_i(U_i) \subset \RR^{k+l}_+ =\{(x_1,\ldots,
  x_k, y_1,\ldots, y_l \ ; \ x_1\geq 0\} \, ,
  $$
whose changes of charts are of the form $\varphi_i\circ
\varphi_j^{-1}(x,y) = (X_{ij}(x,y), Y_{ij}(y))$, where $Y_{ij}(y)$
is of class $C^s$ and $X_{ij}(x,y)$ is of class $C^{r,s}$.

A pre-solenoid with boundary is a pair $(S,W)$ where $W$ is a
foliated manifold with boundary and $S\subset W$ is a compact
subspace which is a collection of leaves.

Two pre-solenoids with boundary $(S,W_1)$ and $(S,W_2)$ are
equivalent if there are open subsets $U_1\subset W_1$, $U_2\subset
W_2$ with $S\subset U_1$ and $S\subset U_2$, and a diffeomorphism
$f: U_1\to U_2$ (preserving leaves, of class $C^{r,s}$) which is the
identity on $S$.

A $k$-solenoid with boundary is an equivalence class of
pre-solenoids with boundary.
\end{definition}

Note that any manifold with boundary is a solenoid with boundary.

The boundary of a $k$-solenoid with boundary $S$ is the
$(k-1)$-solenoid (without boundary) $\bd S$ defined by the foliated
manifold $\bd W$, where $W$ is a foliated manifold with boundary
defining the solenoid structure of $S$.

A $k$-solenoid with boundary $S$ has two types of flow-boxes. If
$p\in S-\bd S$ is an interior point, then there is a flow-box
$(U,\varphi)$ with $p\in U$, of the form $\varphi : U\to D^k\times
K(U)$. If $p\in \bd S$ is a boundary point, then there is a flow-box
$(U,\varphi)$ with $p\in U$ such that $\varphi$ is a homeomorphism
 $$
 \varphi : U\to D^k_+\times K(U) \, ,
 $$
where $D^k_+=\{ (x_1,\ldots, x_k)\in D_k \, ; \, x_1\geq 0\}$, and
$K(U)\subset \RR^l$, $\varphi(p)=(0,\ldots, 0, y_0)$, for some
$y_0\in K$. Note that writing
 $$
 U'=\bd S\cap U =\varphi^{-1}(D^{k-1}\x K(U))\, ,
 $$
where $D^{k-1} = \{ (0,x_2,\ldots, x_k)\in D_k \} \subset D_k^+$,
$(U',\varphi_{|U'})$ is a flow-box for $\bd S$. Therefore, if $T$ is
a transversal for $\bd S$, then it is also transversal for $S$.

For a solenoid with boundary $S$ there is also a well-defined notion
of holonomy pseudo-group. If $T$ is a local transversal for $\bd S$,
and $h:T\to T$ is a holonomy map for $\bd S$ defined by a path in
$\bd S$, then $h$ lies in the holonomy pseudo-group of $S$. So
  $$
  \Hol_{\bd S}(T)\subset \Hol_S(T)\, ,
  $$
but they are in general not equal. In particular, if $S$ is
connected with non-empty boundary then
 $$
 \cM_\cT(S) \subset \cM_\cT(\bd S) \, .
 $$
That is, if $\mu=(\mu_T)$ is a transversal measure for $S$, then it
yields a transversal measure for $\bd S$, by considering only those
transversals $T$ which are transversals for $\bd S$. We denote this
transversal measure by $\mu$ again.

If $S$ comes equipped with an orientation, then $\bd S$ has a
natural induced orientation.
Note that any leaf $l\subset S$ is a manifold with boundary and
each connected component of $\bd l$ is a leaf of $\bd S$.

\begin{theorem}\textbf{\em (Stokes theorem)} \label{thm:Stokes}
  Let $(f,S_\mu)$ be an oriented $(k+1)$-solenoid with boundary,
  endowed with a transversal measure, and immersed into a smooth
  manifold $M$. Let $\omega$ be a $k$-form on $M$. Then
  $$
  \la [f,S_\mu], d\omega \ra = \la [f_{|\bd S},\bd S_\mu] , \omega
  \ra \, .
  $$
\end{theorem}

\begin{proof}
  Let $\{U_i\}$ be a covering of $S$ by flow-boxes, and let $\{\rho_i\}$
  be a partition of unity subordinated to it. Adding up the
  equalities
  $$
  \begin{aligned}
    \int_{K(U_i)}& \left( \int_{L_y} d\rho_i \wedge f^* \omega \right)
    d\mu_{K(U_i)} (y)
  +
    \int_{K(U_i)} \left( \int_{L_y} \rho_i f^* d\omega \right)
    d\mu_{K(U_i)} (y) \\
  &=
    \int_{K(U_i)} \left( \int_{L_y} d (\rho_i f^* \omega) \right)
    d\mu_{K(U_i)} (y)
  =
    \int_{K(U_i)} \left( \int_{\bd L_y} \rho_i f^* \omega \right)
    d\mu_{K(U_i)} (y)\,,
    \end{aligned}
  $$
 for all $i$, and using that $\sum d\rho_i \equiv 0$, we get
  $$
  \begin{aligned}
  \la  [f,S_\mu], d\omega \ra &\, =
   \sum_i \int_{K(U_i)} \left( \int_{L_y} \rho_i f^* d\omega \right)
    d\mu_{K(U_i)} (y) \\
  & \, =
   \sum_i \int_{K(U_i)} \left( \int_{\bd L_y} \rho_i f^* \omega \right)
    d\mu_{K(U_i)} (y)
    =
    \la [f_{|\bd S},\bd S_\mu] , \omega  \ra \, .
  \end{aligned}
  $$
\end{proof}

Let $S$ be a $k$-solenoid of class $C^{r,s}$. We give $S\x I=S\x
[0,1]$ a natural $(k+1)$-solenoid with boundary
structure of the same class, by taking
a foliated manifold $W\supset S$ defining the solenoid structure of
$S$, and foliating $W\x I$ with the leaves $l\x I$,
$l\subset W$ being a leaf of $W$. Then $S\x I \subset W\x
I$ is a $(k+1)$-solenoid with boundary.
The boundary of $S\x I$ is
 $$
 (S\x\{0\}) \sqcup (S\x \{1\})\, .
 $$

If $S$ is oriented then $S\x I$ is naturally oriented and its boundary
consists of  $S\x\{0\} \cong S$ with orientation reversed, and
$S\x\{1\} \cong S$ with orientation preserved.

Moreover if $T$ is a transversal for $S$, then $T'=T \x \{0\}$
is a transversal for $S'=S\x I$.
The following is immediate.

\begin{lemma} \label{lem:SxI}
There is an identification of the holonomies of $S$ and $S\x I$.
More precisely, under the identification $T\cong T' =T \x \{0\}$,
 $$
 \Hol_S(T) =\Hol_{S\x I}(T')\, .
 $$

In particular,
 $$
 \cM_\cT(S)=\cM_\cT(S\x I)\, .
 $$
\end{lemma}

\begin{definition}\textbf{\em (Equivalence of immersions)}\label{def:14.equivalence}
Two solenoid immersions $(f_0,S_0)$ and $(f_1,S_1)$ of class
$C^{r,s}$ in $M$ are immersed equivalent if there is a
$C^{r,s}$-diffeomorphism $h:S_0\to S_1$ such that
 $$
 f_0=f_1\circ h \, .
 $$
Two measured solenoid immersions are immersed equivalent if $h$ can be chosen
to preserve the transversal measures.
\end{definition}

\begin{definition}\textbf{\em (Homotopy of immersions)}\label{def:14.homotopy}
Let $S$ be a $k$-solenoid of class $C^{r,s}$ with $r\geq 1$. A
homotopy between immersions $f_0:S\to M$ and $f_1: S\to M$ is an
immersion of solenoids $f:S\x I\to M$ such that $f_0(x)=f(x,0)$ and
$f_1(x)=f(x,1)$.
\end{definition}

\begin{definition} \textbf{\em (Cobordism of solenoids)}
\label{def:14.cobordism} Let $S_0$ and $S_1$ be two
$C^{r,s}$-solenoids. A cobordism of solenoids $S$ is a
$(k+1)$-solenoid with boundary $\bd S=S_0\sqcup S_1$.

If $S_0$ and $S_1$ are oriented, then an oriented cobordism is a
cobordism $S$ which is an oriented solenoid such that it induces the
given orientation on $S_1$ and the reversed orientation on $S_0$.

If $S_0$ and $S_1$ have transversal measures $\mu_0$ and $\mu_1$,
respectively, then a measured cobordism is a cobordism $S$ endowed
with a transversal measure $\mu$ inducing the measures $\mu_0$ and
$\mu_1$ on $S_0$ and $S_1$, respectively.
\end{definition}

\begin{definition}\textbf{\em (Homology equivalence)}
\label{def:homology}
 Let $(f_0,S_0)$ and $(f_1,S_1)$ be two immersed solenoids in $M$.
 We say that they are homology equivalent if there exists a
 cobordism of solenoids $S$ between $S_0$ and $S_1$ and a solenoid
 immersion $f:S \to M$ with $f_{|S_0}=f_0$, $f_{|S_1}=f_1$.
 We call $(f,S)$ a homology between $(f_0,S_0)$ and $(f_1,S_1)$.

 Let $(f_0,S_{0,\mu_0})$ and $(f_1,S_{1,\mu_1})$ be two immersed
 oriented measured solenoids. They are homology equivalent if
 there exists a immersed oriented measured solenoid $(f,S_\mu)$ such that
 $(f,S)$ is a homology between $(f_0,S_0)$ and $(f_1,S_1)$ and
 $S_\mu$ is a measured oriented cobordism from $S_0$ to $S_1$.
\end{definition}

Clearly two homotopic immersions of a solenoid give homology
equivalent immersions.

\begin{theorem} \label{thm:homeo-RS}
 Suppose that two oriented measured solenoids
 $(f_0,S_{0,\mu_0})$ and $(f_1,S_{1,\mu_1})$ immersed in $M$
 are homology equivalent. Then the
generalized currents coincide
   $$
   [f_0,S_{0,\mu_0}]=[f_1,S_{1,\mu_1}]\, .
   $$

 The same happens if they are immersed equivalent.
\end{theorem}

\begin{proof}
In the first case, let $\omega$ be a closed $k$-form on $M$, then
Stokes' theorem gives
  $$
  \la [f_1,S_{1,\mu_1}] , \omega \ra- \la [f_0,S_{0,\mu_0}] ,
  \omega \ra =
  \la [f_{|\bd S},\bd S_\mu] , \omega \ra = \la [f,S_\mu],d\omega\ra =0\, .
  $$

In the second case, $f_0=f_1\circ h$ implies that the actions of the
generalized currents over a closed form on $M$ coincide, since the
pull-back of the form to the solenoids agree through the
diffeomorphism $h$, and the integrals over the transversal measure
gives the same numbers, since the measures correspond by $h$.
\end{proof}

\begin{remark} \label{rem:que-numero-le-pongo}
 In both definitions \ref{def:14.cobordism} and
 \ref{def:14.homotopy}, we do not need to require that $f$ be an
 immersion. Actually, the generalized current $[f,S_\mu]$ makes
 sense for any measured solenoid $S_\mu$ and any regular map $f:S\to
 M$, of class $C^{r,s}$ with $r\geq 1$. Theorem \ref{thm:homeo-RS}
 still holds with these extended notions.
\end{remark}

\section{Intersection theory of solenoids}\label{sec:intersection}

Let $M$ be a smooth $C^\infty$ oriented manifold.

\begin{definition} \textbf{\em (Transverse intersection)}
\label{def:15.1} Let $(f_1,S_1)$, $(f_2,S_2)$ be two immersed
solenoids in $M$. We say that they intersect transversally if, for
every $p_1\in S_1$, $p_2\in S_2$ such that $f_1(p_1)=f_2(p_2)$, the
images of the leaves through $p_1$ and $p_2$ intersect
transversally.
\end{definition}

If two immersed solenoids $(f_1,S_1)$, $(f_2,S_2)$, of dimensions
$k_1$, $k_2$ respectively, intersect transversally. We define the
intersection of immersed solenoids $(f,S)$ as defined by
  \begin{equation} \label{eqn:intersect-solenoids}
   S=\{ (p_1,p_2) \in S_1\x S_2 \ ; \ f_1(p_1)=f_2(p_2)\} \, .
   \end{equation}
and by the map  $f:S\to M$ given by
  \begin{equation} \label{eqn:immersion-intersect}
 f(p_1,p_2) =f_1(p_1)=f_2(p_2) , \quad (p_1,p_2) \in S\, .
   \end{equation}
We will see that $S$, the intersection solenoid, is indeed a
solenoid. Also the intersection $(f,S)$ of the two immersed
solenoids $(f_1,S_1)$, $(f_2,S_2)$ is an immersed solenoid. In order
to prove this, we consider the intersection of the product solenoid
$(f_1\x f_2,S_1\x S_2)$ in $M\x M$ with the diagonal $\Delta \subset
M\x M$. So we have to analyze first the case of the intersection of
an immersed solenoid with a submanifold. The notion of transverse
intersection given in definition \ref{def:15.1} applies to this case
(a submanifold is an embedded solenoid).

\begin{lemma} \label{lem:intersect-submfd}
  Let $(f,S)$ be an immersed $k$-solenoid in $M$ intersecting
  transversally an embedded closed submanifold $N\subset M$ of
  codimension $q$. Suppose that $S'=f^{-1}(N)\subset S$ is
  non-empty, then $(f_{|S'},S')$ is an immersed $(k-q)$-solenoid in $N$.

  If $S$ and $N$ are oriented, so is $S'$.

  If $S$ has a transversal measure $\mu$, then $S'$ inherits
  a natural transversal measure, also denoted by $\mu$.
\end{lemma}

\begin{proof}
First of all, note that $S'$ is a compact and Hausdorff space.

Let $W$ be a foliated manifold defining the solenoid structure of
$S$ such that there is a smooth map $\hat f:W\to M$ of class
$C^{r,s}$, extending $f:S\to M$, which is an immersion on leaves. By
definition, for any leaf $l\subset S$, $f(l)$ is transverse to $N$.
Thus reducing $W$ if necessary, the same transversality property
occurs for any leaf of $W$. The transversality of the leaves implies
that the map $\hat f:W\to M$ is transversal to the submanifold
$N\subset M$, meaning that for any $p\in W$ such that $\hat f(p)\in
N$,
 $$
 d\hat f(p)(T_pW) + T_{\hat f(p)}N =T_{\hat f(p)}M \, .
 $$
This implies that $W'=\hat f^{-1}(N)$ is a submanifold of $W$ of
codimension $q$ (in particular, $k-q\geq 0$). Moreover, it is
foliated by the connected components $l'$ of $l\cap \hat
f^{-1}(N)=(\hat f_{|l})^{-1}(N)$, where $l$ are the leaves of $W$.
By transversality of $\hat f$ along the leaves, $l'$ is a
$(k-q)$-dimensional submanifold of $l$. So $W'$ is a foliated
manifold with leaves of dimension $k-q$. This gives the required
solenoid structure to $S'=S \cap \hat f^{-1}(N)=f^{-1}(N)$.

Clearly, $f_{|S'}:S'\to N$ is an immersion (of class $C^{r,s}$)
since $\hat f_{|W'}:W'\to N$ is a smooth map which is an immersion
on leaves.

If $S$ and $N$ are oriented, then each intersection $l'=l\cap
\hat f^{-1}(N)$ is also oriented (using that $M$ is oriented as well). Therefore the
leaves of $S'$ are oriented, and hence $S'$ is an oriented solenoid.

  Let $p\in S'$ and let $U\cong D^k\x K(U)$ be a flow-box for $S$ around $p$.
  We can take $U$ small enough so that $f(U)$ is contained in a
  chart of $M$ in which $N$ is defined by functions $x_1=\ldots =x_q
  =0$. By the transversality property, the differentials
  $dx_1,\ldots, dx_q$ are linearly independent on each leaf $f(D^k\x
  \{y\})$, $y\in K(U)$. Therefore, $x_1, \ldots, x_q$ can be completed to a set of
  functions $x_1,\ldots, x_k$ such that $dx_1,\ldots, dx_k$ are a
  basis of the cotangent space for each leaf (reducing $U$ if
  necessary). Thus the pull-back of $x=(x_1,\ldots,x_k)$ to $U$
  give coordinates functions so that, using the coordinate $y\in K(U)$
  for the transversal direction, $(x,y)$ are coordinates for $U$,
  and $f^{-1}(N)$ is defined as $x_1=\ldots =x_q
  =0$. This means that
    $$
    S'\cap U \cong \{ (0,\ldots, 0, x_{q+1},\ldots, x_k,y)\in D^k\x
    K(U)\} \cong D^{k-q} \x K(U)\, .
    $$
Therefore any local transversal $T$ for $S'$ is a local transversal
for $S$, and any holonomy map for $S'$ is a holonomy map for $S$. So
  $$
  \Hol_{S'}(T) \subset \Hol_{S}(T) \, .
  $$
Hence a transversal measure for $S$ gives a transversal measure
for $S'$.
\end{proof}

Now we can address the general case.

\begin{proposition} \label{prop:15.3}
Suppose that $(f_1,S_1)$, $(f_2,S_2)$ are two immersed solenoids in $M$
intersecting transversally, and let $S$ be its intersection solenoid
defined in (\ref{eqn:intersect-solenoids}) and let $f$ be the map
(\ref{eqn:immersion-intersect}). If $S\neq \emptyset$,
then $(f,S)$ is an immersed solenoid of dimension
$k=k_1+k_2-n$ (in particular, $k$ is a non-negative number).

If $S_1$ and $S_2$ are both oriented, then $S$ is also oriented.

If $S_1$ and $S_2$ are endowed with transversal measures $\mu_1$
and $\mu_2$ respectively,
then $S$ has an induced measure $\mu$.
\end{proposition}

\begin{proof}
  The product $S_1\x S_2$ is a
  $(k_1+k_2)$-solenoid and
   $$
   F=f_1\x f_2 :S_1\x S_2 \to M\x M
   $$
  is an immersion. Let $\Delta\subset M\x M$ be the diagonal. There is
  an identification (as sets)
   $$
   S = (S_1\x S_2) \cap F^{-1}(\Delta)\, .
   $$
 The condition that $(S_1,f_1)$, $(S_2,f_2)$ intersect transversally can
 be translated into that
 $(F,S_1\x S_2)$ and $\Delta$ intersect transversally in $M\x M$.

Therefore applying lemma \ref{lem:intersect-submfd},
$(S,F_{|S})$ is an immersed $k$-solenoid,
where $F_{|S}:S \to \Delta$ is defined as $F(x_1,x_2)=f_1(x_1)$.
Using the diffeomorphism $M\cong\Delta$, $x\mapsto (x,x)$,
$F_{|S}$  corresponds to $f:S\to M$. So $(S,f)$ is an immersed $k$-solenoid.

If $S_1$ and $S_2$ are both oriented, then $S_1\x S_2$ is also oriented. By
lemma \ref{lem:intersect-submfd}, $S$ inherits an orientation.

If $S_1$ and $S_2$ are endowed with transversal measures $\mu_1$
and $\mu_2$, then $S_1\x S_2$ has a product transversal measure $\mu$.
For any local transversals $T_1$ and $T_2$ to $S_1$ and $S_2$, respectively,
$T=T_1\x T_2$ is a local transversal to $S_1\x S_2$ (and conversely). We
define
 \begin{equation} \label{eqn:mu-product}
 \mu_T= \mu_{1,T_1} \x \mu_{2,T_2} \, .
 \end{equation}
Now lemma \ref{lem:intersect-submfd} applies to give the transversal measure for $S$.
Note that the local transversals to $S$ are of the form $T_1\x T_2$, for some
local transversals $T_1$ and $T_2$ to $S_1$ and $S_2$.
\end{proof}

\begin{remark} \label{rem:15.4}
If $k_1+k_2=n$ then $S$ is a $0$-solenoid.
For a $0$-solenoid $S$, an orientation is a continuous
assignment $\epsilon:S\to \{\pm 1\}$ of
sign to each point of $S$.

Note also that for a $0$-solenoid $S$, $T=S$ is a transversal and a
transversal measure is a Borel measure on $S$.
\end{remark}

Let $(f_1,S_1)$, $(f_2,S_2)$ be two immersed solenoids in $M$
intersecting transversally, with $(f,S)$ its intersection solenoid.
Let $p=(p_1,p_2)\in S$. Then we can choose flow-boxes $U_1=
D^{k_1}\x K(U_1)$ for $S_1$ around $p_1$ with coordinates
$(x_1,\ldots, x_{k_1}, y)$, and $U_2= D^{k_2}\x K(U_2)$ for $S_2$
around $p_2$ with coordinates $(x_1,\ldots, x_{k_2}, z)$, and
coordinates for $M$ around $f(p)$, such that
  $$
  \begin{aligned}
   f_1(x,y) & \, = (x_1,\ldots,x_{k_1+k_2-n},x_{k_1+k_2-n+1},\ldots,
   x_{k_1}, B_1(x,y),\ldots, B_{n-k_1}(x,y)) \, ,\\
   f_2(x,z) & \, = (x_1,\ldots, x_{k_1+k_2-n}, C_1(x,z),\ldots,
   C_{n-k_2}(x,z), x_{k_1+k_2-n+1},\ldots, x_{k_2})\, .
  \end{aligned}
  $$
Then $S$ is defined
locally as $D^{k_1+k_2-n}\x K(U_1)\x K(U_2)$ with coordinates $(x_1,\ldots,
x_{k_1+k_2-n},y,z)$ and
  $$
  \begin{aligned}
  f( & x_1,\ldots,
  x_{k_1+k_2-n},y,z)= \\ &= (x_1,\ldots, x_{k_1+k_2-n},C_1(x,z),\ldots,
  C_{n-k_2}(x,z),B_1(x,y),\ldots, B_{n-k_1}(x,y) ) \, .
  \end{aligned}
  $$

\medskip

\begin{theorem}\label{thm:product-dual-1}
Let $(f,S_\mu)$ be an oriented measured $k$-solenoid immersed in $M$
intersecting transversally a closed subvariety $i:N\inc M$ of
codimension $q$, such that $S'=f^{-1}(N)\subset S$ is non-empty.
Consider the oriented measured $(k-q)$-solenoid immersed in $N$,
$(f',S')$, where $f'=f_{|S'}$. Then, under the restriction map
  \begin{equation} \label{eqn:restr}
  i^*: H^{n-k}_c(M) \to H^{(n-q)-(k-q)}_c(N)\, ,
  \end{equation}
the dual of the generalized current $[f,S_\mu]^*$ maps to
$[f',S'_\mu]^*$.
\end{theorem}

\begin{proof}
Let $U\subset M$ be a tubular neighbourhood of $N$ with projection
$\pi:U\to N$. Note that $U$ is diffeomorphic to the unit disc bundle
of the normal bundle of $N$ in $M$. Let $\tau$ be a Thom form for
$N\subset M$, that is a closed form $\tau \in \Omega^q(M)$ supported
in $U$, whose integral in any normal space $\pi^{-1}(n)$, $n\in N$,
is one. The dual of the map (\ref{eqn:restr}) under Poincar\'e
duality is the map
 $$
 H^{k-q}(N) \to H^k(M)\, ,
 $$
which sends $[\b]\in H^{k-q}(N)$ to $[\tilde\b]$, where
$\tilde\b=\pi^*\b \wedge \tau$ (this form is extended from $U$ to
the whole of $M$ by zero). So we only need to see that
 $$
 \la [f,S_\mu], \tilde\b\ra = \la [f',S'_\mu], \b\ra\, .
 $$

Take a covering of $S$ by flow-boxes $U_i\cong D^k \x K(U_i)\cong
D^{q}\x D^{k-q}\x K(U_i)$ so that $U_i'=U_i \cap S'$ is given by
$x_1=\ldots =x_q=0$. Making the tubular neighborhood $U\supset N$
smaller if necessary, we can arrange that $f^{-1}(U)\cap U_i$ is
contained in $D^{q}_r\x D^{k-q}\x K(U_i)$, for some $r<1$. It is
easy to construct a map $\tilde\pi: f^{-1}(U) \to f^{-1}(N)$ which
consists on projecting in the normal directions along the leaves.
Then $f\circ \tilde\pi$ and $\pi\circ f$ are homotopic.

Let $S_i'$ be a measurable partition of $S'$ with $S_i'\subset
U_i'$. We may assume that $S_i=\tilde\pi^{-1}(S_i')$ is contained in
$U_i$. 
The sets $S_i$ form a measurable partition containing $f^{-1}(U)$,
the support of $f^*\tilde{\beta}=f^*(\pi^*\b \wedge \tau)$. Then
  $$
  \begin{aligned}
  \la [f,S_\mu], \tilde\b\ra = & \,
  \sum_i \int_{K(U_i)} \left( \int_{S_i\cap (D^k\x \{y\})} f^*(\pi^*\b \wedge \tau) \right)
   d \mu_{K(U_i)}(y) \\
   = &\, \sum_i \int_{K(U_i)} \left( \int_{S_i\cap (D^k\x \{y\})} \tilde\pi^*f^*\b \wedge f^*\tau) \right)
   d \mu_{K(U_i)} (y)\\
 = &\, \sum_i \int_{K(U_i)} \left( \int_{S_i'\cap (D^{k-q}\x \{y\})} f^*\b \right)
   \left( \int_{f(D^q)} \tau \right) d \mu_{K(U_i)} (y)\\
 = &\, \sum_i \int_{K(U_i')} \left( \int_{S_i'\cap (D^{k-q}\x \{y\})} f^*\b \right)
   d \mu_{K(U_i')} (y)\\
 = &\, \la [f',S'_\mu], \b\ra\, .
 \end{aligned}
 $$
\end{proof}

\begin{theorem} \label{thm:product-dual-2}
 Suppose that $(f_1,S_{1,\mu_1})$, $(f_2,S_{2,\mu_2})$ are two immersed solenoids in $M$
 intersecting transversally, and let $(f,S_\mu)$ be its intersection solenoid. Then
 the duals of the generalized currents satisfy
  $$
  [f,S_\mu]^* = [f_1,S_{1,\mu_1}]^* \cup [f_2,S_{2,\mu_2}]^* \, .
  $$
\end{theorem}

\begin{proof}
  Note that $[f_1,S_1]^* \in H^{n-k_1}_c(M)$ and $[f_2,S_2]^* \in H^{n-k_2}_c(M)$, so
  $[f_1,S_1]^* \cup [f_2,S_2]^*$ and $[f,S]^*$ both live in
  $$
  H^{n-k_1+n-k_2}_c(M) =H^{n-k}_c(M)\, .
  $$

 Consider the immersed solenoid $(F, S_1\x
  S_2)$, where $F=f_1\x f_2:S_1\x S_2\to M\x M$ and $S_1\x S_2$
  has the transversal measure $\mu$ given by
  (\ref{eqn:mu-product}). Let us see that the following equality,
  involving the respective generalized currents,
  $$
  [F,(S_1\x S_2)_\mu]=
  [f_1,S_{1,\mu_1}] \ox [f_2,S_{2,\mu_2}] \in H_{k_1+k_2}(M\x M)
  $$
 holds. We prove this by
 applying both sides to $(k_1+k_2)$-cohomology
classes in $M\x M$. Using the K\"unneth decomposition it is enough
to evaluate on a form $\b=p_1^*\beta_1\wedge p_2^*\b_2$, where
$\b_1,\b_2\in H^*(M)$ are closed forms and $p_1,p_2:M\x M\to M$ are
the two projections. Let $\{U_i\}$, $\{V_j\}$ be open covers of
$S_1$, $S_2$ respectively, by flow-boxes, and let $\{\rho_{1,i}\}$,
$\{\rho_{2,j}\}$ be partitions of unity subortinated to such covers.
Then
 $$
 \begin{aligned}
 \la & [F,(S_1\x S_2)_\mu],\b\ra  = \\ &= \sum_{i,j} \int_{K(U_i)\x K(V_j)}
 \left( \int_{L_{y}\x L_z} (p_1^*\rho_{1,i})\,
 (p_2^*\rho_{2,j}) F^*(p_1^*\beta_1\wedge
 p_2^*\b_2) \right) d\mu_{K(U_i)\x K(V_j)} (y,z) \\ &\ = \sum_{i,j} \int_{K(U_i)\x K(V_j)}
 \left( \int_{L_{y}\x L_z} p_1^*(\rho_{1,i} f_1^*\beta_1) \wedge p_2^*( \rho_{2,j}
 f_2^*\b_2) \right) d\mu_{1,K(U_i)}(y) \, d\mu_{2,K(V_j)}(z) \\ &\ =\left( \sum_{i} \int_{K(U_i)}
 \left( \int_{L_{y}} \rho_{1,i} f_1^*\beta_1 \right) d\mu_{1,K(U_i)} (y) \right)
 \left( \sum_{j} \int_{K(V_j)}
 \left( \int_{L_{z}} \rho_{1,j} f_2^*\beta_2 \right) d\mu_{2,K(V_j)} (y) \right) \\
 &\ =\la [f_1,S_{1,\mu_1}],\b_1\ra \, \la [f_2,S_{2,\mu_2}],\b_2\ra\, ,
 \end{aligned}
 $$
as required.

Now we are ready to prove the statement of the theorem. Let
$\varphi:M\to \Delta$ be the natural diffeomorphism of $M$ with the
diagonal $\Delta \subset M\x M$, and let $i:\Delta \inc M\x M$ be
the inclusion. Then, using theorem \ref{thm:product-dual-2},
  $$
  \begin{aligned} \ 
  [f,S_\mu]^* \, &=  [\varphi \circ f, S_\mu]^* = i^* ([F,(S_1\x
  S_2)_\mu]^*)= \\
   &= i^* ([f_1,S_{1,\mu_1}]^* \ox [f_2,S_{2,\mu_2}]^* )
  = [f_1,S_{1,\mu_1}]^* \cup [f_2,S_{2,\mu_2}]^* \, .
  \end{aligned}
  $$
\end{proof}

\medskip

Let us look more closely to the case where $k_1+k_2=n$.
We assume that $(f_1,S_{1,\mu_1})$
and $(f_2,S_{2,\mu_2})$ are two oriented immersed measured solenoids of
dimensions $k_1,k_2$ respectively,
which intersect transversally. 
Let $(f,S_\mu)$ is the intersection $0$-solenoid of $(f_1,S_{1,\mu_1})$
and $(f_2,S_{2,\mu_2})$.

\begin{definition}\textbf{\em (Intersection index)}
\label{def:intersection-index}
At each point $x=(x_1,x_2)\in S$, 
the intersection index $\epsilon(x_1,x_2)\in \{\pm 1\}$
is the sign of the intersection of the leaf of $S_1$ through
$x_1$ with the leaf of $S_2$ through $x_2$. 
The continuous function $\epsilon:S\to
\{\pm 1\}$ gives the orientation of $S$.
\end{definition}

Recall that the $0$-solenoid $(f,S_\mu)$ comes equipped with a natural measure $\mu$
(for a $0$-solenoid the notions of measure and transversal measure coincide).
If $x=(x_1,x_2)\in S$, then locally around $x$, $S$ is homeomorphic to
$T=T_1\x T_2$, where $T_1$ and $T_2$ are small local transversals of $S_1$
and $S_2$ at $x_1$ and $x_2$, respectively. The measure $\mu_T$ is the
product measure  $\mu_{1,T_1} \x \mu_{2,T_2}$.

\begin{definition}\textbf{\em (Intersection measure)}\label{def:intersection-measure}
The intersection measure is the transversal measure $\mu$ of the
intersection solenoid $(f,S_\mu)$, induced by those of
$(f_1,S_{1,\mu_1})$ and $(f_2,S_{2,\mu_2})$.

\end{definition}

\begin{definition}\textbf{\em (Intersection pairing)} \label{def:intersection-pairing}
We define the intersection pairing as the real number
 $$
 (f_1, S_{1,\mu_1}) \cdot (f_2,S_{2,\mu_2})=
 \int_{S } \epsilon  \ d \mu \, .
 $$
\end{definition}

\begin{theorem} \label{thm:product-RS}
If $(f_1,S_{1,\mu_1})$
and $(f_2,S_{2,\mu_2})$ are two oriented immersed measured solenoids of
dimensions $k_1,k_2$ respectively,
which intersect transversally, such that $k_1+k_2=n$. Then
 $$
 (f_1, S_{1,\mu_1}) \cdot (f_2,S_{2,\mu_2})=
 [f_1, S_{1,\mu_1}]^* \cdot [f_2,S_{2,\mu_2}]^*\, .
 $$
\end{theorem}

\begin{proof}
By theorem \ref{thm:product-dual-2},
 $$
 [f_1, S_{1,\mu_1}]^* \cup [f_2,S_{2,\mu_2}]^* =
 [f,S_\mu]^* \in H^n_c(M,\RR)\, .
 $$
The intersection product $ [f_1, S_{1,\mu_1}]^* \cdot
[f_2,S_{2,\mu_2}]^*$ is obtained by evaluating this cup product on
the element $1\in H^0 (M,\RR)$, i.e.
  $$
  [f_1, S_{1,\mu_1}]^* \cdot [f_2,S_{2,\mu_2}]^* =
  \la [f,S_\mu],1\ra = \int_{S} f^*(1) d\mu(x) =
  \int_{S} \epsilon d\mu \, ,
  $$
since the pull-back of a function gets multiplied by the orientation
of $S$, which is the function $\epsilon$.
\end{proof}

\medskip

When the solenoids are uniquely ergodic we can recover this
intersection index by a natural limitting procedure.

\begin{theorem}
Let $(f_1, S_{1,\mu_1})$ and $(f_2, S_{2,\mu_2})$ be two immersed, oriented,
uniquely ergodic solenoids transversally intersecting. 
Let $l_1\subset S_{1}$ and $l_2\subset
S_{2}$ be two arbitrary leaves. Choose two base points $x_1\in
l_1$ and $x_2 \in l_2$, and fix Riemannian exhaustions $(l_{1,n})$
and $(l_{2,n})$. Define
 $$
 (f_1,l_{1,n}) \cdot (f_2,l_{2,n})=\frac{1}{M_n} \sum_{p=(p_1,p_2)\in
 l_{1,n}\times l_{2,n} \atop f_1(p_1)=f_2(p_2)} \epsilon (p)\, ,
 $$
where $M_n =\Vol_{k_1}(l_{1,n})\cdot \Vol_{k_2}(l_{2,n})$.

Then
 $$
 \lim_{n\to +\infty }
 (f_1,l_{1,n}) \cdot (f_2,l_{2,n})= (f_1,S_{1\mu_1}) \cdot (f_2,S_{2,\mu_2}) \, .
 $$
In particular, the limit exists and is independent of the choices of $l_1$,
$l_2$, $x_1$, $x_2$ and the radius of the riemannian exhaustions.
\end{theorem}

\begin{proof}
The key observation is that because of the unique ergodicity, the
atomic transversal measures associated to the normalizad $k$-volume
of the Riemannian exhaustions (name them $\mu_{1,n}$ and
$\mu_{2,n}$) are converging to $\mu_1$ and $\mu_2$, respectively. 
In particular, in each local flow-box we have
 $$
 \mu_{1,n}\times \mu_{2,n}\to \mu_1\times \mu_2 =\mu \, .
 $$
Therefore the average defining $(f_1,l_{1,n})\cdot (f_2,l_{2,n})$
converges to the integral defining
$(f_1,S_{1\mu_1})\cdot (f_2,S_{2,\mu_2})$ since $\epsilon$ is a
continuous and integrable function (indeed bounded by $1$).
\end{proof}

\begin{remark}
The previous theorem and proof work in the same form for
ergodic solenoids, provided that we know that the Schwartzman limit
measure for almost all leaves is the given ergodic measure. This is
simple to prove for ergodic solenoids with trappings regions mapping to a
contractible ball in $M$
(cf. theorem \ref{thm:11.12bis}).
\end{remark}

\section{Almost everywhere transversality} \label{sec:aet}

The intersection theory developed in section \ref{sec:intersection} is not fully satisfactory
since we do have examples of solenoids (e.g. foliations) which do not intersect
transversally, and cannot even be perturbed to it. However, a weaker notion is enough
to develop intersection theory for solenoids. Indeed, the intersection pairing can 
also be defined for oriented, measured 
solenoids $(f_1, S_{1,\mu_1})$ and $(f_2,S_{2,\mu_2})$ immersed in an oriented
$n$-manifold $M$, with $k_1+k_2=n$, $k_1=\dim S_1$, $k_2=\dim S_2$, which
intersect transversally almost everywhere. Let us first define this notion.

\begin{definition}\textbf{\em (Almost everywhere transversality)}\label{def:almost-transverse}
Let $(f_1, S_{1,\mu_1})$ and $(f_2,S_{2,\mu_2})$ be two measured
immersed oriented solenoids. They intersect almost everywhere
transversally if 
the set 
 $$
 \begin{aligned}
 F=\{(p_1,p_2) & \,\in S_1\x S_2 \, ; \\ &\,  f_1(p_1)=f_2(p_2), 
 df_1(p_1) (T_{p_1}S_1)+df_2(p_2) (T_{p_2}S_2) \neq T_{f_1(p_1)}M \}
 \subset S_1\x S_2
  \end{aligned}
 $$ 
of non-transversal intersection points
is null-transverse in
$S_1\x S_2$ (with the natural product transversal measure $\mu$), i.e. 
if the set of leaves of $S_1\x S_2$ intersecting $F$ has zero
$\mu$-measure.
\end{definition}

We recall that a set $F\subset S_\mu$ in a measured solenoid
is null-transverse if for any local
transversal $T$, the set of leaves passing through $F$ 
intersects $T$ is a set of zero $\mu_T$-measure.

Observe that when $N\subset M$ is a closed submanifold of codimension
$k$ and $(f,S_\mu)$ is a measured immersed oriented $k$-solenoid and
$N\subset M$, then they intersect almost everywhere transversally if
and only if the subspace of non-transversal intersection points 
 $$
 F=\{p\in S \, ;\, f(p) \in N, \ 
 df(p) (T_{p}S)+T_{f(p)}N \neq T_{f(p)}M \}
 \subset S
 $$ 
is null-transverse.

Then it is useful to translate to $S_1\times S_2$ the meaning of
almost everwywhere transversality. We have the following
straightforward lemma.

\begin{lemma}
The solenoids $(f_1, S_{1,\mu_1})$ and $(f_2,S_{2,\mu_2})$ are almost
everywhere transversal if and only if $(f_1\x f_2, (S_1\x S_2)_\mu)$ and the
diagonal $\Delta \subset M\x M$ intersect almost everywhere
transversally.
%
\end{lemma}

Let $(f,S_\mu)$ be an immersed solenoid intersecting transversally
almost everywhere a closed submanifold $N\subset M$. Write
$S'=f^{-1}(N)$ and let $F\subset S'$ be the subset of
non-transversal points. Note that $F$ is closed, hence
$S'_{reg}=S'-F$ is open in $S'$. Moreover, $S'_{reg}$ consists of
the transversal intersections, so the intersection index
$\epsilon:S'_{reg} \to \{\pm 1\}$ is well defined and continuous. We
define the intersection number as
 $$
 \int_{S'-F} \epsilon(x) d\mu(x) \, .
 $$

\begin{theorem} \label{thm:intersect-number-almost}
Suppose that  $(f,S_\mu)$ and $N\subset M$ intersect almost
everywhere transversally. Then
 $$
 [f, S_\mu]^* \cdot [N] = \int_{S'-F} \epsilon \ d\mu \, .
 $$
\end{theorem}

\begin{proof}
Fix an accessory Riemannian metric on $M$. Let $(U_n)$ be a nested
sequence of open neighbourhoods of $F$ in $S$ such that
$\bigcap_{n\geq 1} U_n=F$. Then
 $$
 \int_{S'-U_n} \epsilon \ d\mu \to  \int_{S'-F} \epsilon \ d\mu \, .
 $$

In $S'-U_n$ the angle of intersection between $f(S)$ and $N$ is
bounded below, so there is a small $\rho>0$ (depending on $n$)
such that if $U_\rho$ is
the $\rho$-tubular neighbourhood of $N$ in $M$, then for each
intersection point $x\in S'-U_n$, there is a (topological) disc
$D_x$ contained in a local leaf through $x$, which is exactly the
path component of $f^{-1}(U_\rho)$ through $x$.

Let $\tau_\rho$ be a Thom form for $N\subset M$, that is a closed
$k$-form supported in $U_\rho$, whose integral in the
normal space to $N$ is one. Then $\int_{D_x} \tau_\rho=1$ for any $x\in S'-U_n$.
So
 $$
 \int_{S'-U_n} \epsilon \ d\mu = \int_{A_n}  f^*\tau_\rho\, ,
 $$
where $A_n$ consists of those discs $D_x$ with $x\in S'-U_n$, and
$B_n=f^{-1}(U_\rho) -A_n$, so that $f^{-1}(U_\rho) =A_n\cup B_n\subset S$.
Note that $B_n$ is contained in a neighbourhood of $F$ slightly larger
than $U_n$ (say $U_{n-1}$, taking $\rho$ small enough).

On the other hand,
 $$
 [f, S_\mu]^* \cdot [N] = \la  [f, S_\mu] , [\tau_\rho]\ra  = \int_{S_\mu}
 f^*\tau_\rho = \int_{A_n} f^*\tau_\rho  +\int_{B_n} f^*\tau_\rho =
 \int_{S'-U_n} \epsilon \ d\mu +
 \int_{B_n} f^*\tau_\rho
 \, .
 $$
Note that $\int_{B_n} f^*\tau_\rho$ is independent of $\rho$ small enough.

So we need to see that
 $$
 \int_{B_n}  f^*\tau_\rho \to 0\, ,
 $$
when $n\to \infty$. As $B_n\subset U_{n-1}$ and $U_{n-1}$ has very
small transversal measure, everything reduces to bound uniformly
 $$
 \int_{L_y}  f^*\tau_\rho\, ,
 $$
for any leaf $L_y\subset B_n$. If $L_y$ is contained in a flow-box, then
 \begin{equation} \label{eqn:local-intersect}
 \int_{L_y}  f^*\tau_\rho  = [L_y,\bd L_y]\cdot [N]
 \end{equation}
is the intersection number of $(L_y,\bd L_y)$ with $N$. This is
well-defined because the boundary of $L_y$ does not intersect $N$
(as $\bd L_y\subset \bd B_n$, we have that $\bd L_y \subset \bd
f^{-1}(U_\rho)$, so $\bd L_y \cap N= \emptyset$).
Moreover, by compactness of $S$, the quantity
(\ref{eqn:local-intersect}) must be bounded uniformly. If $L_y$ is
not contained in a flow-box, it may be pieced into several
components (the number of them is uniformly bounded). Then there are
contributions of the boundaries of the components to
(\ref{eqn:local-intersect}), but they cancel each other after
addition.
\end{proof}

Consider now two immersed measured oriented solenoids $(f_1,
S_{1,\mu_1})$, $(f_2,S_{2,\mu_2})$ intersecting almost everywhere
transversally. Let $F\subset S_1\x S_2$ be the subspace of
non-transversal intersection points, which has null-transversal
measure in $S_1\x S_2$. Set $S= (S_1 \x S_2) \cap f^{-1}(\Delta)$.
Then there is an intersection index $\epsilon(x)$ for each $x\in S
-F$ and an intersection measure $\mu$ on $S-F$. We define the
intersection product as
 $$
 \int_{S-F} \epsilon  \ d \mu  \, .
 $$
Theorem \ref{thm:intersect-number-almost} implies that
 $$
 [f_1, S_{1,\mu_1}] \cdot [f_2,S_{2,\mu_2}] = \int_{S-F} \epsilon (x) \ d \mu(x) \, .
 $$

\begin{remark} \label{rem:15.13}
Let $(f,S_\mu)$ be an \emph{embedded} $k$-solenoid in a $n$-manifold
with $n=2k$. Assume that the transversal measure has no atoms. Then
the set of non-transversal points $F \subset S\x S$ is the diagonal
$\Delta_S\subset S\x S$. As the transversal measures of $S$ have no
atoms, $F$ has null-transversal measure, so $S$ intersects itself
almost everywhere transversally. Moreover
  $$
  S'=(S\x S)\cap f^{-1}(\Delta)= \Delta_S\, ,
  $$
because $f$ is injective. Therefore
  $$
  [f, S_{\mu}] \cdot [f,S_{\mu}] = \int_{\Delta_S-F} \epsilon  \ d \mu =0 \, .
  $$
This gives another proof of theorem \ref{thm:self-intersection} in the
case $n=2k$.
\end{remark}

\medskip

Moreover, we can always homotop solenoids so that they become almost
everywhere transverse.

\begin{theorem} \label{thm:homotopic}
  Let $(f,S_\mu)$ be an immersed $k$-solenoid in $M$, and let
  $N\subset M$ be a closed submanifold of codimension $k$.
  Then there exists a homotopic immersion
  $f_1:S\to M$ such that $(f_1,S_\mu)$ intersects $N$
  almost everywhere transversally.
\end{theorem}

\begin{proof}
Let $U=D^k\x K(U)$ be a flow-box for $S$, and $(u_1,\ldots,
u_n)$ coordinates on an open subset $W\subset M$ containing 
$\overline{f(U)}$ such that $N$ is defined as $u_{k+1}=\cdots =u_n=0$ in such
coordinates. Take a smooth function $\rho$ on $S$ with $\rho_{|U}>0$ and
$\rho_{|S-U}\equiv 0$.

In these coordinates, the immersion $f:U =D^k\x K(U) \to W$ is
written as $f(x,y)=(f_1(x,y),\ldots, f_n(x,y))$. Take $\epsilon>0$
small so that $f(x,y)+v \in W$, for any $v\in \RR^n$ with
$|v|<\epsilon$. For each $y\in K(U)$, consider the map
  \begin{equation}\label{eqn:fy}
   \begin{aligned}
    f_y  : D^k & \, \to \RR^k \\
     x & \, \mapsto (f_1(x,y), \ldots, f_k(x,y)) \, .
   \end{aligned}
   \end{equation}
and the map $\bar{f}_y=f_y /\rho$.

We say that $a\in D^k_\epsilon$ is a regular value for $\bar{f}_y$
if $\det(d{\bar{f}}_y(x))\neq 0$ for each $x\in D^k$ with
$\bar{f}_y(x)=a$. Consider the set
  $$
  A:=  \{ (y,a)\in K(U)\x D^k_\epsilon \, ;\, a \text{ is a regular value for }
  \bar{f}_y \}
  $$
is open. The set $A_y=(\{y\} \x D_\epsilon^k)\cap A$ is the set of
regular values of $\bar{f}_y$ on $D^k_\epsilon$. By Sard's theorem,
this is of full measure in $D^k_\epsilon$.

Let $\nu$ be the product measure in $K(U)\x D^k_\epsilon$, i.e.
$\nu=\mu_{K(U)}\x \lambda$, where $\lambda$ is the standard Lebesgue
measure in $D^k_\epsilon$. Then the set $A$ is a set of full
$\nu$-measure. This implies that for $\lambda$-almost every $a\in
D^k_\epsilon$, the set
    $$
    C_a=\{ y \in K(U) \, ; \, a \text{ is not a regular value for } \bar{f}_y
    \} \subset K(U)
    $$
is of zero $\mu_{K(U)}$-measure. Fix one such $a\in D^k_\epsilon$.
Then the map
   $$
   \begin{aligned}
     & \tilde{f} : U=D^k \x K(U) \to W \\
     & \tilde{f}(x,y) = (f_1(x,y)-a_1\rho(x,y), \ldots, 
     f_k(x,y)-a_k\rho(x,y), f_{k+1}(x,y), \ldots, f_n(x,y) )
   \end{aligned}
   $$
is transverse to $N$ for all $y\in K(U)-C_a$: take $x\in D^k$ such that
$\tilde{f}(x,y)\in N$. Then $\bar{f}_y(x)=a$, so $\det(d\bar{f}_y(x))\neq 0$.
But setting  $\tilde{f}_y(x)= (f_1(x,y)-a_1\rho(x,y), \ldots, 
f_k(x,y)-a_k\rho(x,y))$, we have 
 $$
 d\tilde{f}_y(x) = d f_y(x) - a \, d_x\rho = d(\rho \bar{f}_y) (x) - a \, d_x\rho =
 \rho(x,y) \, d\bar{f}_y(x)\,.
 $$
So  $\det(d\tilde{f}_y(x))\neq 0$ (since $\rho(x,y)>0$) and $\tilde{f}$ is transverse
to $N$ at $(x,y)$.

We extend $\tilde{f}$ to the rest of $S$ by setting 
$\tilde{f}_{|S-U}={f}_{|S-U}$. Then $\tilde{f}:S\to M$
is an immersion homotopic to $f$ and close to it, and it satisfies
that it is almost everywhere transversal to $N$ over $U$.

Fixing a finite covering of $S$ by flow-boxes, and repeating this
process, we obtain an immersion homotopic to $f$ and almost
everywhere transversal to $N$. 
\end{proof}

\begin{corollary} \label{cor:homotopic}
  Let $(f_1,S_{1,\mu_1})$ and $(f_2,S_{2,\mu_2})$ be two immersed solenoids in $M$,
  with $k_1+k_2=n$.
  Then we may homotop the immersion $f_1:S_1\to M$ so that they intersect
  almost everywhere transversally.
\end{corollary}

\begin{proof}
Consider the solenoid $S=S_1\x S_2$, and flow-boxes $U_1$ for $S_1$
and $U_2$ for $S_2$. Fix a smooth function $\rho_1$ on $S_1$ 
such that $\rho_{1|U_1}>0$ and $\rho_{1|S_1-U_1}\equiv 0$. We have no need
of introducing a smooth function for $S_2$ since we shall only perturb $S_1$.

Using natural coordinates $(z_1,\ldots, z_n,w_1,\ldots, w_n)$
for $M\x M$, the diagonal is defined by the equations $z_1-w_1=0,
\ldots, z_n-w_n=0$. The immersion $F:S_1\x S_2 \to M\x M$ is given
in coordinates as
  $$
  F(x_1,y_1,x_2,y_2)=(f_1(x_1,y_1), f_2(x_2,y_2))\, ,
  $$
for $(x_1,y_1)\in U_1=D^{k_1}\x K(U_1)$, $(x_2,y_2)\in U_2=D^{k_2}\x K(U_2)$.
So the map (\ref{eqn:fy}) is in this situation
 \begin{equation}\label{eqn:fy-2}
   (x_1,x_2) \mapsto  f_1(x_1,y_1) - f_2(x_2,y_2).
 \end{equation}

As in theorem \ref{thm:homotopic} we find $a\in \RR^n$ small enough
which is a regular value of (\ref{eqn:fy-2}) divided by $\rho_1(x_1,y_1)$, for $\mu_T$-almost
every $(y_1,y_2)\in K(U_1)\x K(U_2)$.

The map $\tilde{f}_1: U_1\to M$ defined as $\tilde{f}_1(x_1,y_1)=
f_1(x_1,y_1) - a\rho_1(x_1,y_1)$ is almost everywhere transversal to $U_2$ along
$U_1$. The rest of the argument is analogous to that at the end of
the proof of theorem \ref{thm:homotopic}.
\end{proof}

\bigskip

Observe that the intersection pairing can be defined for solenoids
endowed with a transversal signed measure (resp.\ complex
measure). The space of such solenoids is a real (resp.\ complex)
vector space, where the sum of immersed solenoids is their disjoint
union, and the zero is the empty set.
With this vector space structure, we have

\begin{proposition}
The intersection pairing is bilinear and graded commutative. 
\end{proposition}

This is a consequence of theorem \ref{thm:product-RS}. Now using
theorem \ref{thm:homeo-RS}, we get

\begin{proposition}
The intersection pairing is invariant by cobordism of immersed solenoids.
\end{proposition}

\begin{definition}
We define the space $S_k(M,\RR)$ of immersed measured solenoids modulo cobordism.
\end{definition}

\begin{theorem}
The intersection pairing of solenoids,
 $$
 S_k(M,\RR)\x S_{n-k}(M,\RR) \to \RR \, ,
 $$
is well-defined, bilinear and graded-commutative. It factors, via the
Ruelle-Sullivan map, through real homology.
\end{theorem}

Whenever two solenoids intersect almost everywhere transversally, the intersection pairing
is defined geometrically. This extends the case of manifolds and homology with
integer coefficients.

\newpage

\setcounter{section}{0}

\renewcommand{\thesection}{{Appendix \Alph{section}}}
\section{Norm on the homology}\renewcommand{\thesection}{\Alph{section}}\label{sec:appendix1-norm}

Let $M$ be a compact $C^\infty$ Riemannian manifold. For each
$a\in H_1(M,\ZZ)$ we define
 $$
 l(a)=\inf_{[\g ]=a} l(\g ) \, ,
 $$
where $\g$ runs over all closed loops in $M$ with homology class
$a$ and $l(\g)$ is the length of $\g$,
 $$
 l(\g) =\int_{\g} \ ds_g \, .
 $$

By application of Ascoli-Arzela it is classical to get

\begin{proposition} \label{prop:A.1}
For each $a\in H_1(M,\ZZ)$ there exists a minimizing geodesic loop
$\g$ with $[\g]=a$ such that
 $$
  l(\g)=l(a) \, .
 $$
\end{proposition}

Note that the minimizing property implies the geodesic character
of the loop. We also have

\begin{proposition} \label{prop:A.2}
There exists a universal constant $C_0=C_0(M)>0$ only depending on
$M$, such that for $a,b \in H_1(M,\ZZ)$ and $n\in \ZZ$, we have
 $$
 l(n\cdot a)\leq |n| \ l(a) \, ,
 $$
and
 $$
 l(a+b)\leq l(a)+l(b)+C_0 \, .
 $$
(We can take for $C_0$ twice the diameter of $M$.)
\end{proposition}

\begin{proof}
Given a loop $\g$, the loop $n\g$ obtained from $\g$ running
through it $n$ times (in the direction compatible the sign of $n$)
satisfies
 $$
 [n\g]=n \ [\g],
 $$
and
 $$
 l(n\g)=|n|\, l(\g)\, .
 $$
Therefore
 $$
 l(n\cdot a)\leq l(n\g)=|n|\, l(\g) \, ,
 $$
and we get the first inequality taking the infimum over $\g$.

Let $C_0$ be twice the diameter of $M$. Any two points of $M$ can
be joined by an arc of length smaller than or equal to $C_0/2$. Given
two loops $\a$ and $\b$ with $[\a ]=a$ and $[\b ]=b$, we can
construct a loop $\g$ with $[\g ]=a+b$ by picking a point in $\a$
and another point in $\b$ and joining them by a minimizing arc
which pastes together $\a$ and $\b$ running through it back and
forth. This new loop satisfies
 $$
 l(\g )=l(\a)+l(\b )+C_0 \, ,
 $$
therefore
 $$
 l(a+b)\leq l(\a)+l(\b )+C_0 \, .
 $$
and the second inequality follows.
\end{proof}

\begin{remark} \label{rem:A.3}
It is not true that $l(n\cdot a)=n\, l(\g )$ if $l(a)=l(\g)$. To
see this take a surface $M$ of genus $g\geq 2$ and two elements
$e_1,e_2\in H_1(M,\ZZ)$ such that
 $$
 l(e_1)+l(e_2)<l(e_1+e_2)\, .
 $$
(For instance we can take $M$ to be the connected sum of a large
sphere with two small $2$-tori at antipodal points, and let $e_1$,
$e_2$ be simple closed curves, non-trivial in homology, inside
each of the two tori.) Let $a=e_1+e_2$. Then
 $$
 l(n \cdot a)=l(n \cdot (e_1+e_2))\leq n\, l(e_1)+n\, l(e_2)+ C_0 \, ,
 $$
we get for $n$ large
 $$
 l(n\cdot a)<n\, l(a) \, .
 $$
\end{remark}

\begin{theorem}\textbf{\em (Norm in homology)} \label{thm:A.4}
Let $a\in H_1(M,\ZZ)$. The limit
 $$
 ||a||=\lim_{n\to +\infty } \frac{l(n\cdot a) }{n} \ \, ,
 $$
exists and is finite. It satisfies the properties
\begin{enumerate}
 \item[(i)] For $a\in H_1(M,\ZZ)$, we have $||a||=0$ if and only if
 $a$ is torsion.
 \item[(ii)] For $a\in H_1(M,\ZZ)$ and $n\in \ZZ$, we have
 $||n\cdot a||=|n|\, ||a|| $ .
 \item[(iii)] For $a,b\in H_1(M,\ZZ)$, we have
 $$
 ||a+b||\leq ||a||+||b|| \, .
 $$
\item[(iv)] $||a||\leq l(a)$.
\end{enumerate}
\end{theorem}

\begin{proof}
Let $u_n=l(n\cdot a)+C_0$. By the properties proved before, the
sequence $(u_n)$ is sub-additive
 $$
 u_{n+m}\leq u_n+u_m \, ,
 $$
therefore
 $$
 \limsup_{n\to +\infty } \frac{u_n}{n}=\liminf_{n\to +\infty }
 \frac{u_n}{n} \, .
 $$
Moreover, we have also
 $$
 \frac{u_n}{n} \leq l(a) <+\infty \, ,
 $$
thus the limit exists and is finite. Property (iv) holds.

Property (ii) follows from
 $$
 ||n\cdot a||= \lim_{m\to \infty} \frac{l(mn\cdot a)}{m} =
 |n|\, \lim_{m\to \infty} \frac{l(m|n|\cdot a)}{m|n|} = |n|\
 ||a||\, .
 $$
Property (iii) follows from
 $$
 l(n\cdot (a+b))\leq l(n\cdot a)+l(n\cdot b) +C_0\leq n\, l(a)+n\, l(b) +C_0 \, ,
 $$
dividing by $n$ and passing to the limit.

Let us check property (i). If $a$ is torsion then $n\cdot a=0$, so
$||a||=\frac1n ||n\cdot a||=0$. If $a$ is not torsion, then there
exists a smooth map $\phi:M\to S^1$ which corresponds to an
element $[\phi]\in H^1(M,\ZZ)$ with $m=\la [\phi],a\ra > 0$. Then
for any loop $\gamma:[0,1]\to M$ representing $n\cdot a$, $n>0$,
we take $\phi\circ\gamma$ and lift it to a map
$\tilde{\g}:[0,1]\to \RR$. Thus
 $$
 \tilde{\g}(1)-\tilde{\g}(0)=\la [\phi],n\cdot a\ra = m\, n\, .
 $$
Now let $C$ be an upper bound for $|d\phi|$. Then
 $$
 m\, n=|\tilde{\g}(1)-\tilde{\g}(0)|=l(\phi\circ \gamma)\leq C\,
 l(\gamma)\, ,
 $$
so $l(\g)\geq m\, n/C$, hence $l(n\cdot a)\geq m\, n/C$ and
$||a||\geq m/C$.
\end{proof}

Now we can define a norm in $H_1(M,\QQ)=\QQ \otimes H_1(M,\ZZ)$ by
 $$
 ||\lambda \otimes a ||=|\lambda|\cdot ||a|| \, ,
 $$
and extend it by continuity to $H_1(M,\RR)=\RR\otimes H_1(M,\ZZ)$.

\renewcommand{\thesection}{{Appendix \Alph{section}}}
\section{De Rham cohomology of solenoids}\renewcommand{\thesection}{\Alph{section}}\label{sec:appendix2-generalities}

In this appendix, we present the definition of the De Rham
cohomology groups for solenoids. The general theory for foliated
spaces from \cite[chapter 3]{MoSch} can be applied to our
solenoids. In \cite{MoSch}, the required regularity is of class $C^{\infty,0}$,
but it is easy to see that the arguments extends to the case of regularity
of class $C^{1,0}$.

Let $S$ be a $k$-solenoid of class $C^{r,s}$ with $r\geq 1$. The space of
$p$-forms $\Omega^p(S)$ consists of $p$-forms on leaves $\a$, such
that $\a$ and $d\a$ are of class $C^{1,0}$. 
Note that the exterior differential
 $$
 d:\Omega^p(X) \to \Omega^{p+1}(X)
 $$
is the differential in the leaf directions. We can define the De
Rham cohomology groups of $S$ as usual,
  $$
  H^p_{DR}(S):= \frac{\ker (d:\Omega^p(S)\to \Omega^{p+1}(S))}{\im
  (d: \Omega^{p-1}(S)\to \Omega^{p}(S))}\, .
  $$
The natural topology of the spaces $\Omega^p(X)$ gives a topology
on $H^p_{DR}(S)$, so this is a topological vector space, which is
in general non-Hausdorff. Quotienting by $\overline{\{0\}}$, the
closure of zero, we get a Hausdorff space
  $$
  \bar H^p_{DR}(S)= \frac{H^p_{DR}(S)}{\overline{\{0\}}} =
  \frac{\ker (d:\Omega^p(S)\to \Omega^{p+1}(S))}{\overline{\im
  (d: \Omega^{p-1}(S)\to \Omega^{p}(S))}}\, .
  $$

We define the solenoid homology as
 $$
 H_k(S):=\Hom_{cont}(H^k_{DR}(S),\RR)= \Hom_{cont}(\bar H^p_{DR}(S),\RR)\, ,
 $$
the continuous homomorphisms from the cohomology to $\RR$. For a
manifold $M$, $H^k_{DR}(M)$ and $H_k(M)$ are equal to the usual
cohomology and homology with real coefficients.

\begin{definition}{\textbf{\em (Fundamental class)}}
\label{def:B.1} Let $S$ be an oriented $k$-solenoid with a
transversal measure $\mu=(\mu_T)$. Then there is a well-defined
map giving by integration of $k$-forms
 $$
 \int_{S_\mu} (\ \cdot \ ) :\Omega^k(S) \to \RR\, ,
 $$
by assigning to any $\alpha\in \Omega^k(S)$ the number
 $$
  \int_{S_\mu}\alpha :=\sum_i \int_{K(U_i)} \left( \int_{L_y\cap S_i}
 \alpha(x,y) dx \right) d\mu_{K(U_i)} (y)\, ,
 $$
where $S_i$ is a finite measurable partition of $S$ subordinated
to a cover $\{U_i\}$ by flow-boxes. It is easy to see, as in section
\ref{sec:Ruelle-Sullivan}, that $\int_{S_\mu} d\beta
=0$ for any $\beta\in \Omega^{k-1}(S)$. Hence $\int_{S_\mu}$ gives
a well-defined map
 $$
 H^k_{DR}(S) \to \RR\, .
 $$
Moreover, this is a continuous linear map, hence it defines an
element
 $$
 [S_\mu]\in H_k(S)\, .
 $$
We shall call $[S_\mu]$ the fundamental class of $S_\mu$.
\end{definition}

The following result is in \cite[theorem 4.27]{MoSch}.

\begin{theorem} \label{thm:Hk}
 Let $S$ be a compact, oriented $k$-solenoid. Then the map
  $$
  \cV_\cT(S) \to H_k(S)\, ,
  $$
 which sends $\mu\mapsto [S_\mu]$, is an isomorphism from the
 space of all signed transversal measures to the $k$-homology of $S$.
\end{theorem}

The set of transversal measures $\cM_\cT(S)$ is a cone, which
generates $\cV_\cT(S)$. Its extremal points are the ergodic
transversal measures. These ergodic measures are linearly
independent. Therefore, the dimension of $H_k(S)$ coincides with
the number of ergodic transversal measures of $S$. Hence, if $S$
is uniquely ergodic, then $H_k(S)\cong \RR$, and $S$ has a unique
fundamental class (up to scalar factor). The uniquely ergodic
solenoids are the natural extension of compact manifolds without
boundary. For a compact, oriented, uniquely ergodic $k$-solenoid
$S$, there is a Poincar{\'e} duality pairing,
  $$
  H^d_{DR}(S)\otimes H^{k-d}_{DR}(S) \to H^k_{DR}(S)
  \stackrel{\int_{S_\mu}}{\too} \RR\, ,
  $$
where $\mu$ is the transversal measure (unique up to scalar).

\medskip

A real vector bundle of rank $m$ over a solenoid $S$ is defined as
follows. A rank $m$ vector bundle over a pre-solenoid $(S,W)$ is a
rank $m$ vector bundle $\pi:E_W\to W$ whose transition functions
 $$
 g_{\a\b}: U_\a\cap U_\b \to \GL(m,\RR)
 $$
are of class $C^{r,s}$. We denote $E=\pi^{-1}(S)$, so that there
is a map $\pi:E\to S$. Let $(S,W_1)$ and $(S,W_2)$ be two
equivalent pre-solenoids, with $f:U_1\to U_2$ a diffeomorphism of
class $C^{r,s}$, $S\subset U_1\subset W_1$, $S\subset U_1\subset
W_1$ and $f_{|S}=\id$, then we say that $\pi_1:E_{W_1}\to W_1$ and
$\pi_2:E_{W_2}\to W_2$ are equivalent if
$\pi_1^{-1}(S)=\pi_2^{-1}(S)=E$ and there exists a vector bundle
isomorphism $\hat{f}:E_{W_1}\to E_{W_2}$ covering $f$ such that
$\hat{f}$ is the identity on $E$. Finally a vector bundle
$\pi:E\to S$ over $S$ is defined as an equivalence class of such
vector bundles $E_W\to W$ by the above equivalence relation.

Note that the total space $E$ of a rank $m$ vector bundle over a
$k$-solenoid $S$ inherits the structure of a $(k+m)$-solenoid
(although non-compact).

A vector bundle $E\to S$ is oriented if each fiber
$E_p=\pi^{-1}(p)$ has an orientation in a continuous manner. This
is equivalent to ask that there exist a representative $E_W\to W$
(where $W$ is a foliated manifold defining the solenoid structure
of $S$) which is an oriented vector bundle over the
$(k+l)$-dimensional manifold $W$.

\medskip

Let $S$ be a solenoid of class $C^{r,s}$ with $r\geq 1$,
and let $E\to S$ be a vector bundle. We may define forms on the total space $E$. A
form $\a\in \Omega^p(E)$ is of vertical compact support if the
restriction to each fiber is of compact support. The space of such
forms is denoted by $\Omega^p_{cv}(E)$. Note that this condition is
preserved under differentials, so it makes sense to talk about the
cohomology with vertical compact supports,
 $$
 H_{cv}^*(E) =\frac{\ker (d:\Omega^p_{cv}(E)\to \Omega^{p+1}_{cv}(E))}{\im
  (d: \Omega^{p-1}_{cv}(E)\to \Omega^{p}_{cv}(E))}\,  .
 $$

\begin{definition} \textbf{\em (Thom form)}
A Thom form for an oriented vector bundle $E\to S$ of rank $m$
over a solenoid $S$ is an $m$-form
 $$
 \Phi\in \Omega_{cv}^m (E)\, ,
 $$
such that $d\Phi=0$ and $\Phi|_{E_p}$ has integral $1$ for each
$p\in S$ (the integral is well-defined since $E$ is
oriented).
\end{definition}

By the results of \cite{MoSch}, Thom forms exist. They represent a
unique class in $H_{cv}^m(E)$, i.e. if $\Phi_1$ and $\Phi_2$ are
two Thom forms, then there is a $\beta\in \Omega_{cv}^{m-1}(E)$
such that $\Phi_2-\Phi_1=d\beta$. Moreover, the map
  $$
  H^{k}(S) \to H^{m+k}_{cv}(E)
  $$
given by
  $$
  [\alpha] \mapsto [\Phi\wedge\pi^*\alpha]\, ,
  $$
is an isomorphism.

\end{document}